\documentclass{amsart}
\usepackage[mathscr]{eucal}
\usepackage{amssymb, amsmath,array, amscd}
\usepackage[OT2,T1]{fontenc}

\newtheorem{thm}{Theorem}
\newtheorem{thma}{Theorem}[section]

\newtheorem{lemma}{Lemma}[section]
\newtheorem{claim}{Claim}[section]

\newtheorem{definition}{Definition}

\newtheorem{rem}{Remark}[section]

\newtheorem{cor}{Corollary}[section]

\newtheorem{Cor}{Corollary}

\newtheorem{prob}{Problem}

\newtheorem{ex}{Example}[section]

\def\wh{\widehat}
\def\wt{\widetilde}

\def\s{\mathcal{S}}

\def\V{\mathcal{V}}
\def\de{\mathcal{D}}

\def\fol{\mathbb{F}\text{ol}}

\def\p{\mathbb{P}}

\def\Te{\Theta}
\def\te{\theta}
\def\lim{\underset{z\to 0}{lim}\,}
\def\pa{\partial}
\def\sup{\supset}
\def\sub{\subset}

\def\I{\mathcal{I}}
\def\d{\delta}

\def\emp{\emptyset}
\def\ov{\overline}
\def\om{\omega}
\def\Om{\Omega}
\def\fa{\mathcal{F}}
\def\a{\alpha}
\def\be{\beta}

\def\ep{\epsilon}
\def\C{\mathbb{C}}

\def\{{\lbrace}
\def\}{\rbrace}

\def\la{\lambda}

\def\U{\mathcal{U}}
\def\H{\mathcal{H}}
\def\L{\mathcal{L}}
\def\g{\gamma}

\def\N{\mathbb{N}}
\def\Si{\Sigma}
\def\te{\theta}
\def\si{\sigma}
\def\G{\mathcal G}

\def\Z{\mathbb{Z}}
\def\O{\mathcal{O}}
\def\*{\star}

\def\Q{\mathbb{Q}}

\def\X{\mathcal{X}}
\def\s{\mathcal{S}}

\begin{document}
\title[Logarithmic foliations]{Logarithmic foliations}



\author[D. Cerveau]{D. Cerveau}
\address{Inst. Math\'ematique de Rennes, Campus de Beaulieu, 35042 RENNES Cedex, Rennes, France}

\email{dominique.cerveau@univ-rennes1.fr}

\author[A. Lins Neto]{A. Lins Neto}
\address{IMPA, Est. D. Castorina, 110, 22460-320, Rio de Janeiro, RJ, Brazil}

\email{alcides@impa.br}

\thanks{}

\begin{abstract}
The purpose of this paper is to study singular holomorphic foliations of arbitrary codimension defined by logarithmic forms on projective spaces.
\end{abstract}



\maketitle

\tableofcontents

\section{Basic definitions and results}\label{ss:1}

Recall that a logarithmic form on a complex manifold $M$ is a meromorphic q-form $\eta$ on $M$ such that the pole divisors of $\eta$ and $d\eta$ are reduced. It is known that a holomorphic form on a compact Kähler manifold is closed. This statement were generalized by Deligne in the context of logarithmic forms as follows:

\begin{thma}\label{t:11}
Let $\eta$ be a logarithmic $q$-form on a compact Kähler manifold $M$. Assume that the pole divisor $(\eta)_\infty$ of $\eta$ is an hypersurface with normal crossing singularities. Then $\eta$ is closed.
\end{thma}

In the case of germs of closed meromorphic 1-forms there are "normal forms" describing them in terms of the poles and residues (cf. \cite{ce-ma}). These normal forms can be translated to the projective spaces and in the logarithmic case they are of the type
\[
\eta=\sum_j\la_j\,\frac{df_j}{f_j}\,\,,\,\,\la_j\in\C^*\,\,,\,\,f_j\,\,\text{holomorphic}\,\,.
\]

One of our purposes is to generalize the above normal form for $p$-forms, $p\ge2$, in a special case. We need a definition. Let $X\sub(\C^n,0)$ be a germ at $0\in\C^n$ of holomorphic hypersurface and $f\in\O_n$ be a reduced germ $f=f_1...f_r$, defining $X$: $X=(f=0)$.

\begin{definition}\label{d:1}
{\rm We say that the hypersurface $X$ has {\it strictly ordinary singularities outside $0$} (briefly s. o. s) if $0\in\C^n$ is an isolated singularity of $f_{i}$ (i.e. $(f_i=0)\setminus\{0\}$ is smooth), $1\le i\le r$, and $X$ is normal crossing outside the origin.}
\end{definition}
The two conditions in definition \ref{d:1} can be expressed as follows: for any sequence of indexes $1\le i_1<...<i_k\le r$ we have
\begin{itemize}
\item{} if $k\ge n$ then $(f_{i_1}=...=f_{i_k}=0)=\{0\}$;
\item{} if $1\le k<n$ then
\[
df_{i_1}(z)\wedge...\wedge df_{i_k}(z)\ne0\,\,,\,\,\forall\,z\in(f_{i_1}=...=f_{i_k}=0)\setminus\{0\}\,\,.
\]
\end{itemize}
In our first result we describe the germs of closed logarithmic p-forms with poles along a hypersurface $X$ with strictly ordinary singularities outside $0$.
With this purpose we introduce a notation that will be used along the paper.
Given $r\in\N$ and $1\le s\le r$ we denote by $\s^r_s$ the set of sequences $I=(i_1,...,i_s)$, where $1\le i_1<...<i_s\le r$.

\begin{thm}\label{t:1}
Let $\eta$ be a germ at $0\in\C^n$ of closed logarithmic $p$-form with poles along a hypersurface $X=(f_1...f_r=0)$ with s.o.s outside $0$. Assume that $n\ge p+2$. Then:
\begin{itemize}
\item[(a).] If $r<p$ then $\eta$ is exact; $\eta=d\Te$, where $\Te$ is logarithmic non-closed and has the same pole divisor as $\eta$.
\item[(b).] If $r\ge p$ then there are numbers $\la_I\in\C$, $I\in \s_p^r$, such that
\end{itemize}
\begin{equation}\label{eq:1}
\eta=\underset{I=(i_1<...<i_p)}{\sum_{I\in\,\s_r^p}}\la_I\,\frac{df_{i_1}}{f_{i_1}}\wedge...\wedge\frac{df_{i_p}}{f_{i_p}}\,+\,d\Te\,\,,
\end{equation}
where, either $\Te=0$, or $\Te$ is logarithmic non-closed and has pole divisor contained in $X$.
\end{thm}

\begin{rem}\label{r:11}
{\rm In the above statement, if $r=0$ then $X=\emp$ and $\eta$ is holomorphic and closed. In this case it can be written as $\eta=d\Te$, where $\Te$ is a holomorphic $(p-1)$-form, by Poincaré lemma.
On the other hand, if $p=1$ and $r\ge1$ then $\eta$ can be written as
\[
\eta=\sum_j\la_j\,\frac{df_j}{f_j}+dg\,\,,\,\,g\in\O_n\,\,,
\]
whereas when $p=2$ and $r\ge2$ then it can be proved that
\[
\eta=\sum_{i<j}\la_{ij}\,\frac{df_i}{f_i}\wedge\frac{df_j}{f_j}+\sum_jdg_j\wedge\frac{df_j}{f_j}+d\a\,\,,
\]
where $g_1,...,g_r\in\O_n$ and $\a\in\Om^1(\C^n,0)$.}
\end{rem}

\begin{rem}\label{r:12}
{\rm The numbers $\la_I$ in (\ref{eq:1}), $I\in\s_p^r$, are called the numerical residues of $\eta$. Given $I=(i_1<...<i_p)$ then $\la_I$ can be calculted by integrating $\eta$ as follows: since $1\le p< n$ the germ of analytic set $X_I:=(f_{i_1}=...=f_{i_p}=0)$ has dimension $n-p\ge1$. 
Moreover, by the normal crossing condition the set $\wt{X}_I:=X_I\setminus\bigcup_{j\notin I}(f_j=0)$ is not empty. If we fix $m\in \wt{X}_I$ then there are local coordinates $z=(z_1,...,z_n)$ such that $z(m)=0$ and $f_{i_j}=z_j$ for all $j=1,...,p$. Given $\ep>0$ small, consider the real p-dimensional torus
\[
T_\ep^p=\left\{z\,;\,|z_j|=\ep\,\text{if $1\le j\le p$, and $z_j=0$ if $j>p$}\right\}\,.
\]
It follows from (\ref{eq:1}) that
\[
\la_I=\frac{1}{(2\pi i)^p}\int_{T_\ep^p}\eta\,\,.
\]}
\end{rem}

As a consequence of theorem \ref{t:1} and \ref{r:11} we can state the following in the case of logarithmic $p$-forms on $\p^n$:

\begin{Cor}\label{c:1}
Let $\eta$ be a logarithmic $p$-form on $\p^n$, $p\le n-1$. Assume that the divisor of poles $(\eta)_\infty$ is given in homogeneous coordinates by $f_1...f_r$, where the $f_{i's}$ are irreducible homogeneous polynomials. Furthermore suppose that the hypersurface $X=(f_1...f_r=0)$ has s.o.s outside $0\in\C^{n+1}$. Then $r\ge p+1$ and
there are numbers $\la_I$, $I\in\s_r^p$, such that in homogeneous coordinates we have
\begin{equation}\label{eq:2}
\eta=\underset{I=(i_1<...<i_p)}{\sum_{I\in \s_r^p}}\la_I\,\frac{df_{i_1}}{f_{i_1}}\wedge...\wedge\frac{df_{i_p}}{f_{i_p}}\,\,,
\end{equation}
where $i_R\eta=0$.
\end{Cor}
In the above statement $i_R$ means the interior product.
\vskip.1in
\begin{ex}
{\rm Theorem \ref{t:1} is false if $p=n-1$ as shows the following example in $\C^n$:
\begin{equation}\label{eq:1'}
\eta=\frac{i_R\,(dz_1\wedge...\wedge dz_n)}{P(z_1,...,z_n)}=\frac{\sum_{j=1}^n(-1)^{j-1}\,z_j\,\,dz_1\wedge...\wedge\wh{dz_j}\wedge...\wedge dz_n}{P(z_1,...,z_n)}\,\,,
\end{equation}
where in (\ref{eq:1'}) $P$ is an irreducible homogeneous polynomial of degree n and $\wh{dz_j}$ means omission of $dz_j$ in the product.

We would like to observe that the same example shows that Corollary \ref{c:1} is false in $\p^m=\p^{n-1}$ if $p=m$: the form $\eta$ represents in homogeneous coordinates a closed logarithmic p-form on $\p^m$ which is not like in the statement of the corollary.} 
\end{ex}

{\bf Notation.} Let us fix homogeneous polynomials $f_1,...,f_r\in\C[z_0,...,z_n]$. The projectivization of the vector space of p-forms $\eta$ that can be written as in (\ref{eq:2}) (not satisfying $i_R\eta=0$ necessarily) will be denoted by $\L^p(f_1,...,f_r)$.
The subspace of forms $\eta\in\L^p(f_1,...,f_r)$ such that $i_R\eta=0$ will be denoted by $\L^p_R(f_1,...,f_r)$. Note that $\L_R^p(f_1,...,f_r)\subsetneq\L^p(f_1,...,f_r)$.

\vskip.1in 

We now turn our atention to $p$-forms defining codimension $p$ foliations. A holomorphic $p$-form $\om$, on an open subset $U\sub\C^n$, defines a codimension $p$ distribution, outside its singular set $Sing(\om)=\{z\in U\,|\,\om(z)=0\}$, if it is {\it locally totally decomposable} (briefly l.t.d) on $U\setminus Sing(\om)$. This means that for any $z\in U\setminus Sing(\eta)$ there are holomorphic 1-forms $\om_1,...,\om_p$, defined in some neighborhood $V$ of $z$, such that $\om|_V=\om_1\wedge...\wedge\om_p$. The distribution $\de$ is then defined on $U\setminus Sing(\om)$ by the codimension $p$ planes
\[
\de_z=ker(\om(z)):=\{v\in T_zU\,|\,i_v\,\om(z)=0\}=\bigcap_{1\le j\le p}\,ker(\om_j(z))\,.
\]

\begin{definition}\label{d:2}
{\rm A holomorphic p-form $\om$ will be said integrable if it is l.t.d outside its singular set and satisfies Frobenius integrability condition. In this context it means that, if $\om|_V=\om_1\wedge...\wedge\om_p$ as above then $d\om_j\wedge\om=0$ for all $j=1,...,p$.}
\end{definition}

Remark that if $\om$ is closed and l.t.d then the Frobenius condition is automatic:
\[
\om_j\wedge\om=0\,\,,\,\,\forall\,j\,\,\implies\,\,d\om_j\wedge\om=d(\om_j\wedge\om)=0\,\,,\,\,\forall\,j\,\,.
\]
In particular, if $\om$ is a closed logarithmic $p$-form then it is integrable if and only if it is l.t.d outside $(\om)_\infty\cup Sing(\om)$.

\begin{ex}
{\rm Let $f_1,...,f_r$ be irreducible homogeneous polynomials on $\C^{n+1}$. Then any 1-form $\te\in\L^1_R(f_1,...,f_r)$ defines a logarithmic codimension one foliation on $\p^n$, denoted by $\fa_\te$. Let $\te_1,...,\te_p\in\L^1_R(f_1,...,f_r)$ and $\eta:=\te_1\wedge...\wedge\te_p$. If $\eta\not\equiv0$ then $\eta\in\L_R^p(f_1,...,f_r)$ and defines a singular codimension p foliation on $\p^n$, denoted by $\fa_\eta$. The leaves of $\fa_\eta$, outside the pole divisor $f_1...f_r=0$, are contained in the intersection of the leaves of $\fa_{\te_1}$,..., $\fa_{\te_p}$. By this reason, $\fa_\eta$ is called the intersection of the foliations $\fa_{\te_1}$,..., $\fa_{\te_p}$.}
\end{ex}

{\bf Notation.} We will use the notation $\L^p_\fa(f_1,...,f_r)=\{\eta\in\L^p_R(f_1,...,f_r)\,|\,\eta$ is integrable$\}$.

\begin{rem}
{\rm We would like to observe that $\L^p_\fa(f_1,...,f_r)$ is an algebraic subset of $\L^p_R(f_1,...,f_r)$.
The proof is left as an exercise to the reader.}
\end{rem}

\begin{definition}\label{d:3}
{\rm We say that $\eta\in\L^p(f_1,...,f_r)$ is {\it totally decomposable into logarithmic forms} (biefly t.d.l.f) if $\eta=\te_1\wedge...\wedge\te_p$, where $\te_1,...,\te_p\in\L^1(f_1,...,f_r)$.
We will use the notation
\[
\L^p_{td}(f_1,...,f_r)=\{\eta\in\L^p_R(f_1,...,f_r)\,|\,\eta\,\,\text{is t.d.l.f}\}\,\,.
\]}
\end{definition}

Observe that $\L_{td}^p(f_1,...,f_r)$ is an irreducible algebraic subset of $\L^p_\fa(f_1,...,f_r)$.
\begin{prob}\label{p:1}
When $\L^p_{td}(f_1,...,f_r)=\L^p_\fa(f_1,...,f_r)$ ?
\end{prob}

A partial answer for problem \ref{p:1} is the following:
\begin{thm}\label{t:2}
Let $f_1,...,f_r$ be homogeneous polynomials on $\C^{n+1}$ and assume that the pole divisor $f_1...f_r=0$ has s.o.s outside $0\in\C^{n+1}$.
Then:
\begin{itemize}
\item[(a).] If $p=2$, or $r\in\{p+1,p+2\}$ then $\L^p_{td}(f_1,...,f_r)=\L^p_\fa(f_1,...,f_r)$.
\item[(b).] If $2<p\le n-2$ and $r>p+2$ then $\L^p_{td}(f_1,...,f_r)$ is an irreducible component of $\L^p_\fa(f_1,...,f_r)$. In particular, if $\L^p_\fa(f_1,...,f_r)$ is irreducible then $\L^p_{td}(f_1,...,f_r)=\L^p_\fa(f_1,...,f_r)$.
\end{itemize}
\end{thm}

An interesting consequence of theorem \ref{t:2} is the following:
\begin{Cor}\label{c:2}
In the hypothesis of theorem \ref{t:2} if $r=p+1$ and $\eta\in\L^p_\fa(f_1,...,f_{p+1})$ then the foliation $\fa_\eta$ in $\p^n$ is a rational fibration of codimension $p$ on $\p^n$. In other words, $\fa_\eta$ has a rational first integral $F\colon\p^n-\to\p^p$ that in homogeneous coordinates can be written as
\[
F=\left(f_1^{k_1},...,f_{p+1}^{k_{p+1}}\right)\,\,,
\]
where $k_1.\,deg(f_1)=...=k_{p+1}.\,deg(f_{p+1})$ and $gcd(k_1,...,k_{p+1})=1$.
\end{Cor}

Let us recall the definition of the {\it Kupka set} of a foliation. Let $\om$ be an 
integrable p-form defining the foliation $\fa$ on an open subset $U\sub\C^n$. We say that $p\in U$ is a singularity of Kupka type of $\om$ if $\om(p)=0$ and $d\om(p)\ne0$. The Kupka set of $\om$, denoted by
$K(\om)$, is the set of singularities of Kupka type of $\om$. Let us note the following facts:
\begin{itemize}
\item[1.] If $f\in\O^*(U)$ then $K(\om)=K(f.\,\om)$. In particular, the Kupka set depends only of the foliation $\fa$ and not on the p-form defining it. The Kupka set of the foliation $\fa$ will be denoted $K(\fa)$.
\item[2.] If $dim(\fa)=n-p>1$ then $\fa$ has the {\it local product property} near any point $p\in K(\fa)$: the germ of $\fa$ at $p$ is holomorphically equivalent to the product of a singular one dimensional foliation at $0\in\C^{p+1}$ by a regular foliation of dimension $n-p-1$ (cf. \cite{ku} and \cite{med}).
This means that there are local coordinates around $p$, $z=(x,y)\in\C^{p+1}\times\C^{n-p-1}$ and a germ of vector field $Y=\sum_{j=1}^{p+1}Y_j(x)\frac{\pa}{\pa x_j}$ such that the germ of $\fa$ at $p$ is equivalent to the foliation defined the germ of p-form at $0\in\C^{p+1}\times\C^{n-p-1}$
\[
\om=i_Y\,dx_1\wedge...\wedge dx_{p+1}\,.
\] 
\end{itemize}

If $Sing(\fa)$ has an irreducible component $Z$ entirely contained in $K(\fa)$ we say that $Z$ is a {\it Kupka component} of $\fa$.

\begin{ex}
{\rm In the case $r=p+1$ (corollary \ref{c:2}) the set $(f_1=...=f_{p+1}=0)\sub Sing(\fa_\eta)$ is a  Kupka component of $\fa_\eta$ (cf. \cite{celn}).}
\end{ex}

\begin{ex}
{\rm In the case $p=n-1\ge2$ then $\L^{\,n-1}_\fa(f_1,...,f_r)=\L^{\,n-1}_R(f_1,...,f_r)$, because any $(n-1)$-form in $\C^n$ is integrable. Moreover, if $r\ge p+2$ then $\L^{\,n-1}_{td}(f_1,...,f_r)$ is a proper algebraic subset of $\L^{\,n-1}_\fa(f_1,...,f_r)$. The reason is that if $\eta$ is decomposable, $\eta=\te_1\wedge...\wedge\te_p$, where $\te_1,...,\te_p$ are logarithmic 1-forms as in theorem \ref{t:2}, then $\eta$ cannot have isolated singularities outside its pole divisor. A specific example on $\p^3$ is given in homogeneous coordinates by the logarithmic 2-form
\[
\eta=\sum_{1\le i<j\le6}
\la_{ij}\,\frac{d\ell_i}{\ell_i}\wedge\frac{d\ell_j}{\ell_j}\,\,,
\]
where $\la_{ij}\in\C$, $1\le i<j\le 6$, and $\ell_j\in\C[z_0,...,z_3]$ is homogeneous of degree one, $1\le j\le6$.
If we choose the $\ell_{j\,'s}$ and $\la_{ij\,'s}$ generic then the foliation $\fa_\eta$ defined by $\eta$ has degree three and $40=3^3+3^2+3+1$ isolated singularities. Each plane $\ell_j$ is $\fa_\eta$-invariant and the reatriction $\fa_\eta|_{\ell_j}$ also defines a degree three foliation and so contains $13=3^2+3+1$ singularities, $1\le j\le6$, each line $\ell_i\cap\ell_j$ contains $4=3+1$ singularities, $1\le i<j\le6$, and each point $\ell_i\cap\ell_j\cap\ell_k$ one singularity, $1\le i<j<k\le6$.
In particular, there are $13\times 6-4\times\#(\ell_i\cap\ell_j)+\#(\ell_i\cap\ell_j\cap\ell_k)=38$ singularities contained in $\bigcup_{j=1}^6\ell_j$ and so $2=40-38$ singularities not contained in the pole divisor. If $\eta$ was decomposable as in theorem \ref{t:2} then these two singularities could not be isolated.

As a consequence of theorem \ref{t:2} we can assert that if $\G$ is a codimension two foliation on $\p^4\sup\p^3$ such that $\G|_{\p^3}=\fa_\eta$ then $\G$ cannot be tangent to $\p^3$ outside the pole divisor $\bigcup_j\ell_j$. In fact, we have the following result:}
\end{ex}

\begin{thm}\label{t:3}
Let $\G$ be a codimension $p$ foliation on $\p^n$, where $n\ge p+2$. Assume that there is a $p+1$ plane $\Si\simeq\p^{p+1}$ such that the foliation $\fa:=\G|_\Si$ has singular set $Sing(\fa)$ of codimension $\ge3$. Then $\G$ is the pull-back of $\fa$ by some linear projection $T\colon\p^n-\to\Si$. In particular, there exists an affine coordinate system $(z,w)\in\C^{p+1}\times\C^{n-p-1}=\C^n\sub\p^n$ such that
$\G$ is represented in these coordinates by a $p$-form depending only of $z$ and $dz$:
\[
\eta=\sum_{j=1}^{p+1}P_j(z)\,dz_1\wedge...\wedge \wh{dz_j}\wedge...\wedge dz_{p+1}=i_X\,dz_1\wedge...\wedge dz_n\,\,,
\]
where $X=\sum_{j=1}^{p+1}(-1)^{j-1}P_j(z)\frac{\pa}{\pa z_j}$.
\end{thm}

In the proof we will use the local version of theorem \ref{t:3}.
We will consider the following situation: let $Z_o\not\equiv0$ be a germ at $(\C^{p+1},0)$ of holomorphic vector field, where $p+1\ge3$. The germ of foliation defined by $Z_o$ is also defined by the germ of $p$-form $\eta_o=i_{Z_o}\nu$, where $\nu=dz_1\wedge...\wedge dz_{p+1}$.
If $Z_o=\sum_{j=1}^{p+1}f_j(z)\frac{\pa}{\pa z_j}$ then
\[
\eta_o=\sum_{j=1}^{p+1}(-1)^{j-1}f_j(z)\,dz_1\wedge...\wedge\wh{dz_j}\wedge...\wedge dz_{p+1}\,\,.
\]
We will assume that there is a germ of integrable holomorphic $p$-form $\eta$ at $0\in\C^n$, where $\C^n=\C^{p+1}\times\C^{n-p-1}$, $n>p+1$, such that $\eta_o=i^*\,\eta$, where $i$ is the inclusion $\C^{p+1}\mapsto\C^{p+1}\times\C^{n-p-1}$.

\begin{thm}\label{t:4}
In the above situation, assume that $cod\left(Sing(Z_o)\right)\ge3$.
Then there exists a local coordinate system $(z,w)\in(\C^{p+1}\times\C^{n-p-1},(0,0))$ and an unity $\phi\in\O^*_n$ such that
\[
\eta=\phi\,\sum_{j=1}^{p+1}(-1)^{j-1}\,f_j(z)\,dz_1\wedge...\wedge\wh{dz_j}\wedge...\wedge dz_{p+1}=\phi\,i_{Z_o}dz_1\wedge...\wedge dz_{p+1}\,.
\]
In particular, the foliation generated by $\eta$ is equivalent to the product of the singular one dimensional foliation generated by $Z_o$ by the non-singular foliation of dimension $n-p-1$ with leaves $z=constant$.
\end{thm}

Another kind of result that we will prove concerns the "stability" of logarithmic foliations on $\p^n$, $n\ge3$. 
In order to precise this phrase we recall the definition of the degree of a foliation on $\p^n$. 

\begin{definition}\label{d:4}
{\rm Let $\fa$ be a holomorphic foliation of codimension $p$ on $\p^n$, $1\le p<n$. The degree of $\fa$, $deg(\fa)$, is defined as the degree of the divisor of tangencies of $\fa$ with a generic plane of complex dimension $p$ of $\p^n$.}
\end{definition}

\begin{rem}\label{r:13}
{\rm In the particular case of codimension one foliations the degree is the number of tangencies of the foliation with a generic line $\p^1\sub\p^n$. More generally, a codimension $p$ foliation $\fa$ on $\p^n$ can be defined by a meromorphic integrable $p$-form on $\p^n$, say $\eta$, with  $Cod_\C(Sing(\eta))\ge2$. If we consider a generic $p$-plane $\Si\simeq \p^p\sub\p^n$ then the degree of $\fa$ is the degree of the divisor of zeroes of $\eta|_\Si$.

Note that, if $\Pi\colon\C^{n+1}\setminus\{0\}\to\p^n$ is the canonical projection, then the foliation $\Pi^*(\fa)$ can be extended to a foliation $\fa^*$ on $\C^{n+1}$. This foliation is represented by a holomorphic $p$-form $\eta$ whose coefficients are homogeneous polynomials of degree $deg(\fa)+1$ and such that $i_R\eta=0$, where $R$ is the radial vector field on $\C^{n+1}$. We say that the form $\eta$ represents $\fa$ in homogeneous coordinates.

A consequence of the definition, is that if $T\colon\p^m-\to\p^n$ is a linear map of maximal rank, where $m>p$, then $deg(T^*(\fa))=deg(\fa)$. In particular, if $\p^m\sub\p^n$ is a generic $m$-plane, where $m>p$, then the degree of $\fa|_{\p^m}$ is equal to the degree of $\fa$.}
\end{rem}

The space of dimension $k$ (codimension $p=n-k$) foliations on $\p^n$ of degree $d$ will be denoted by $\fol(d;k,n)$. 
Note that $\fol(d;k,n)$ can be identified with the subset of the projectivisation of the space of $(n-k)$-forms $\eta$ on $\C^{n+1}$ such that: $\eta$ is integrable, $\eta$ has homogeneous coefficients of degree $d+1$, $cod_\C(Sing(\eta))\ge2$ and $i_R\eta=0$.

When $k=1$ the integrability condition is automatic and $\fol(d;1,n)$ is a Zariski open and dense subset of some projective space $\p^N$.
However, if $k\ge2$ then the integrability condition is non-trivial and $\fol(d;k,n)$ is an algebraic subset of some Zariski open and dense subset of a projective space.

\begin{ex}\label{ex:14}
{\rm Let $\fa$ be the logarithmic foliation on $\p^n$ defined in homogeneous coordinates by an integrable $p$-form $\eta$ on $\C^{n+1}$ as below:
\begin{equation}\label{eq:3}
\eta=\underset{I=(i_1<...<i_p)}{\sum_{I\in \s_r^p}}\la_I\,\frac{df_{i_1}}{f_{i_1}}\wedge...\wedge\frac{df_{i_p}}{f_{i_p}}\,\,,
\end{equation}
where $f_1,...,f_r$ are homogeneous polynomials on $\C^{n+1}$ with $deg(f_j)=d_j$, $1\le j\le r$.
We assume also that $f_1,...,f_r$ are normal crossing outside the origin and $\la_I\ne0$, $\forall I\in\s_r^p$. With these conditions then the holomorphic form $\wt\eta:=f_1...f_r\,\eta$ has singular set of codimension $\ge2$ and so defines $\fa$ in homogeneous coordinates.
Since the degree of the coefficients of $\wt\eta$ is $\sum_{j=1}^rd_j-p$ we obtain
\[
deg(\fa)=\sum_{j=1}^rd_j-p-1:=D(d_1,...,d_r,p)\,\,\implies\,\,\fa\in\fol(D(d_1,...,d_r,p);n-p,n)
\]
{\bf Notation.} The space of dimension $k=n-p$ logarithmic foliations of $\p^n$ defined by a closed $p$-form as in (\ref{eq:3}) will be denoted by
$\L_\fa(d_1,...,d_r;k,n)$. Note that $\L_\fa(d_1,...,d_r;k,n)\sub\fol(D(d_1,...,d_r,p);k,n)$.

The sub-space of $\L_\fa(d_1,...,d_r;k,n)$ of foliations that are defined by t.d.l.f p-forms will be denoted by $\L_{td}(d_1,...,d_r;k,n)$.}
\end{ex}
The next result generalizes a theorem by Calvo-Andrade \cite{CA}\,:
\begin{thm}\label{t:5}
If $k\ge2$ and $r\ge p+2=n-k+2$ then $\ov{\L_{td}(d_1,...,d_r;k,n)}$ is an irreducible component of $\fol\,(D(d_1,...,d_r,p);k,n)$ for all $r>p$ and $d_1,...,d_r\ge1$.
\end{thm}

\begin{rem}
{\rm The above result is also true in the case $r=p+1$. In fact, in \cite{CJI} it is proven the stability of foliations induced by rational maps. On the other hand, by corollary \ref{c:2} the set $\L_\fa(d_1,...,d_{p+1};n-p,n)$ coincides with the set of foliations induced by a rational map 
\[
F=(f_1^{k_1},...,f_{p+1}^{k_{p+1}})\colon\C^{n+1}\to\C^{p+1}\,\,,
\]
where $deg(f_j)=d_j$ and $k_1.\,d_1=...=k_{p+1}.\,d_{p+1}$.}
\end{rem}

The proof of theorem \ref{t:5} will be done first in the case of foliations of dimension two. The general case will be reduced to this one by using the following result:

\begin{thm}\label{t:6}
Let $\fa$ be a codimension $p$ holomorphic foliation on $\p^n$, $n>p+1$. Assume that there is an algebraic smooth submanifold $M\sub\p^n$, $dim_\C(M)=m$, where $p+1\le m<n$, such that:
\begin{itemize}
\item{} The set of tangencies of $\fa$ with $M$ has codimension $\ge2$ on $M$.
\item{} $\fa|_M$ can be defined by a closed meromorphic $p$-form on $M$, say $\eta$.
\end{itemize}
Then $\eta$ can be extended to a closed meromorphic $p$-form $\wt\eta$ on $\p^n$ defining $\fa$.
Moreover, if $\eta$ is logarithmic so is $\wt\eta$.
\end{thm}
In fact, theorem \ref{t:5} is a generalization of a result in \cite{CLS} (see also \cite{ln2}).
\vskip.1in
Theorem \ref{t:5} and problem \ref{p:1} motivate the following question:

\begin{prob}
When $\L_\fa(d_1,...,d_r;k,n)=\L_{td}(d_1,...,d_r;k,n)$?
\end{prob}

\begin{rem}
{\rm Just before finishing this paper we have found a work by Javier Gargiulo Acea \cite{jga} in which he studies some of the problems that we have treated in our paper.
For instance, he obtains the same results (decomposability and stability) of our theorems \ref{t:2} and \ref{t:5} in the case $p=2$ (2-forms). He also proves the normal form for logarithmic p-forms on $\p^n$ if the pole divisor is normal crossing and $p\le n-1$ (our corollary \ref{c:1}). The local case and the logarithmic foliations of codimension $\ge3$ are not treated by him.
We would like to observe that his proof of the stability of logarithmic 2-forms is purely algebraic: he computes the Zariski tangent space at a generic point.} 
\end{rem}

\section{Normal forms}\label{ss:2}
The aim of this section is to prove theorem \ref{t:1} and its corollary.
\subsection{Preliminaries.}
Let $\eta$ be a germ at $0\in\C^n$ of meromorphic $p$-form with reduced pole divisor $X=(f_1...f_r=0)$, $r\ge1$. At the begining we will not assume that $\eta$ is closed.
\begin{rem}\label{r:21}
{\rm It follows from the definition that a germ at $0\in\C^n$ of meromorphic $p$-form $\eta$ is logarithmic if its pole divisor is reduced, $(\eta)_\infty=(f_1...f_r)$, and $f_1...f_r.\,d\eta$ is holomorphic. Since $(\eta)_\infty=(f_1...f_r)$ we can write $\eta=\frac{1}{f_1...f_r}\,\,\om$ where $\om\in\Om^p_n$ is a germ of holomorphic $p$-form. We would like to observe that the following assertions are equivalent:
\begin{itemize}
\item[(a).] $\eta=\frac{1}{f_1...f_r}\,\,\om$ is logarithmic.
\item[(b).] $f_j$ divides $df_j\wedge\om$, for all $1\le j\le r$.
\end{itemize}
In particular, we have:
\begin{itemize}
\item[(c).] $\frac{1}{f_1...f_s}\,\,\om$ is logarithmic, for all $s\le r$.
\end{itemize}

In fact:
\[
\text{(a) $\iff$ $f_1...f_r.\,d\eta=d\om-\frac{d(f_1...f_r)}{f_1...f_r}\wedge\om=\mu$ is holomorphic $\iff$}
\]
\[
f_1...f_r.\,\sum_j\frac{df_j}{f_j}\wedge\om=f_1...f_r\,(d\om-\mu)\,\,\iff\,\,f_j\,\,\text{divides $df_j\wedge\om$, $1\le j\le r$.}
\]}
\end{rem}

The proof of theorem \ref{t:1} will be based in the following:

\begin{lemma}\label{l:21}
Let $\eta=\frac{1}{f_1...f_r}\,\,\om$ be a germ of at $0\in\C^n$ of logarithmic $p$-form, where $1\le p\le n-2$.
Assume that the pole divisor of $\eta$ is $X=(f_1...f_r=0)$, $r\ge1$, has s.o.s outside $0$.
Then $\eta$ can be written as
\begin{equation}\label{eq:4}
\eta=\a_0+\sum_{s=1}^{p-1}\left({\sum_{I\in\s_r^s}}\a_I\wedge\frac{df_{i_1}}{f_{i_1}}\wedge...\wedge\frac{df_{i_s}}{f_{i_s}}\right)+{\sum_{J\in\s_r^p}}g_J.\,\frac{df_{j_1}}{f_{j_1}}\wedge...\wedge\frac{df_{j_p}}{f_{j_p}}\,\,,
\end{equation}
where $\a_0\in\Om_n^p$, $\a_I\in\Om^{p-s}_n$ if $I\in\s_r^s$, $s<p$, and $g_J\in\O_n$ if $J\in\s_r^p$.
\end{lemma}

{\it Proof.}
In the proof we will use the well known concept of residue of a logarithmic form along an irreducible pole (cf. \cite{d}).
Let $\eta=\frac{1}{f_1...f_r}\,\om$ be a germ at $0\in\C^n$ of logarithmic $p$-form as above.
Let us define its residue along $Y_k:=(f_k=0)$, $1\le k\le r$. Fix representatives of $f_1,...,f_r$ and $\eta$, denoted by the same symbols, on some polydisc $Q$. We will assume that the $f_{j's}$ are irreducible in $Q$, and that the divisor $f_1...f_r$  has s.o.s on $Q\setminus\{0\}$. In particular, the $f_{j's}$ have isolated singularity at $0\in Q$.
We have seen that $f_k$ divides $df_k\wedge\om$. In particular, we can write $df_k\wedge\om=f_k.\,\te$ where $\te\in\Om_n^{p+1}$. 
This implies that $df_k\wedge\te=0$. Since $df_k$ has an isolated singularity at $0\in Q$ and $p+1\le n-1$, it follows from de Rham's division theorem \cite{dr} that $\te=df_k\wedge\be_k$, where $\be_k\in\Om^p(Q)$. Therefore, we can write $df_k\wedge(\om-f_k.\,\be_k)=0$ which implies, via the division theorem \cite{dr}, that
there exists $\a_k\in\Om^{p-1}(Q)$ such that $\om=\a_k\wedge df_k+f_k.\,\be_k$.
The residue of $\frac{1}{f_k}\,\om$ along $Y_k$ is the $(p-1)$-form along $Y_k$ defined as
$Res\left(\frac{1}{f_k}\,\om,Y_k\right):=\a_k|_{Y_k}$.
Finaly, the residue of $\eta=\frac{1}{f_1...f_r}\,\om$ along $Y_k$ is defined as
$Res(\eta,Y_k):=\frac{1}{f_1...\wh{f_k}...f_r}\,\a_k|_{Y_k}$, where $\wh{f_k}$ means omission of the factor $f_k$ in the product.

\begin{rem}\label{r:22}
{\rm Let $\eta$ and $Y_k$ be as above. It is well known that $Res\left(\eta,Y_k\right)$ does not depend on the particular decomposition $\om=\a_k\wedge df_k+f_k\,\be_k$ and on the particular equation of $Y_k$ (cf. \cite{d}).}
\end{rem}
The above remark allow us to define the residue of a logarithmic form $\eta$ on a arbitrary complex manifold $M$ along any codimension one smooth irreducible submanifold $Y$ contained in the pole divisor of $\eta$. In particular, we can define the iterated residue. Given $I=(i_1<...<i_k)\in \s_r^k$, set $X_I=(f_{i_1}=...=f_{i_k}=0)$ and $X_I^*=X_I\setminus\{0\}$. We define $Res(\eta,X_I)$ inductively. If $k=1$ then $Res(\eta,X_I)=Res(\eta,Y_{i_1})$ and for $k\ge2$, $Res(\eta,X_I)=Res\left(Res(\eta,Y_{i_k}),X_{I\setminus\{i_k\}}\right)$. This definition depends only of the ordering of the $f_{j's}$, that we will assume fixed.

\begin{ex}\label{ex:21}
{\rm If $\eta=\a\wedge\frac{df_{i_1}}{f_{i_1}}\wedge...\wedge\frac{df_{i_k}}{f_{i_k}}$, where $\a$ is holomorphic, then $Res(\eta,X_I)=\a|_{X_I}$, $I=(i_1<...<i_k)$. We leave the proof to the reader.}
\end{ex}

\begin{rem}\label{r:23}
{\rm Let $\eta=\frac{1}{f_1...f_k}\,\om$ be logarithmic as above, $Y_k=(f_k=0)$. We would like to observe the following facts:
\begin{itemize}
\item[(a).] If $Res(\eta,Y_k)=0$ then $f_k$ divides $\om$, or equivalently $f_k$ is not contained in the pole divisor of $\eta$. 
\item[(b).] If $p=1$ then $Res(\eta,Y_k)$ is a holomorphic function on $Y_k$.
\item[(c).] If $p\ge2$ then $Res(\eta,Y_k)$ is logarithmic on $Y_k$. Moreover, the pole divisor of $Res(\eta,Y_k)$ is $f_1...\wh{f}_k...f_r|_{Y_k}$.
\item[(d).] $d\eta$ is logarithmic and $Res(d\eta,Y_k)=d\,Res(\eta,Y_k)$. In particular, if $\eta$ is closed then $Res(\eta,Y_k)$ is closed.
\end{itemize}

We will prove (b) and (c). The proofs of (a) and (d) will be left to the reader. Write $\om=\a_k\wedge df_k+f_k\,\be_k$ (resp. $\om=g_k\,df_k+f_k\,\be_k$ if $p=1$) as before. It is sufficient to prove that if $\ell\ne k$ then $f_\ell|_{Y_k}$ divides $\a_k\wedge df_\ell|_{Y_k}$ (resp. $f_\ell|_{Y_k}$ divides $g_k|_{Y_k}$ if $p=1$). Note that $dim(Y_k\cap Y_\ell)\ge1$, because $n\ge p+2\ge3$. Therefore we can fix a point $m\in Y_k\cap Y_\ell$ where $df_k(m)\wedge df_\ell(m)\ne0$. Let $(U,z=(z_1,...,z_n))$ be a coordinate system around $m$ such that $f_k|_U=z_1$ and $f_\ell|_U=z_2$. 
Write $\om=\om_1\wedge dz_1+\om_2\wedge dz_2+\om_{12}\wedge dz_1\wedge dz_2+\te$, where $\om_1$ does not contain terms with $dz_2$, $\om_2$ does not contain terms with $dz_1$ and $\te$ does not contain terms in $dz_1$ or $dz_2$ (resp. $\om=\sum_jh_j\,dz_j$ if $p=1$).

Let us consider the case $p>1$. In this situation, $\om\wedge dz_1=\om_2\wedge dz_2\wedge dz_1+\te\wedge dz_1$ and $\te\wedge dz_1$ does not contain terms with $dz_1\wedge dz_2$, so that $z_1$ divides $\om_2$ and $\te$. Similarly, $z_2$ divides $\om_1$ and $\te$. Therefore, we can write $\om=z_2\,\wt\om_1\wedge dz_1+z_1\,\wt\om_2+\om_{12}\wedge dz_1\wedge dz_2+z_1\,z_2\,\wt\te$, which implies
\[
\om=\left(z_2\,\wt\om_1\wedge-\,\om_{12}\wedge dz_2\right)\wedge dz_1+z_1\,\left(\wt\om_2+z_2\,\wt\te\right)\,\,\implies
\]
\[
\a_k|_{Y_k\cap U}=\left(z_2\,\wt\om_1\wedge-\,\om_{12}\wedge dz_2\right)|_{Y_k\cap U}\implies\a_k\wedge dz_2|_{Y_k\cap U}=z_2\,\wt\om_2\wedge dz_2|_{Y_k\cap U}\,\,,
\]
which implies (c). 

In the case $p=1$, since $z_j$ divides $\om\wedge dz_j$, $\forall j$, then $z_1$ divides $h_j$ if $j>1$ and $z_2$ divides $h_1$, so that $\om|_U=z_2\,\wt{h}_1\,dz_1+z_1\,\sum_{j>1}\wt{h}_j\,dz_j$. Hence, $g_k|_{Y_k\cap U}=z_2\,\wt{h}_1|_{Y_k\cap U}$, which implies (b).}
\end{rem}

Let us prove lemma \ref{l:21} in the case $p=1$. As before, write $\eta=\frac{1}{f_1...f_r}\,\om$. The proof will be by induction on the number $r$ of components of the pole divisor.

{\it Formula {\rm(\ref{eq:4})} is true if $r=1$ and $p\ge1$.} When $\eta=\frac{1}{f_1}\,\om$ is logarithmic we have seen that $\om=\a_1\wedge df_1+f_1.\,\be_1$ (resp. $\om=g_1\,df_1+f_1\,\be_1$ if $p=1$). Hence $\eta=\a_1\wedge\frac{df_1}{f_1}+\be_1$ (resp. $\eta=g_1\,\frac{df_1}{f_1}+\be_1$ if $p=1$), as we wished.

{\it If $p=1$ and} (\ref{eq:4}) {\it is true for $r-1\ge1$ then it is true for $r$.} Let $\eta=\frac{1}{f_1...f_r}\,\om$ and $Q\sub\C^n$ be a polydisc where $f_1,...,f_r$ and $\om$ have representatives as before. As before we set $Y_j=(f_j=0)\sub Q\sub\C^n$.
We will use the following well known result in the case $n\ge3$:
\begin{lemma}\label{l:22}
Any holomorphic function $h\in\O(Y_j\setminus\{0\})$ has an extension $g\in\O(Q)$.
\end{lemma}
In fact, lemma \ref{l:22} is a particular case of theorem \ref{l:23} stated below and that will proved in \S\,\ref{s:6} (see remark \ref{r:24}).

Let $h=Res(\eta,Y_r)\in\O(Y_r\setminus\{0\})$. By lemma \ref{l:22}, $h$ has an extension $g_r\in\O(Q)$. The form $g_r\,\frac{df_r}{f_r}$ is logarithmic and $Res\left(g_r\,\frac{df_r}{f_r},Y_r\right)=h$. Therefore, the form $\wt\eta=\eta-g_r\,\frac{df_r}{f_r}$ is also logarithmic and $Res\left(\wt\eta,Y_r\right)=0$. In particular, $f_r$ is not a pole of $\wt\eta$ by remark \ref{r:23}. Since the pole divisor of $\wt\eta$ has $r-1$ irreducible components, by the induction hypothesis we can write
\[
\eta-g_r\,\frac{df_r}{f_r}=\wt\eta=\a_0+\sum_{j=1}^{r-1}g_j\,\frac{df_j}{f_j}\,\,\implies\,\,(\ref{eq:4})\,\,\text{in the case $p=1$}.
\]

The case $p\ge2$ is more involved, but the idea of the proof is the same as in the case $p=1$.
Before given the details let us sketch the proof.

Given $s\in\{0,1,...,r\}$ set $Y_s=(f_s=0)$ if $s\ge1$, $X_0=Q$ and $X_s=Y_1\cap...\cap Y_s$ if $s\ge1$. Set also $X_s^*=X_s\setminus\{0\}$, $0\le s\le r$.
Note that $X_s=\{0\}$ if $s\ge n$. On the other hand, if $1\le s\le n-1$ then $X_s$ is an analytic reduced germ of codimension $s$ and $X_s^*$ is a complex smooth manifold of dimension $n-s$.
The proof will involve two induction arguments. In order to state properly these arguments we need a definition.

Given $1\le s\le p\le n-2$ and $q\ge1$, we will say that $X_s^*$ satisfies the {\it $q$ decomposition property} if any logarithmic $q$-form $\te$ on $X_s^*$ with pole divisor on the zeroes of $f_{s+1}...f_r|_{X_s^*}:=\wt{f}_{s+1}...\wt{f}_r$ can be decomposed as in formula (\ref{eq:4}):
\[
\te=\a_0+\sum_{\ell=1}^{q-1}\left({\sum_{I\in\s_r^{\ell}}}\a_I\wedge\frac{d\wt{f}_{i_1}}{\wt{f}_{i_1}}\wedge...\wedge\frac{d\wt{f}_{i_\ell}}{\wt{f}_{i_\ell}}\right)+{\sum_{J\in\s_r^{\,q}}}g_J.\,\frac{d\wt{f}_{j_1}}{\wt{f}_{j_1}}\wedge...\wedge\frac{d\wt{f}_{j_q}}{\wt{f}_{j_q}}\,\,,
\]
where $\a_0$ is a holomorphic q-form, the $\a_{I's}$ are holomorphic $(q-\ell)$-forms on $X_s^*$ and the $g_{J's}$ are holomorphic functions on $X_s^*$.
We resume below the main steps in the arguments.
\begin{itemize}
\item[$1^{st}$ step.] If $0\le s\le p-1$ then $X_s^*$ satisfies the 1 decomposition property.
\item[$2^{nd}$ step.] If $2\le q\le p-s$, where $s\ge0$, and $X_{s+1}^*$ satisfies the $q-1$ decomposition property then $X_s^*$ satisfies the $q$ decomposition property.
\end{itemize}
The $1^{st}$ and $2^{nd}$ steps above will be proved by induction on the number of $r\ge1$ of factors in the pole divisor $f_1...f_r$.
In the proof we will use the following result:

\begin{thm}\label{l:23}
Let $X$ be a germ at $0\in\C^n$ with a s.o.s at $0\in X$ and $dim_\C(X)=k$, where $2\le k\le n-1$. Set $X^*:=X\setminus\{0\}$.
Then there are representatives of $X$ and $X^*$ in a polydisc $Q\sub\C^n$, denoted by the same letters, such that:
\begin{itemize}
\item[(a).] If $0\le q\le k-2$ then any form $\wt\a\in\Om^q(X^*)$ can be extended to a form $\a\in\Om^q(Q)$.
\item[(b).] If $q\ge1$, $\ell\ge0$ and $1\le q+\ell\le k-2$ then $H^q(X^*,\Om^\ell)=0$.
\end{itemize}
\end{thm}

\begin{rem}\label{r:24}
{\rm Note that lemma \ref{l:22} is a particular case of theorem \ref{l:23}.}
\end{rem}

Theorem \ref{l:23} implies that, if $X_s^*$ is as before and $0\le q\le n-s-2$, then any holomorphic $q$-form on $X_s^*$ can be extended to a holomorphic $q$-form on $Q$.
The proof of theorem \ref{l:23} will be done in \S\, \ref{s:6}.
Let us finish the proof of lemma \ref{l:21} assuming theorem \ref{l:23}.

\vskip.1in

{\it Proof of the 1$^{st}$ step.} Is similar to the case $p=1$ done above with lemma \ref{l:22}, which corresponds to the case $s=0$. Therefore, we will assume $1\le s\le p-1$. Note that the assertion is trivially true if $r=s$, because in this case the 1-form is holomorphic.

Assume that the assertion is true for any logarithmic 1-form on $X_s^*$ with pole divisor containing $k-1\ge0$ functions in the set $\{\wt{f_j}=f_j|_{X_s^*}\,|\,s+1\le j\le r\}$. 
Let $\te$ be a logarithmic $1$-form on $X_s^*$ with pole divisor $\wt{f}_{s+1}...\wt{f}_{s+k}$. By remark \ref{r:23}, $\wt{g}_{s+1}:=Res(\te,X_{s+1}^*)\in\O(X_{s+1}^*)$. By theorem \ref{l:23}, $\wt{g}_{s+1}$ admits an extension $g_{s+1}\in\O(Q)$. In particular, $\wh{g}_{s+1}:=g_{s+1}|_{X_s^*}$ is a holomorphic extension of $\wt{g}_{s+1}$ on $X_s^*$. 
Let $\wt{\te}:=\wh{g}_{s+1}\,\frac{d\wt{f}_{s+1}}{\wt{f}_{s+1}}\,$. Then $\wt{\te}$ is logarithmic and $Res(\wt\te,X_{s+1}^*)=Res(\te,X_{s+1}^*)$. In particular, $\wt{f}_{s+1}$ is not contained in the pole divisor of $\te-\wt\te$, by (a) of remark \ref{r:23}. By the induction hypothesis, $\te-\wt\te$ can be decomposed as in (\ref{eq:4}):
\[
\te-\wt\te=\a_0+\sum_{j=2}^k\wh{g}_{s+j}\,\frac{d\wt{f}_{s+j}}{\wt{f}_{s+j}}\,\,\implies\,\,\text{\it $1^{st}$ step.}
\]

{\it Proof of the 2$^{nd}$ step.} The proof is again by induction on the number $r-s$ of factors of the pole divisor.
The assertion is trivially true if $r=s$. Assume that the assertion is true for any logarithmic $q$-form, $2\le q\le p-s$, on $X_s^*$ with pole divisor containing $k-1\ge0$ functions in the set $\{\wt{f_j}=f_j|_{X_s^*}\,|\,s+1\le j\le r\}$. Let $\te$ be a logarithmic $q$-form on $X_s^*$ with pole divisor $\wt{f}_{s+1}...\wt{f}_{s+k}$. By remark \ref{r:23} the $(q-1)$-form $\mu:=Res(\te,X_{s+1}^*)$ is logarithmic on $X_{s+1}^*$ with pole divisor $\wh{f}_{s+2}...\wh{f}_{s+k}:=\wt{f}_{s+2}...\wt{f}_{s+k}|_{X_{s+1}^*}$ (or holomorphic if $k=1$). Since $X_{s+1}^*$ satisfies the $q-1$ decompostion property, we can write $\mu=$
\[
\a_0+\sum_{t=1}^{q-2}{\sum_{I\in\s_r^t}}\a_I\wedge\frac{d\wh{f}_{s+i_1-1}}{\wh{f}_{s+i_1-1}}\wedge...\wedge\frac{d\wh{f}_{s+i_t-1}}
{\wh{f}_{s+i_t-1}}+{\sum_{J\in\s_r^{q-1}}}g_J.\,\frac{d\wh{f}_{s+j_1-1}}{\wh{f}_{s+j_1-1}}\wedge...\wedge\frac{d\wh{f}_{s+j_{(q-1)}-1}}
{\wh{f}_{s+j_{(q-1)}-1}}
\]
where $\a_0$ and the $\a_{I's}$ are holomorphic forms on $X_{s+1}^*$ and the $g_{J's}$ are holomorphic functions on $X_{s+1}^*$. By theorem \ref{l:23} each $\a_{I}$ (resp. each $g_J$) has a holomorphic extension $\wt\a_I$ (resp. $\wt{g}_J$) on $X_s^*$. Therefore, $\mu$ has a logarithmic extension $\wt\mu$ on $X_s^*$, $\wt\mu=$
\[
\wt\a_0+\sum_{t=1}^{q-2}{\sum_{I\in\s_r^t}}\wt\a_I\wedge\frac{d\wt{f}_{s+i_1-1}}{\wt{f}_{s+i_1-1}}\wedge...\wedge\frac{d\wt{f}_{s+i_t-1}}
{\wt{f}_{s+i_t-1}}+{\sum_{J\in\s_r^{q-1}}}\wt{g}_J.\,\frac{d\wt{f}_{s+j_1-1}}{\wt{f}_{s+j_1-1}}\wedge...\wedge\frac{d\wt{f}_{s+j_{(q-1)}-1}}
{\wt{f}_{s+j_{(q-1)}-1}}
\]
Therefore, $\te_1:=\wt\mu\wedge\frac{d\wt{f}_{s+1}}{\wh{f}_{s+1}}$ is logarithmic on $X_s^*$ and
\[
Res(\te_1,X_{s+1}^*)=Res(\te,X_{s+1}^*)\,\,\implies\,\,Res(\te-\te_1,X_{s+1}^*)=0\,\,.
\]
Hence, $\wt{f}_{s+1}$ is not contained in the pole divisor of $\te-\te_1$. By the induction hypothesis, $\te-\te_1:=\te_2$ admits a decomposition as in (\ref{eq:4}), and so $\te=\te_1+\te_2$ admits a decomposition as in (\ref{eq:4}).
This finishes the proof of lemma \ref{l:21}.
\qed

\subsection{Proof of theorem \ref{t:1}}\label{ss:22}

In the proof of theorem \ref{t:1} we will use theorem \ref{l:23} and Hamm's generalization of Milnor's theorem (cf. \cite{ha}, \cite{ha1}, \cite{mi} and \cite{js}):

\begin{thma}\label{t:21}
Let $X=(f_1=...=f_\ell=0)$ be a germ at $0\in\C^m$ of a complete intersection with an isolated singularity at $0$, so that $dim_\C(X)=m-\ell:=n$. Then there exist representatives of $f_1,...,f_\ell$ and $X$ defined in a ball $B_\ep=B(0,\ep)$, denoted by the same letters, such that:
\begin{itemize}
\item[(a).] $X^*=X\setminus\{0\}$ is rectratible to the link $K:=X\cap \mathbb{S}_\ep^{2m-1}$, $\mathbb{S}_\ep^{2m-1}=\pa B_\ep$.
\item[(b).] If $n\ge3$ then $K$ is $(n-2)$-connected. In particular, $X^*$ is connected and $H_{DR}^k(X^*)=\{0\}$ if $1\le k\le n-2$.
\item[(c).] If $n=2$ then $X^*$ is connected.
\end{itemize}
\end{thma}

When $n=1$, $X^*$ is not necessarily connected, as shows the example $X=(x^2+y^2+z^2=z=0)\sub\C^3$.

\vskip.1in

Let $\eta$ be a germ at $0\in\C^n$ of a closed logarithmic $p$-form, $1\le p\le n-2$, with pole divisor $f_1...f_r$ with a strictly ordinary singularity outside $0$. According to lemma \ref{l:21} we can write $\eta$ as a sum of a holomorphic $p$-form $\a_0$, and "monomial" $p$-forms of the type $\a_I\wedge\frac{df_{i_1}}{f_{i_1}}\wedge...\wedge\frac{df_{i_s}}{f_{i_s}}$, or $g_J\,\frac{df_{j_1}}{f_{j_1}}\wedge...\wedge\frac{df_{j_p}}{f_{i_p}}$, where $I\in\s_r^s$ and $J\in\s_r^p$.

Given a monomial $\mu=\a_I\wedge\frac{df_{i_1}}{f_{i_1}}\wedge...\wedge\frac{df_{i_s}}{f_{i_s}}$ we define the {\it pseudo depth} of $\mu$ as $\wt{Dp}\,(\mu)=s$. Given $\eta=\sum_{j=1}^m\mu_j$, where the $\mu_j$ are monomials as above, we set $\wt{Dp}\,(\eta)=max\{\wt{Dp}\,(\mu_j)\,|\,1\le j\le m\}$.

Observe that $\wt{Dp}$, as defined above, is not well defined. For instance, if $g\in\left<f_1,...,f_p\right>$, the ideal generated by $f_1,...,f_p$, $g=\sum_{j=1}^ph_j.\,f_j$, then
\[
g\,\frac{df_{1}}{f_{1}}\wedge...\wedge\frac{df_{p}}{f_{p}}=\sum_j\a_j\wedge\frac{df_{1}}{f_{1}}\wedge...\wedge\wh{\frac{df_{j}}{f_{j}}}\wedge\frac{df_{p}}{f_{p}}\,\,,
\]
where $\a_j=\pm\,h_j\,df_j$, $1\le j\le p$. Therefore, if $\eta$ is a logarithmic form as above, then we define its depth as
\[
Dp\,(\eta)=min\left\{\wt{Dp}\,\left(\sum_j\mu_j\right)\,|\,\eta=\sum_j\mu_j\,,\,
\text{where the $\mu_{j's}$ are monomials}\right\}\,.
\]
When $\eta$ is holomorphic we define $Dp\,(\eta)=0$.

\begin{claim}\label{cl:21}
Let $\eta$ be a germ at $0\in\C^n$ of logarithmic closed $p$-form, $1\le p\le n-2$. If $Dp\,(\eta)=p$ there exists a collection $(\la_J)_{J\in\s_r^p}$, $\la_J\in\C$, such that
\[
Dp\,\left(\eta-\sum_{J\in\s_r^p}\la_J\,\frac{df_{j_1}}{f_{j_1}}\wedge...\wedge\frac{df_{j_p}}{f_{j_p}}\right)\le p-1\,\,.
\]
\end{claim}

{\it Proof.}
If $Dp\,(\eta)=p$ then the decomposition of $\eta$ as in (\ref{eq:4}) contains at least one monomial of the form $\mu_J=g_J\,\frac{df_{j_1}}{f_{j_1}}\wedge...\wedge\frac{df_{j_p}}{f_{j_p}}$, where $g_J\notin\left<f_{j_1},...,f_{j_p}\right>$.
As before, set $X_J:=(f_{j_1}=...=f_{j_p}=0)$ and $X_J^*=X_J\setminus\{0\}$. 
We assert that $g_J|_{X_J^*}$ is a constant $\la_J\in \C^*$.

In fact, since $dim_\C X_J=n-p\ge2$, $X_J^*$ is connected, by theorem \ref{t:21}. Note that $Res(\mu,X_J)=Res(\eta,X_J)=g_J|_{X_J^*}$ (see example \ref{ex:21}).
Since $\eta$ is closed, we have $Res(d\eta,X_J)=dg_J|_{X_J^*}=0$.
Hence, $g_J|_{X_J^*}=\la_J\in\C$. On the other hand, if $\la_J=0$ then $g_J|_{X_J}=0$ and since $X_J$ is a complete intersection we get $g_J\in\left<f_{j_1},...,f_{j_p}\right>$, a contradiction.

Let $\mu_J:=\la_J\,\frac{df_{j_1}}{f_{j_1}}\wedge...\wedge\frac{df_{j_p}}{f_{j_p}}$. Note that $\eta-\mu_J$ is still logarithmic, closed and does not contain terms multiples of $\frac{df_{j_1}}{f_{j_1}}\wedge...\wedge\frac{df_{j_p}}{f_{j_p}}$. By repeating this procedure finitely many times we can find the collection $(\la_J)_{J\in\s_r^p}$ as in the statement of the claim.
\qed

\begin{claim}\label{cl:22}
{\rm Let $\eta$ be logarithmic closed p-form with pole divisor $f_1...f_r=0$, with a s.o.s at $0\in\C^n$. If $Dp\,(\eta)<p$ then $\eta$ is exact: $\eta=d\Te$, where $\Te$ is either zero, or is logarithmic with pole divisor contained in $f_1...f_r=0$.}
\end{claim}

{\it Proof.}
The proof will be by induction on the depth of $\eta$. If $Dp\,(\eta)=0$ then $\eta$ is holomorphic and so it is exact by Poincaré lemma.

Assume that any closed logarithmic $p$-form $\om$ with $Dp\,(\om)\le q-1<p-1$ is exact: $\om=d\te$ with $\te$ logarithmic as above. 
Let $\eta$ be a logarithmic $p$-form with pole divisor $f_1...f_r=0$ with $Dp\,(\eta)=q<p$. By the definition of depth, when we write $\eta$ as in (\ref{eq:4}) then we get
\[
\eta=\sum_{I\in\s_r^q}\a_I\wedge\frac{df_{i_1}}{f_{i_1}}\wedge...\wedge\frac{df_{i_q}}{f_{i_q}}+\be\,\,,
\]
where $\be$ is logarithmic and $Dp\,(\be)<q$. Recall that, if $I=(i_1,...,i_q)\in\s_r^q$ then $X_I=(f_{i_1}=...=f_{i_q}=0)$ and $X_I^*=X_I\setminus\{0\}$. As the reader can check, we have
\[
Res\,(\eta,X_I)=\a_I|_{X_I}:=\wt\a_I\in\Om^{p-q}(X_I^*)\,\,,\,\,\forall I\in\s_r^q\,\,,
\]
where $\wt\a_I$ is closed, by remark \ref{r:23}.
Now, we use theorem \ref{t:21} and (c) of theorem \ref{l:23}: since $dim(X_I^*)=n-q$ we get $H^k_{DR}(X^*)=0$ if $1\le k\le n-q-2$. But, $p\le n-2$ and so $p-q\le n-q-2=dim_\C(X^*)-2$ which implies that $\wt\a_I\in\Om^{p-q}(X_I^*)$ is exact: $\wt\a_I=d\wt\be_I$, where $\wt\be_I$ is in principle a $C^\infty$ $(p-q-1)$-form. However, the fact that $H_{\ov\pa}^{s\,r}(X_I^*)\simeq H^r(X_I^*,\Om^s)=0$ if $r+s=p-q-1$ implies that we can assume $\wt\be_I\in \Om^{p-q-1}(X^*)$ (cf. \cite{gr}).

Therefore, there are $(p-q-1)$-forms $\wt\be_I\in\Om^{p-q-1}(X_I^*)$ such that $\wt\a_I=d\wt\be_I$, $\forall I\in\s_r^q$. By (b) of theorem \ref{l:23} each form $\wt\be_I$ admits an extension $\be_I\in\Om^{p-q-1}(Q)$, where $Q$ is some polydisc of $\C^n$ where $X_I$ has a representative.
Define a logarithmic form $\mu$ by
\[
\mu=\sum_{I\in\s_r^q}\be_I\wedge \frac{df_{i_1}}{f_{i_1}}\wedge...\wedge\frac{df_{i_q}}{f_{i_q}}
\]
so that
\[
d\mu=\sum_{I\in\s_r^q}d\be_I\wedge \frac{df_{i_1}}{f_{i_1}}\wedge...\wedge\frac{df_{i_q}}{f_{i_q}}\,\,\implies\,\,Res(d\mu,X_I)=d\be_I|_{X_I}=Res(\eta,X_I)\,\,,\,\,\forall I\in\s_r^q\,\,.
\]
In particular, $Res(\eta-d\mu,X_I)=0$ for all $I\in\s_r^q$, and this implies that $Dp\,(\eta-d\mu)<q$. 
Finally, since $\eta-d\mu$ is closed the induction hypothesis implies that $\eta-d\mu=d\te$, where either $\te=0$, or $\te$ is logarithmic with pole divisor contained in $f_1...f_r=0$.
This finishes the proof of claim \ref{cl:22} and of theorem \ref{t:1}.
\qed

\subsection{Proof of corollary \ref{c:1}}\label{ss:23}

Let $\eta$ be a logarithmic $p$-form on $\p^n$, where $p\le n-1$. Assume that the pole divisor of $\eta$ is a hypersurface with normal crossing singularities and smooth irreducible components, so that by Deligne's theorem (theorem \ref{t:11}) $\eta$ is closed.

Let $\Pi\colon\C^{n+1}\setminus\{0\}\to\p^n$ be the canonical projection. The pull-back $\Pi^*(\eta)$ can be extended to a closed logarithmic $p$-form on $\C^{n+1}$ which is called the expression of $\eta$ in homogeneous coordinates. We will denote this extension by $\wt\eta$. The pole divisor of $\wt\eta$ is of course the pull-back of the pole divisor of $\eta$, so that it can be represented in $\C^{n+1}$ by $f_1...f_r$, where $f_j$ is a homogeneous polynomial of degree $d_j$, $1\le j\le r$.
In particular, we can write
\[
\wt\eta=\frac{1}{f_1...f_r}\,\sum\,g_J\,dz^J
\]
where $dz^J=dz_{j_1}\wedge...\wedge dz_{j_p}$ and $g_J$ is a homogeneous polynomial.
Using that $\wt\eta$ is invariant by any homothety $H_t(z)=t.\,z$, $\forall t\in\C^*$, and with a straighforward computation we see that $g_J$ is homogeneous of degree $deg(g_J)=deg(f_1...f_r)-p$. This implies that the coefficients of $\wt\eta$ are meromorphic homogenous of degree $-p$.

Now, the hypothesis on the pole divisor of $\eta$ implies that the pole divisor of $\wt\eta$, $f_1...f_r$, has s.o.s outside $0\in\C^{n+1}$. 
Therefore, by theorem \ref{t:1} we have
\[
\wt\eta=\sum_{I\in\s_r^p}\la_I\,\frac{df_{i_1}}{f_{i_1}}\wedge...\frac{df_{i_p}}{f_{i_p}}\,+\,d\Te\,\,,
\]
where $\Te$ is logarithmic with pole divisor contained in $(\eta)_\infty$. It is enough to prove that $d\Te=0$.

The proof of theorem \ref{t:1} implies that the monomials of $\Te$ have depth $<p$ and are, either of the form $\a\wedge\frac{df_{j_1}}{f_{j_1}}\wedge...\wedge\frac{df_{j_q}}{f_{j_q}}$, where $\a$ is a $(p-q-1)$-form, or of the form $g.\,\frac{df_{j_1}}{f_{j_1}}\wedge...\wedge\frac{df_{j_{p-1}}}{f_{j_{p_1}}}$, where $g$ is a holomorphic function. In particular, the monomials of $d\Te$ are, either of the form $d\a\wedge\frac{df_{j_1}}{f_{j_1}}\wedge...\wedge\frac{df_{j_q}}{f_{j_q}}$, or of the form $dg\wedge\frac{df_{j_1}}{f_{j_1}}\wedge...\wedge\frac{df_{j_{p-1}}}{f_{j_{p_1}}}$. In both cases, the meromorphic degree of the coefficients of the monomial is $>-p$ and this implies that $d\Te=0$.
The proof that $i_R\,\wt\eta=0$ follows from the fact that $D\Pi(z).R(z)=0$ for all $z\in\C^{n+1}\setminus\{0\}$. Finally, $i_R\,\wt\eta=0$ implies that $r\ge p+1$, as the reader can check.
\qed

\section{Decomposition of logarithmic foliations}\label{ss:3}
In this section we will prove theorem \ref{t:2} and corollary \ref{c:2}. Since the case $r=p+1$ is the easier one, we will do it, together with the proof of corollary \ref{c:2}, in \S \ref{ss:31}.
In \S \ref{ss:32} we prove the theorem in the case $p=2$. In \S \ref{ss:33} we will see that the proof of (b) can be reduced to the case of 2-dimensional foliations (in which $p=n-2$). The proof of the case $r=p+2$ will be done in \S \ref{ss:34} and the proof of (b) in \S \ref{ss:35}.

\subsection{Proof of the case $r=p+1$ and of corollary \ref{c:2}.}\label{ss:31}

The proof will be based in the remark that a $p$-vector $\Om$ in a vector space $V$ of dimension $p+1$ is always decomposable. In fact, if $\{v_1,...,v_{p+1}\}$ is a basis of $V$, then we can write
\[
\Om=\sum_{j=1}^{p+1}a_j\,v_1\wedge...\wedge\wh{v_j}\wedge...\wedge v_{p+1}\,\,,\,\,a_j\in K\,\,,\,\,1\le j\le p+1\,\,.
\]
Since $\Om\ne0$, we can assume that $a_1\ne0$. Dividing $\Om$ by $a_1$ if necessary, we can assume that $a_1=1$.

Let $\{g_1,...,g_{p+1}\}$ be dual basis of the basis $\{v_1,...,v_{p+1}\}$; $g_j(v_i)=\d_{ij}$. If $X=g_1+\sum_{j=2}^{p+1}(-1)^{j-1}a_j\,g_j$ then $\Om=i_X\,v_1\wedge...\wedge v_{p+1}$. Now, if we set $\te_j:=v_j+(-1)^ja_j\,v_1$, $2\le j\le p+1$, then $i_X\,\te_j=0$ and the reader can verify that $\Om=\te_2\wedge...\wedge \te_{p+1}$.

Let $\wt\eta$ be the extension of $\Pi^*(\eta)$ to $\C^{n+1}$, as in corollary \ref{c:1}. Let $f_1...f_{p+1}$ be the pole divisor $(\wt\eta)_\infty$, so that
\begin{equation}\label{eq:5}
\wt\eta=\sum_{j=1}^{p+1}\la_j\frac{df_1}{f_1}\wedge...\wedge\wh{\frac{df_j}{f_j}}\wedge...\wedge\frac{df_{p+1}}{f_{p+1}}\,\,.
\end{equation}
By the above remark $\wt\eta$ is decomposable: if we assume $\la_1\ne0$ then there exist $\mu_2,...,\mu_{p+1}\in\C$ such that, if we set $\te_j=\frac{df_j}{f_j}-\mu_j\,\frac{df_1}{f_1}$ then
\[
\wt\eta=\la_1\,\te_2\wedge...\wedge\te_{p+1}\,\,.
\]

We assert that $\mu_j\in\Q_+$, $2\le j\le p+1$.
In fact, from $i_R\,\wt\eta=0$ we get
\[
i_R\,(\te_1\wedge...\wedge\te_{p+1})=0\,\,\implies\,\,\sum_{j=2}^{p+1}\,(-1)^j\,i_R\,(\te_j).\,\te_2\wedge...\wedge\wh{\te_j}\wedge...\wedge\te_{p+1}=0\,\,\implies
\]
\[
i_R\,(\te_j)=0\,\,,\,\,2\le j\le p+1\,\,\implies\,\,\mu_j=\frac{deg(f_j)}{deg(f_1)}:=\frac{d_j}{d_1}\,\in\Q_+\,\,.
\]
In particular, the rational function $f_j^{d_1}/f_1^{d_j}$ is a first integal of $\te_j$; $d\left(f_j^{d_1}/f_1^{d_j}\right)\wedge\te_j=0$, $2\le j\le p+1$. This of course implies that $F=\left(f_1^{k_1},...,f_{p+1}^{k_{p+1}}\right)$ is a first integral of $\wt\eta$ if $k_j:=d_1...d_{p+1}/d_j$.
\qed

\subsection{Proof of theorem \ref{t:2} in the case $p=2$: foliations of codimension two.}\label{ss:32}

Let $\fa$ be a logarithmic foliation of codimension two on $\p^n$ defined by a logarithmic 2-form $\wt\eta\in\L^2_\fa(f_1,...,f_r)$. Note that the hypothesis $p=2\le n-2$ implies that $n\ge4$.
\begin{rem}\label{r:31}
{\rm The condition of local decomposability of $\wt\eta$ outside the singular set is equivalent to $\wt\eta\wedge\wt\eta=0$.
This is a consequence of the fact that a two vector $\te$ on a complex vector space is decomposable if, and only if, $\te\wedge\te=0$.

In particular, we have
\[
\L^2_\fa(f_1,...,f_r)=\{\om\in\L^2_R(f_1,...,f_r)\,|\,\om\wedge\om=0\}\,\,.
\]}
\end{rem}

As we have seen, a form $\om\in\L^2_R(f_1,...,f_r)$ can be written as
\begin{equation}\label{eq:eta}
\om=\sum_{1\le i<j\le r}\,\mu_{ij}\,\frac{df_i}{f_i}\wedge\frac{df_j}{f_j}\,\,.
\end{equation}
As the reader can check,
\[
\om\wedge\om=\sum_{1\le i<j<k<\ell\le r}
\,2\,\Psi(\mu_{ij},\mu_{k\ell},\mu_{ik},\mu_{j\ell},\mu_{i\ell},\mu_{jk})\,\frac{df_i}{f_i}\wedge\frac{df_j}{f_j}\wedge\frac{df_k}{f_k}\wedge\frac{df_\ell}{f_\ell}\,\,,
\] 
where $\Psi(a,b,c,d,e,f)=a\,b-c\,d+e\,f$. If $\om\wedge\om=0$ their numerical residues must vanish (see remark \ref{r:12}). This implies that $\L^2_\fa(f_1,...,f_r)$ is isomorphic to the algebraic subset $\mathcal{A}$ of $\C^{r(r-1)/2}$ defined by
\[
\mathcal{A}=\left\{(\la_{ij})_{1\le i<j\le r}\,|\,\Psi(\la_{ij},\la_{k\ell},\la_{ik},\la_{j\ell},\la_{i\ell},\la_{jk})=0\,\,,\,\,\forall\,1\le i<j<k<\ell\le r\right\}\,\,,
\]
where the isomorphism is given by
\[
\left(\la_{ij}\right)_{1\le i<j\le r}\in\mathcal{A}\,\mapsto\,\sum_{1\le i<j\le r}\,\la_{ij}\,\frac{df_i}{f_i}\wedge\frac{df_j}{f_j}
\]
On the other hand, if we fix a base $\{e_1,...,e_r\}$ of $\C^r$, a 2-vector $\te$ on $\C^r$ can be written as
\[
\te=\sum_{1\le i<j\le r}\,a_{ij}\,\,e_i\wedge e_j\,\,.
\]
Since
\[
\te\wedge\te=\sum_{1\le i<j<k<\ell\le r}
\,2\,\Psi(a_{ij},a_{k\ell},a_{ik},a_{j\ell},a_{i\ell},a_{jk})\,\,e_i\wedge e_j\wedge e_k\wedge e_\ell\,\,,
\]
we obtain $\te\wedge\te=0$ if, and only if, $\left(a_{ij}\right)_{1\le i<j\le r}\in\mathcal{A}$. Now, if $\te\wedge\te=0$ then $\te$ is decomposable: $\te=\a\wedge\be$, where $\a,\be\in\C^r$. In fact, if $\te\ne0$ let $u,v$ be in the dual of $\C^r$ and such that $\te(u,v)\ne0$. Then
\[
0=i_u\,(\te\wedge\te)=2\,i_u(\te)\wedge\te\,\,\implies\,\,\te=c.\,i_u(\te)\wedge i_v(\te)\,,\,c=1/\te(u,v)\,.
\]
Finally, if $\om$ is as in (\ref{eq:eta}) and satisfies $\om\wedge\om=0$ then the 2-vector $\te=\sum_{i<j}\mu_{ij}\,e_i\wedge e_j$ is decomposable: $\te=\a\wedge\be$, $\a=\sum_ia_i\,e_i$ and $\be=\sum_jb_j\,e_j$, so that $\om=\om_1\wedge\om_2$, $\om_1=\sum_i\,a_i\,\frac{df_i}{f_i}$ and $\om_2=\sum_j\,b_j\,\frac{df_j}{f_j}$. Moreover, if $i_R\,\om=0$ then $i_R\,\om_1=i_R\,\om_2=0$ because
\[
0=i_R\,(\om_1\wedge\om_2)=i_R\,\om_1.\,\om_2-i_R\,\om_2.\,\om_1\,\,\implies\,\,i_R\,\om_1=i_R\,\om_2=0\,\,\qed
\]

\subsection{Some remarks.}\label{ss:33}
From now on, we fix homogeneous polynomials $f_1,...,f_r\in \C[z_0,...,z_n]$, where $r>p+1$, the divisor $f_1,...,f_r$ has s.o.s outside $0\in\C^{n+1}$ and $deg(f_j)=d_j$,
$1\le j\le r$.
Recall that $\L^p(f_1,...,f_r)$ denotes the set of logarithmic p-forms that can be written as below:
\begin{equation}
\wt\eta=\sum_{I\in\s_r^p}\la_I\,\frac{df_{i_1}}{f_{i_1}}\wedge...\wedge\frac{df_{i_p}}{f_{i_p}}\,\,,\,\,\la_I\in\C\,,\,\forall I\in\s_r^p\,.
\end{equation}\label{eq:7.1}

Given a base $\{du_1,...,du_r\}$ of $\C^{r*}$ there exists an unique linear map $\Phi^p\colon\bigwedge^p(\C^{r*})\to\L^p(f_1,...,f_r)$ such that
\[
\Phi^p(du_{i_1}\wedge...\wedge du_{i_p})=\frac{df_{i_1}}{f_{i_1}}\wedge...\wedge\frac{df_{i_p}}{f_{i_p}}\,\,.
\]
\begin{lemma}\label{l:31}
$\Phi^p$ is an isomorphism for all $p\ge1$.
Moreover, if $\a\in\bigwedge^p(\C^{r*})$ and $\be\in\bigwedge^q(\C^{r*})$ then
\begin{equation}\label{eq:8}
\Phi^{p+q}(\a\wedge\be)=\Phi^p(\a)\wedge\Phi^q(\be)\,\,.
\end{equation}
\end{lemma}

{\it Proof.}
On one hand, it is clear that $\Phi^p$ is surjective. On the other hand, if $\wt\eta=\sum_{I\in \s_r^p}\la_I\frac{df_{i_1}}{f_{i_1}}\wedge...\wedge\frac{df_{i_p}}{f_{i_p}}$ then each numerical residue $\la_I$, $I\in\s_r^p$, can be calculated by an integral as in remark \ref{r:12}:
\[
\la_I=\frac{1}{(2\pi i)^p}\int_{T_\ep^p}\eta\,\,.
\]
It follows that
\[
\sum_{I\in \s_r^p}\la_I\frac{df_{i_1}}{f_{i_1}}\wedge...\wedge\frac{df_{i_p}}{f_{i_p}}\equiv0\,\,\iff\,\,\la_I=0\,\,,\,\,\forall I\in\s_r^p
\]
and so $\Phi^p$ is injective.

Finally, formula (\ref{eq:8}) is consequence of
\[
\Phi^{p+q}\left((du_{i_1}\wedge...\wedge du_{i_p})\wedge(du_{j_1}\wedge...\wedge du_{j_q})\right)=\frac{df_{i_1}}{f_{i_1}}\wedge...\wedge\frac{df_{i_p}}{f_{i_p}}\wedge\frac{df_{j_1}}{f_{j_1}}\wedge...\wedge\frac{df_{j_q}}{f_{j_q}}=
\]
\[
=\Phi^p(du_{i_1}\wedge...\wedge du_{i_p})\wedge\Phi^q(du_{j_1}\wedge...\wedge du_{j_q})\,\,\qed
\]

\begin{rem}\label{r:32'}
{\rm Given a p-form $\a\in\bigwedge^p(\C^{r*})$ its kernel is defined as
\[
ker(\a)=\{v\in\C^r\,|\,i_v\,\a=0\}\,.
\]
We say that $\a\in\bigwedge^p(\C^{r*})$ is totally decomposable if there are p 1-forms $\a_1,...,\a_p$ such that $\a=\a_1\wedge...\wedge\a_p$.
It is well known that:
\begin{itemize}
\item[(I).] $\a=\a_1\wedge...\wedge\a_p$ is totally decomposable if, and only if, $dim(ker(\a))=r-p$.
\item[(II).] If $\a=\a_1\wedge...\wedge\a_p$ then $ker(\a)=\bigcap_{j=1}^p ker(\a_j)$.
\item[(III).] The projectivization of the set of totally decomposable p-forms of $\bigwedge^p(\C^{r*})$ is isomorphic to the grassmanian of p planes through the origin in $\C^r$. In particular, it is an algebraic subset of $\p\left(\bigwedge^p(\C^{r*})\right)$.
\end{itemize}}
\end{rem}

Recall that $\wt\eta\in\L^p_{td}(f_1,...,f_r)$ if it is totally decomposable into logarithmic forms (t.d.l.f). 
An easy consequence of lemma \ref{l:31} and of (III) of remark \ref{r:32'} is the following:

\begin{cor}\label{c:31}
Let $p\ge2$. A p-form $\wt\eta\in\L^p(f_1,...,f_r)$ is {\rm t.d.l.f} if, and only if, there are 1-forms $\a_1,...,\a_p\in\C^{r*}$ such that
\[
\wt\eta=\Phi^p(\a_1\wedge...\wedge \a_p)\,\,.
\]
In particular, $\L^p_{td}(f_1,...,f_r)$ is an algebraic irreducible subset of $\L^p_R(f_1,...,f_r)$.
\end{cor}

Another consequence of lemma \ref{l:31} is that (b) in the statement of theorem \ref{t:2} can be reduced to the case of 2-dimensional foliations.
Let $\Si\simeq\p^q$ be a $q$-plane linearly embedded in $\p^n$. We say that $\Si$ is in general position with respect to the divisor $f_1...f_r$ if for all $J=\{j_1,...,j_k\}\sub\{1,...,r\}$ then $\Si$ is transverse to $\bigcap_{j\in J}\Pi(f_j=0)$. By transversality theory, the set of $q$-planes of $\p^n$ in general position with respect to $f_1...f_r$ is a Zariski open and dense subset of the grassmanian of $q$-planes on $\p^n$.

\begin{rem}\label{r:32}
{\rm Let $\wt\eta\in\L^p_R(f_1,...,f_r)$.  Let $\Si$ be $(p+k)$-plane of $\p^n$ in general position with respect to $f_1...f_r$, $p<p+k<n$, and $\wt\Si$ be a $p+k+1$ plane through $0\in\C^{n+1}$ such that $\Pi(\wt\Si)=\Si$. Then $\wt\eta|_{\wt\Si}$ is a logarithmic $p$-form on $\wt\Si$. An easy consequence of lemma \ref{l:31} and corollary \ref{c:31} is the following:}
\end{rem}

\begin{cor}\label{c:32}
Let $\wt\eta$, $\Si$ and $\wt\Si$ be as in remark \ref{r:32}. Then $\wt\eta$ is {\rm t.d.l.f} if, and only if, $\wt\eta|_{\wt\Si}$ is {\rm t.d.l.f.}
\end{cor}

Taking $k=2$ in the above statement, we reduce the proofs of the case $r=p+2$ and of (b) in theorem \ref{t:2} to the case of 2-dimensional foliations.
From now on, we will assume that $\wt\eta=\Pi^*(\eta)\in\L_\fa^p(f_1,...,f_r)$ and that $n=p+2$. By \S \ref{ss:31} we will assume also that $r\ge p+2$.
As we have seen, we can write
\begin{equation}\label{eq:la}
\wt\eta=\sum_{I\in\s_r^p}\la_I\frac{df_{i_1}}{f_{i_1}}\wedge...\wedge\frac{df_{i_{n-2}}}{f_{i_{n-2}}}\,\,.
\end{equation}
The foliation $\fa_{\eta}$ is defined in homogeneous coordinates by the (n-2)-form $\om=f_1...f_r\,\wt\eta$.
As a consequence, the part of $Sing(\fa_{\wt\eta})$ contained in the pole divisor contains an union of curves:
given $J=(j_1,...,j_{n-1})\in\s_r^{n-1}$ let $S_J=\Pi(f_{j_1}=...=f_{j_{n-1}}=0)$. By the assumption on the pole divisor $f_1...f_r$, $S_J$ is a smooth complex curve and
\[
Sing(\fa_{\eta})\cap\,\Pi(f_1...f_j...f_r=0)\sup\bigcup_{J\in\s_r^{n-1}}S_J\,.
\]

A point $q=\Pi(p)\in S_J$, for a fixed $J\in\s_r^{n-1}$, will be said {\it generic} if for all $i\notin J$ then $f_i(p)\ne0$. Otherwise, $q$ will be said {\it non-generic}.
By the assumption on the pole divisor, if $q=\Pi(p)$ is non-generic and $f_i(p)=0$ then $f_\ell(p)\ne0$ for all $\ell\notin J\cup\{i\}$.

Let us fix $J=(j_1<...<j_{n-1})\in\s_r^{n-1}$ and a point $q=\Pi(q)\in S_J$. After an automorphism of $\p^n$ we can assume that $q=(0,...,0)$ in the affine chart $(x_0=1)\simeq\C^n$.
In this chart, the pole divisor of $\eta$ is $g_1...g_r$, where $g_j(x)=g_j(x_1,...,x_n)=f_j(1,x_1,...,x_n)$.
Since the equation of the curve $S_J$ is $(g_{j_1}=...=g_{j_{n-1}}=0)$, there exists a holomorphic coordinate system $(U,z=(z_1,...,z_n))$ around $q$ such that $g_{j_i}|_U=z_i$, $1\le i\le n-1$.

\begin{rem}\label{r:34}
{\rm Let $q\in S_J$ and $(U,z)$ be as above. We would like to observe that:
\begin{itemize}
\item[(a).] If $q$ is a generic point of $S_J$ then we can write
\begin{equation}\label{eq:10'}
\eta|_U=\sum_{j=1}^{n-1}\mu_j\,\frac{dz_1}{z_1}\wedge...\wedge\wh{\frac{dz_j}{z_j}}\wedge...\wedge\frac{dz_{n-1}}{z_{n-1}}+d\Te\,\,,
\end{equation}
where, either $\Te=0$, or $\Te$ is a non-closed logarithmic (n-3)-form with pole divisor contained in $x=z_1...z_{n-1}$, and $\mu_i=\la_{I_i}$, $I_i=J\setminus\{j_i\}$. 
\item[(b).] If $q\in S_J$ is a non-generic point then there exists $j\notin J$ such that $g_j(q)=0$ and $g_i(q)\ne0$ if $i\notin J\cup\{j\}$.
In this case, we can assume that $g_j|_U=z_n$. Moreover, we can write
\begin{equation}\label{eq:11}
\eta|_U=\sum_{1\le k<\ell\le n}\,\mu_{k\ell}\,\frac{dz_1}{z_1}\wedge...\wedge\wh{\frac{dz_k}{z_k}}\wedge...\wedge\wh{\frac{dz_\ell}{z_\ell}}\wedge...\wedge\frac{dz_n}{z_n}\,+\,d\Te\,\,,
\end{equation}
where $\Te$ is as in (a) and $\mu_{k\ell}=\la_{I_{k\ell}}$, $I_{k\ell}=J\cup\{j\}\setminus\{j_k,j_\ell\}$ if $\ell<n$, $\mu_{kn}=J\setminus\{j_k\}$.
\end{itemize}

The proof can be done directly by using (\ref{eq:la}) or theorem \ref{t:1}.}
\end{rem}

\subsection{Proof of the case $r=p+2$.}\label{ss:34}
In this case $r=p+2=n$ and the non generic points of $Sing(\fa_\eta)\cap\Pi(f_1...f_n=0)$ are in the finite set $\Pi(f_1=...=f_n=0)$. In particular, if we fix a non-generic point $q\in\Pi(f_1=...=f_n=0)$ there exists a local coordinate system $(U,z=(z_1,...,z_n))$ around $q$ such that $g_j|_U=z_j$, $1\le j\le n$.
In particular, by (\ref{eq:la}) we have
\[
\eta|_U=\sum_{1\le k<\ell\le n}\mu_{k\ell}\,\frac{dz_1}{z_1}\wedge...\wedge\wh{\frac{dz_k}{z_k}}\wedge...\wedge\wh{\frac{dz_\ell}{z_\ell}}\wedge...\wedge\frac{dz_n}{z_n}\,\,.
\]
Since $\eta\in\L^{n-2}_\fa(f_1,...,f_n)$ then $\eta|_U$ is locally decomposable outside the polar set $z_1...z_n=0$.
The foliation $\fa_\eta$ is defined in $U$ by the holomorphic form
\[
\tag{*}\om:=z_1...z_n\,\eta|_U=\sum_{1\le k<\ell\le n}\mu_{k\ell}\,z_k\,z_\ell\,dz_1\wedge...\wedge\wh{dz_k}\wedge...\wedge\wh{dz_\ell}\wedge...\wedge dz_n\,\,.
\]

\begin{rem}\label{r:35}
{\rm Let $\a$ be a holomorphic (n-2)-form on an open subset $V\sub\C^n$. Given $1\le j\le n$ and $p\in V$ such that $\a|_{(z_j=z_j(p))}\not\equiv0$ we can define a vector field $X^j_p$, tangent to the hyperplane $(z_j=z_j(p))$, by
\[
\a|_{(z_j=z_j(p))}=i_{X^j_p}\nu_j\,\,,\,\,\nu_j=dz_1\wedge...\wedge\wh{dz_j}\wedge...\wedge dz_n\,\,.
\]
This procces defines a holomorphic vector field $X^j$ on $V$, tangent to the fibration $(z_j=cte)$, by
\[
X^j(p)=X_p^j(p)\,\,,\,\,p\in V\,\,.
\]
Altough $i_{X^j_p}\,\a|_{(z_j=z_j(p))}=0$, in general $i_{X^j}\,\a\not\equiv0$. However, if the form $\a$ is locally decomposable outside its singular set then $i_{X^j}\,\a\equiv0$, so that $X^j$ is tangent to the distribution defined by $\a$. The proof is straightforward and is left to the reader.}
\end{rem}

If we apply remark \ref{r:35} to the (n-2)-form $\om$ in (*) we obtain $X^j=z_j.\,Y^j$, where
\[
Y^j=\sum_{k\ne j}\,\rho_k^j\,z_k\,\frac{\pa}{\pa z_k}\,\,,
\] 
and $\rho_k^j=(-1)^{k-1}\,\mu_{kj}$, with the convention $\mu_{rs}=-\mu_{sr}$, $\forall\,r,s$.
Since $\fa_\eta$ has dimension two, at least two of the linear vector fields above, that we can suppose to be $Y^1$ and $Y^2$, are not identically zero and generically linearly independent.
In this case, the form $\wt\om=i_{Y^1}i_{Y^2}\,\nu$, $\nu=dz_1\wedge...\wedge dz_n$, defines the same distribution as $\om$.
The reader can check that
\[
\wt\om=\sum_{1\le k<\ell\le n}(-1)^{k+\ell}\,(\rho^1_k\,\rho^2_\ell-\rho^1_\ell\,\rho^2_k)\,z_k\,z_\ell\,dz_1\wedge...\wedge\wh{dz_k}\wedge...\wedge\wh{dz_\ell}\wedge...\wedge dz_n\,\,.
\]
Since the coefficients of $\om$ and $\wt\om$ are homogeneous of degree two it follows that $\wt\om=c.\,\om$, where $c\in\C^*$.
From $\rho_k^j=(-1)^{k-1}\,\mu_{kj}$ and the above expression, we get
\[
\mu_{k1}\,\mu_{\ell2}-\mu_{\ell1}\,\mu_{k2}=c.\,\mu_{k\ell}\,\underset{k=2}{\implies}\,c=\mu_{12}\,\,.
\]

Now, consider the n-2 closed logarithmic 1-forms $\wt\te_3,...,\wt\te_n$ defined by
\[
\wt\te_j=\rho_j^2\,\frac{dz_1}{z_1}-\rho_j^1\,\frac{dz_2}{z_2}+\rho_2^1\frac{dz_j}{z_j}\,\,.
\]
Using that $\rho_1^2=-\rho_2^1$ we get $i_{Y^1}\,\wt\te_j=i_{Y^2}\,\wt\te_j=0$, $3\le j\le n$, and this implies that $\wt\te_3\wedge...\wedge\wt\te_n=k.\,\eta|_U$, $k\in\C^*$.
Comparing the coefficients of $\frac{dz_3}{z_3}\wedge...\wedge\frac{dz_n}{z_n}$ of the two members of the relation we get $k=(\rho_2^1)^{n-3}=\mu_{12}^{n-3}$.
Finally, if we define $\te_j=\rho_j^2\,\frac{df_1}{f_1}-\rho_j^1\,\frac{df_2}{f_2}+\rho_2^1\frac{df_j}{f_j}$ then $\te_3\wedge...\wedge\te_n$ then
\[
\te_3\wedge...\wedge\te_n=\mu_{12}^{n-3}\,\wt\eta\,\,,
\]
which proves that $\L^{n-2}_\fa(f_1,...,f_n)=\L^{n-2}_{td}(f_1,...,f_n)$.

\subsection{Proof of (b).}\label{ss:35}
Let us prove that if $r>n$ then $\L^{n-2}_{td}(f_1,...,f_r)$ is an irreducible component of $\L_\fa^{n-2}(f_1,...,f_r)$.
Since $\L_{td}^{n-2}(f_1,...,f_r)$ is irreducible and
$\L_{td}^{n-2}(f_1,...,f_r)\sub\L_\fa^{n-2}(f_1,...,f_r)$, it is clear that it is sufficient to prove that there exists $\fa_o\in\L_{td}(f_1,...,f_r)$ such that for any germ of curve $\tau\in(\C,0)\mapsto\fa_\tau\in\L_\fa^{n-2}(f_1,...,f_r)$, with $\fa_\tau|_{\tau=0}=\fa_o$, then $\fa_\tau\in\L_{td}^{n-2}(f_1,...,f_r)$ for all
$\tau\in(\C,0)$.

In the proof we will use the normal form of $\eta|_U$ near a generic point $p\in S_J$, $J=(j_1<...<j_{n-1})$ (see (\ref{eq:10'}) in remark \ref{r:34}).
The foliation $\fa_\eta$ is represented in $U$ by the (n-2)-form
\[
\om=z_1...z_{n-1}\eta|_U=\sum_{j=1}^{n-1}\mu_j\,z_j\,dz_1\wedge...\wedge\wh{dz_j}\wedge...\wedge dz_{n-1}+z_1...z_{n-1}\,d\Te\,\,.
\]
The linear part of $\om$ at the point $p=(0,...,0,c)$ is
\[
\om_1=\sum_{j=1}^{n-1}\mu_j\,z_j\,dz_1\wedge...\wedge\wh{dz_j}\wedge...\wedge dz_{n-1}=i_L\,dz_1\wedge...\wedge dz_{n-1}\,\,,
\]
where $L=\sum_{j=1}^{n-1}\rho_j\,z_j\,\frac{\pa}{\pa z_j}$ and $\rho_j=(-1)^{j-1}\,\mu_j$.
The numbers $\rho_1,...,\rho_{n-1}$ are the eigenvalues of $L$. If we fix a form
$\eta\in\L^{n-2}_\fa(f_1,...,f_r)$ then they depend only of $J$. 
\vskip.1in
The proof of (b) will be based in theorem \ref{t:31} that will be stated below. It will be used also in the proof of theorem \ref{t:5}.
In order to state it properly let us recall the definition of a generalized Kupka singularity for two dimensional foliations (see also \cite{ln}).

Let $\om$ be a germ at $p\in\C^n$ of integrable $(n-2)$-form with $\om(p)=0$. Recall that the rotational of $\om$ is the vector field $X=rot(\om)$ defined by
\[
d\om=i_X\,\nu\,,\,\nu=dz_1\wedge...\wedge dz_n\,.
\]
The singularity $p$ of $\om$ is of Kupka type if $X(p)\ne0$ and it is of {\it generalized Kupka type} (briefly g.K) if $X(p)=0$ and $p$ is an isolated singularity of $X$. When $X(p)=0$ and the linear part of $X$ at $p$ is non singular ($det(DX(p))\ne0$) we say that $p$ is non degenerated g.K (briefly n.d.g.K).
If $p$ is of Kupka or g.K type then the division theorem \cite{dr} implies that there exists another germ of holomorphic vector field, say $Y$, such that $\om=i_Yi_X\nu$.

\begin{rem}\label{r:36.a}
{\rm If $p$ is of Kupka type then there exists a local coordinate system $z=(z_1,...,z_n)$ around $p$ such that $z(p)=0$, $X=\frac{\pa}{\pa z_n}$ and $Y=\sum_{j=1}^{n-1}Y_j(z_1,...,z_{n-1})\frac{\pa}{\pa z_j}$, $Y(0)=0$. In particular, the foliation $\fa_\om$ has the structure of a local product, the germ of curve $\g=(z_1=...=z_{n-1}=0)$ is contained in the Kupka set of $\fa_\om$ and the vector field $Y$ defines the normal type of
$\fa_\om$ \cite{ln}.}
\end{rem}

In the next result we will consider the following situation: let $\fa$ be a two-dimensional foliation on $\p^n$, $n\ge4$. Assume that $Sing(\fa)$ contains a smooth irreducible curve, say $S$, with the following properties
\begin{itemize}
\item[(I).] There is a finite subset $F=\{p_1,...,p_k\}\sub S$ such that
$S\setminus F\sub K(\fa)$, the Kupka set of $\fa$. Since $S\setminus F$ is connected, the normal type of $\fa$ is the same at all points of $S\setminus F$.
We will denote by $Y$ a germ at $0\in\C^{n-1}$ of holomorphic vector representing this normal type.
\item[(II).] The eigenvalues of the linear part of $Y$, say $\rho_1,...,\rho_{n-1}$, are in the Poincaré domain and satisfy the following non-resonant conditions
\begin{itemize}
\item[($*$)] $\rho_j\ne\sum_{i\ne j}m_i\,\rho_i$ for all
$m=(m_1,...,\wh{m_j},...,m_{n-1})\in\Z_{\ge0}^{n-2}$ with $\sum_im_i\ge1$.
\end{itemize} 
In particular, we have $\rho_i\ne\rho_j$ if $i\ne j$.
Recall that $\rho_1,...,\rho_{n-1}$ are in the Poincaré domain if there exists $a\ne0$ such that $Re(a.\,\rho_j)>0$, $1\le j\le n-1$.
With these conditions the germ of vector field $Y$ is linearizable and semi-simple (cf. \cite{ar} and \cite{ma}).
 
\item[(III).] Given $p\in F$ let $\om$ be a germ of (n-2)-form defining the germ of $\fa$ at $p$. We will assume that there is a local coordinate system $(U,z=(z_1,...,z_n))$ around $p$ with the following properties:
\begin{itemize}
\item[(i).] $z(p)=0$ and $S\cap U=(z_1=...=z_{n-1}=0)$.
\item[(ii).] Set $X=rot(\om)$, so that $d\om=i_X\nu$, $\nu=dz_1\wedge...\wedge dz_n$. Let $\la_1,...,\la_n$ be the eigenvalues of the linear part $DX(p)$.
We will assume that there exists $a\ne 0$ such that $Re(a.\,\la_n)<0$ and $Re(a.\,\la_j)>0$, $\forall$ $1\le j\le n-1$.
Moreover, we will assume that the eingenspace of $DX(p)$ associated to the eigenvalue $\la_n$ is the tangent space $T_p\,S$.
\item[(iii).] Setting $\rho_n=0$, we will assume that $\la_i\,\rho_j-\la_j\,\rho_i\ne0$, $\forall$ $1\le i<j\le n$.
\end{itemize}
\end{itemize}

\begin{thm}\label{t:31}
If $\fa$ satisfies conditions {\rm (I), (II) and (III)} above then there exist homogeneous polynomials $g_1,...,g_r$ such that $\fa\in\L_{td}^{n-2}(g_1,...,g_r)$.
\end{thm}

Theorem \ref{t:31} will be used also in the proof of theorem \ref{t:5} and will be proved in \S \ref{ss:51}.
Let us finish the proof of (b) using theorem \ref{t:31}. In the case of the foliations $\fa_\eta\in\L^{n-2}_\fa(f_1,...,f_r)$ the curve $S$ will be
$S=\Pi(f_1=...=f_{n-1}=0)$ and the set $F$ will be
\[
F=\bigcup_{j\ge n}\left(S\cap\Pi(f_j=0)\right)\,.
\]
The idea is to find
$\eta_o\in\L^{n-2}_{td}(f_1,...,f_r)$ such that $\fa_{\eta_o}$ satisfies hypothesis (I) and (II) at the points of $S\setminus F$ and (III) at all points of $F$.
After that, we note that the set
\[
B:=\{\eta\in\L_\fa^{n-2}(f_1,...,f_r)\,|\,\fa_\eta\,\,\text{satisfies (I) and (II) on $S\setminus F$ and (III) at $F$}\}
\] 
is open in $\L^{n-2}_\fa(f_1,...,f_r)$. This is true, because (I) and (III) are open conditions and if $\rho=(\rho_1,...,\rho_{n-1})\in\C^{n-1}$ and $\rho_1,..,\rho_{n-1}$ are in the Poincaré domain then the number of possible resonances like in ($*$) invoving $\rho$ is finite. We leave the details for the reader.

In particular, there is an open neighborhood $V$ of $\eta_o$ in
$\L^{n-2}_\fa(f_1,...,f_r)$ such that $V\sub B$.
But, by theorem \ref{t:31} we must have $V\sub \L_{td}^{n-2}(f_1,...,f_r)$. Therefore $\L_{td}^{n-2}(f_1,...,f_r)$ is an irreducible component of
$\L^{n-2}_\fa(f_1,...,f_r)$.
\vskip.1in
In the next claim we will construct $\eta_o$ as above. It will also be used in the proof of theorem \ref{t:5}.

\begin{claim}\label{cl:31}
If $r\ge n$ then there are $\te_2,...,\te_{n-1}\in\L^1_\fa(f_1,...,f_r)$ such that, if
$\eta_o:=\te_2\wedge...\wedge\te_{n-1}$, then $\fa_{\eta_0}$ satisfies {\rm (I), (II) and (III)} of theorem \ref{t:31} along the curve $S=\Pi(f_1=...=f_{n-1}=0)$.
\end{claim}

{\it Proof.}
When $r=n-1$ and $S=\Pi(f_1=...=f_{n-1}=0)$ we have seen in the proof of corollary \ref{c:1} in \S \ref{ss:23} that the foliation $\fa_0$ is unique and defined by
$\te_0^2\wedge...\wedge\te_0^{n-1}$, where
\begin{equation}\label{eq:25}
\te_0^j=\frac{df_{j}}{f_{j}}-A_j\,\frac{df_1}{f_1}\,\,,\,A_j=\frac{d_{j}}{d_1}\,\,,\,\,2\le j\le n-1\,,
\end{equation}
with $d_j=deg(f_j)$.
In a neighborhood $U$ of any point $p\in S$ we can find local coordinates $z=(z_1,z_2,...,z_n)$ such that $f_j|_U=z_j$, $1\le j\le n-1$, $\implies$ $\te_0^j=\frac{dz_j}{z_j}-A_j\,\frac{dz_1}{z_1}$
and $\fa_{\eta_0}$ is defined by $\om=z_1...z_{n-1}\,\te_0^2\wedge...\wedge\te_0^{n-1}|_U$. Since the $\te_{j\,'s}^0$ are closed we get
\[
d\om=z_1...z_{n-1}\,\left(\frac{dz_1}{z_1}+...+\frac{dz_{n-1}}{z_{n-1}}\right)\wedge\te_0^2\wedge...\wedge\te_0^{n-1}|_U=\rho\,dz_1\wedge...\wedge dz_{n-1}\,,
\]
where $\rho=\frac{1}{d_1}\sum_jd_j\ne0$, as the reader can check. Therefore, all points in $S$ are of Kupka type and $S$ is a Kupka component of $\fa_{\eta_0}$. The normal type of
$\fa_{\eta_0}$ at $p$ can be defined by the linear vector field $Y=\sum_{j=1}^{n-1}d_j\,z_j\frac{\pa}{\pa z_j}$ because it satisfies $i_Y\te_0^j=0$, $\forall\,2\le j\le n-1$. In particular, the eigenvalues $\rho_j=d_j$, $1\le j\le n-1$, are in the Poincaré domain: $Re(\rho_j)>0$, $\forall\,1\le j\le n-1$.
\vskip.1in
When $r\ge n$ we will consider small deformations of the forms $\te_0^j$ above.
For instance, if $r=n$ then the non-generic points of $S$ are the points of the set $F=S\cap\Pi(f_n=0)$.

Let us consider the case $r=n$.
Given $\tau=(t_2,...,t_{n-1})\in\C^{n-2}$ consider the family of 1-forms
\[
\te_\tau^j=\frac{df_{j}}{f_{j}}-A_j(\tau)\,\frac{df_1}{f_1}-B_j(\tau)\,\frac{df_n}{f_n}\,,\,2\le j\le n-1\,,
\]
where $A_j(\tau)=\frac{d_j}{d_1}-t_j$ and $B_j(\tau)=\frac{t_j\,d_1}{d_n}$.
Note that $i_R\,\te_\tau^j=0$ $\forall\,2\le j\le n-1$, so that $\te_\tau^j\in\L^1_\fa(f_1,...,f_n)$,
$\forall\,\tau\in\C^{n-2}$, $\forall\,2\le j\le n-1$. Let $\fa_\tau$ be the foliation defined by
$\eta_\tau=\te_\tau^2\wedge...\wedge\te_\tau^{n-1}$.

If $p\in S$ is a generic point then $f_n(p)\ne0$ and there are local coordinates at $p$, $(U,z=(z_1,...,z_{n-1},z_n))$,
such that $f_j|_U=z_j$, $1\le j\le n-1$. In these coordinates we have $S\cap U=(z_1=...=z_{n-1}=0)$ and the normal type can be calculated by considering the restriction of
$\fa_\tau$ to a normal section, for instance $\Si:=(f_n=f_n(p))\cap U$. The foliation $\fa_\tau|_\Si$ is 
defined by the (n-2)-form
\[
z_1...z_{n-1}\,\eta_\tau|_\Si=z_1...z_{n-1}\,\left(\frac{dz_2}{z_2}-A_2(\tau)\frac{dz_1}{z_1}\right)\wedge...\wedge\left(\frac{dz_{n-1}}{z_{n-1}}-A_{n-1}(\tau)\frac{dz_1}{z_1}\right)\,.
\]
In particular, the normal type can be defined by the vector field
\[
Y_\tau=\sum_{j=1}^{n-1}\rho_j(\tau)\,z_j\frac{\pa}{\pa z_j}\,,
\]
where $\rho_1(\tau)=d_1$ and $\rho_j(\tau)=d_j-t_j\,d_1$, because $i_{Y_\tau}\te_\tau^i=0$, $\forall\,2\le i\le n-1$ 

If $|\tau|$ is small enougth then the genereric points of $S$ are of Kupka type and the eigenvalues of the normal type are in the Poincaré domain (these are open conditions). Moreover, the parameter $\tau$ can be chosen in such a way that the eigenvalues $\rho_1(\tau),...,\rho_{n-1}(\tau)$ satisfy the non-resonance conditions ($*$) of (II). This is a consequence of the fact that the set
$\{(\rho_2...,\rho_{n-1})\in \C^{n-2}\,|\,d_1=\rho_1,\rho_2,...,\rho_{n-1}$ satisfy conditions ($*$)$\}$ is dense in $\C^{n-2}$.

At a point $p\in F$ we can find local coordinates $(U,z=(z_1,...,z_n))$ such that
$z(p)=0$ and the foliation is defined by the form
$\om_\tau=z_1...z_n\,\wt{\te}_\tau^2\wedge...\wedge\wt{\te}_\tau^{n-1}$, where
\[
\wt\te_\tau^j=\frac{dz_{j}}{z_{j}}-A_j(\tau)\,\frac{dz_1}{z_1}-B_j(\tau)\,\frac{dz_n}{z_n}\,.
\]
Since the forms $\wt\te_\tau^j$ are closed, we get
\[
d\om_\tau=\left(\frac{dz_1}{z_1}+...+\frac{dz_n}{z_n}\right)\wedge\om_\tau\,.
\]

The rotational $X_\tau$ of $\om_\tau$ is defined by $d\om_\tau=i_{X_t}\,dz_1\wedge...\wedge dz_n$ and so $X_\tau=\sum_{j=1}^n\la_j(\tau)\,z_j\,\frac{\pa}{\pa z_j}$ is linear and 
must satisfy $i_{X_\tau}\left(\frac{dz_1}{z_1}+...+\frac{dz_n}{z_n}\right)=0$ and $i_{X_\tau}\wt\te_\tau^j=0$, $\forall\,2\le j\le n-1$.
It follows that the eigenvalues $\la_1(\tau),...,\la_n(\tau)$ must satisfy the homogeneous system
\begin{equation}\label{eq:26}
\left\{
\begin{matrix}
x_1+...+x_n=0\\
x_j-A_j(\tau)\,x_1-B_j(\tau)\,x_n=0\,,\,1\le j\le n-2\\
\end{matrix}
\right.
\end{equation}
When $\tau=0$ we are in the situation of the case $r=n-1$ and a solution of (\ref{eq:26}) is $x_j=d_j>0$, if $1\le j\le n-1$, and $x_n=-(d_1+...+d_{n-1})<0$. Therefore, $\la_1(0),...,\la_n(0)$ satisfy condition $Re(a.\la_n(0))<0$ and $Re(a.\la_j(0))>0$, $1\le j\le n-1$, for some $a\ne0$.
Of course, this implies that for small $|\tau|$ the eigenvalues of $X_\tau$ has eigenvalues that satisfy $Re(a_\tau.\la_n(\tau))<0$ and $Re(a_\tau.\la_j(\tau))>0$, $1\le j\le n-1$, for some $a_\tau\ne0$.
It remains to verify that $\fa_\tau$ satisfies condition (iii) of (III) near $p$.

First of all, recall that $X_\tau$ and $Z_\tau=Y_\tau=\sum_{j=1}^{n-1}\rho_j(\tau)\,z_j\frac{\pa}{\pa z_j}$ are tangent to $\fa_\tau$. Moreover, since $X_\tau\wedge Z_\tau\not\equiv0$ these vector fields generate the foliation in a neighborhood of $p=0$. In particular, we must have $\om_\tau=b.\,i_{X_\tau}i_{Z_\tau}\nu$ for some $b\ne0$.
If we set $\rho_n(\tau)=0$ then the coefficient of $\frac{\pa}{\pa z_i}\wedge\frac{\pa}{\pa z_j}$ in $X_\tau\wedge Z_\tau$ is  $\la_i(\tau)\,\rho_j(\tau)-\la_j(\tau)\,\rho_i(\tau)\ne0$ if $i<j$.
Therefore (iii) is equivalent to prove that all coefficients of $\om_\tau$ are not zero.

Set $\a=z_1...z_{n-1}\,\wt\te_0^2\wedge...\wedge\wt\te_0^{n-1}$. By the case $r=n-1$ we have $Sing(\a)=S\cap U=(z_1=...=z_{n-1}=0)$.
Since $\om_0=z_n.\,\a$ we have
\[
Sing(\om_0)=(z_n=0)\cup(z_1=...=z_{n-1}=0)\,.
\]
On the other hand, if $\tau\ne0$ then the form $\om_\tau$ can be written as
\[
\om_\tau=z_n.\,\a_\tau+dz_n\wedge\be_\tau
\]
where $\a_0=\a$, $\a_\tau$ has linear coefficients and
\[
\be_\tau=\sum_{1\le i<j\le n-1}A_{ij}(\tau)\,z_i\,z_j\,dz_1\wedge...\wedge\wh{dz_i}\wedge...\wedge\wh{dz_j}\wedge...\wedge dz_{n-1}\,,
\]
where $A_{ij}(\tau)=\pm(A_i(\tau)\,B_j(\tau)-A_j(\tau)\,B_i(\tau))$, if $i,j>1$ and $A_{1j}(\tau)=\pm\,B_j(\tau)$, if $2\le j\le n-1$.
We leave this computation for the reader. If $t_j\ne0$ then $A_{1j}(\tau)=\pm B_j(\tau)=\pm t_j.\,d_1/d_n\ne0$. If $i,j>1$ then
\[
A_{ij}(\tau)=\pm\,\frac{d_1}{d_n}\,(t_j\,d_i-t_i\,d_j)\,.
\]
Hence, we can choose $\tau$ small so that $t_j\,d_i-t_i\,d_j\ne0$, $\forall\,1<i<j\le n-1$.

In the case $r>n$ we consider the parameter $\tau=(t_{ji})^{2\le j\le n-1}_{n\le i\le r}$ and
\[
\te_{\tau}^{rj}=\frac{df_j}{f_j}-A_j(\tau)\frac{df_1}{f_1}-\sum_{i=n}^rB_{ji}(\tau)\frac{df_i}{f_i}\,,
\]
where $A_j(\tau)=\frac{d_j}{d_1}-\sum_{i=n}^rt_{ji}$ and $B_{ji}(\tau)=\frac{d_1}{d_i}\,t_{ji}$, $2\le j\le n-1$. It can be checked directly that $i_R\te_\tau^{rj}=0$, $\forall\,j$.
The proof of the claim in this case can be done by induction on $r\ge n$. We leave the details for the reader.
\qed

\vskip.1in

\section{Proof of theorems \ref{t:3} and \ref{t:4}.}\label{ss:4}

\subsection{Proof of theorem \ref{t:4}.}

Let $\eta_o=i_{Z_o}\,\nu$ be the germ of p-form on $(\C^{p+1},0)$ which can be extended to a germ of integrable p-form $\eta$ on $(\C^n=\C^{p+1}\times\C^{n-p-1},(0,0))$, as in the hypothesis of theorem \ref{t:4}: $cod(Sing(Z_o))\ge3$.

The points in $\C^n=\C^{p+1}\times\C^{n-p-1}$ will be denoted $(z,y)$, where $z\in\C^{p+1}$ and $y\in\C^{n-p-1}$.
We will consider representatives of $Z_o$, $\eta_o$ and $\eta$, denoted by the same letters, the first two defined in a neighborhood $V\sub\C^{p+1}$ of a closed polydisc $\ov{U}$ and the last defined in a neighborhood of $\ov{U}\times\{0\}$ in $\C^n$, so that
\[
\eta_o=\eta|_{V\times\{0\}}=i_{Z_o}\,dz_1\wedge...\wedge dz_{p+1}\,\,.
\]
We will assume $cod_{V}(Sing(Z_o))\ge3$.
Define a holomorphic vector field $Z$ in a neighborhood of $\ov{U}\times\{0\}$ by $Z(z,y_o)=Z_{y_o}(z)=\sum_{j=1}^{p+1}g_j(z,y_o)\,\frac{\pa}{\pa z_j}$, where
\[
\eta|_{(y=y_o)}=i_{Z_{y_o}}\,dz_1\wedge...\wedge dz_{p+1}\,\,.
\]
Since $\eta$ is l.t.d. outside its singular set, we have $i_Z\,\eta=0$.

Note that $Z(z,0)=Z_o$. Therefore, the hypothesis implies that there is a neighborhood $W$ of $\ov{U}\times\{0\}$ in $\C^n$ such that $cod(Sing(Z))\ge3$ and $W\cap(y=0)=V\times\{0\}$.

We assert that there are holomorphic vector fields $X_1,...,X_{n-p-1}$ defined in a smaller neighborhood of $\ov{U}\times\{0\}$, such that $i_{X_j}\,\eta=0$ and
\begin{equation}\label{eq:20}
X_j(x,y)=\frac{\pa}{\pa y_j}+\sum_{i=1}^{p+1}h_{ji}(z,y)\,\frac{\pa}{\pa z_i},\,\forall\,1\le j\le n-p-1\,.
\end{equation}

First of all, we note that the above assertion is true in a neighborhood of any point $(z_o,0)\in (V\times\{0\})\setminus Sing(Z_o)$.
This is true because for $(z,y)$ in a neighborhood $U_\a$ of $(z_o,0)$ some component of $Z(z,y)$ does not vanishes, say $g_{p+1}(z,y)\ne0$, so that
\[
\frac{(-1)^p}{g_{p+1}}\,\eta|_{U_\a}=dz_1\wedge...\wedge dz_p+\wt\Te\,,
\]
where $\wt\Te\wedge dz_{p+1}\wedge dy_1\wedge...\wedge dy_{n-p-1}\equiv0$. As the reader can check this implies the existence of holomorphic vector fields $X_{j\a}$ on $U_\a$ as in (\ref{eq:20}), $1\le j\le n-p-1$. It follows that there are:
\begin{itemize}
\item{} a polydisc $Q=Q_1\times Q_2\sub\C^{p+1}\times\C^{n-p-1}$, with $Q_1\sup\ov{U}$ and $0\in Q_2$.
\item{} a covering $\U=(U_\a)_{\a\in A}$ of $Q\setminus Sing(\eta)$ by polydiscs,
\item{} n-p-1 collections of holomorphic vector fields $(X_{j\a})_{\a\in A}$, $1\le j\le n-p-1$, $X_{j\a}\in \X(U_\a)$,
\end{itemize}
such that
\begin{itemize}
\item[(i).] $cod_Q(Sing(Z|_Q))\ge3$.
\item[(ii).] $X_{j\a}=\frac{\pa}{\pa y_j}+\sum_{i=1}^{p+1}g_{i\a}(z,y)\frac{\pa}{\pa z_i}$.
\item[(iii).] $i_{X_{j\a}}\eta=0$, $\forall 1\le j\le n-p-1$, $\forall\a\in A$.
\item[(iv).] for all $q=(z,y)\in U_\a$ then $ker(\eta(q))=\left<Z(q),X_{1\a}(q),...,X_{n-p-1\,\a}(q)\right>_\C$.
\end{itemize}
If $1\le j\le n-p-1$ and $U_\a\cap U_\be\ne\emp$ then
\[
X_{j\a}-X_{j\be}=\sum_{i=1}^{p+1}(g_{j\a}-g_{j\be})\frac{\pa}{\pa z_i}=h_{\a\be}^j.\,Z\,,
\]
where $h_{\a\be}^j\in \O(U_\a\cap U_\be)$. The collection $(h_{\a\be}^j)_{U_{\a\be}\ne\emp}$ is an additive cocycle. Since $cod(Sing(Z))\ge3$ by Cartan's theorem (cf. \cite{ct} and \cite{gr}) the cocycle is trivial; $h_{\a\be}^j=h_\a^j-h_\be^j$, $h_\a^j\in\O(U_\a)$. Hence, there exists a holomorphic vector field $X_j$ on $Q\setminus Sing(Z)$ as in (\ref{eq:20}) such that $i_{X_j}\eta=0$; $X_j|_{U_\a}=X_{j\a}-h_\a^j\,Z$. By Hartog's theorem $X_j$ can be extended to a holomorphic vector field on $Q$, denoted by the same letter.
In particular, we have
\begin{equation}\label{eq:21'}
ker(\eta)|_Q=\left<Z,X_1,...,X_{n-p-1}\right>_{\O(Q)}\,.
\end{equation}
Finally, (\ref{eq:21'}) and theorem 11 of \cite{ce} imply the theorem:
\begin{itemize}
\item[(I).] There exists a smaller polydisc $\ov{U}\times\{0\}\sub Q'\sub Q$ and holomorphic vector fields $Z',Y_1,...,Y_{n-p-1}\in\left<Z,X_1,...,X_{n-p-1}\right>_{\O(Q')}$ such that $[Y_i,Y_j]=0$, $[Z',Y_j]=0$, $\forall 1\le j\le n-p-1$, and
\[
\left<Z',Y_1,...,Y_{n-p-1}\right>_{\O(Q')}=\left<Z|_{Q'},X_1|_{Q'},...,X_{n-p-1}|_{Q'}\right>_{\O(Q')}\,.
\]
\item[(II).] There are coordinates $(z,w)=(z,w_1,...,w_{n-p-1})$ in $Q'$ such that $Y_j=\frac{\pa}{\pa w_j}$, $\forall\, 1\le j\le n-p-1$. 
\end{itemize}

This finishes the proof of theorem \ref{t:4}.
\qed

\vskip.1in
A simple consequence of theorem \ref{t:4} is the following:
\begin{cor}\label{c:41}
$Sing(\eta)$ is biholomorphic to $Sing(Z_o)\times(\C^{n-p-1},0)$.
\end{cor}

\subsection{Proof of theorem \ref{t:3}.}\label{ss:42}

In this section we consider a holomorphic codimension p foliation $\G$ on $\p^n$, $2\le p\le n-2$. We assume that there is a $p+1$ plane $\p^{p+1}=\Si_o\sub\p^n$ such that $cod_\Si(\G|_{\Si_o})\ge3$.
We want to prove that there is a linear projection $T\colon\p^n-\to\Si_o$ such that $\G=T^*(\G|_{\Si_o})$.
We will prove theorem \ref{t:3} in the case $n=p+2$, or equivalently, when the foliation is two-dimensional.
The general case will be reduced to this case using \S 3.4 of \cite{lw}.

The foliation $\G|_{\Si_o}$ is one dimensional and so it can be defined by a finite covering $(Q_\a)_{\a\in A}$ of $\Si_o$ by polydiscs of $\Si_o$, a collection $(X_\a)_{\a\in A}$ of holomorphic vector fields
$X_\a\in\X(Q_\a)$, and a multiplicative cocycle $(g_{\a\be})_{Q_\a\cap Q_\be\ne\emp}$ such that $X_\a=g_{\a\be}.\,X_\be$ on $Q_\a\cap Q_\be\ne\emp$.
A consequence of theorem \ref{t:4} is the following:

\begin{cor}\label{c:42}
There is a finite covering of $\Si_o$ by polydiscs of $\p^{p+2}$, say $(U_\a)_\a$, and two collections of holomorphic vector fields $(Z_\a)_{\a\in A}$ and $(Y_\a)_{\a\in A}$, $Z_\a,Y_\a\in\X(U_\a)$, with the following properties:
\begin{itemize}
\item[(a).] $U_\a\cap\Si_o=Q_\a$ and $Z_\a$ is an extension of $X_\a$ to $U_\a$. In particular, $Z_\a$ is tangent to $\Si_o$ along $Q_\a$.
\item[(b).] $Sing(Y_\a)=\emp$ and $Y_\a$ is transverse to $\Si_o$ along $Q_\a$.
\item[(c).] If $z\notin Sing(\G)\cap U_\a$ then $T_z\G=\left<Z_\a(z),Y_\a(z)\right>_\C$.
\item[(d).] If $z\in Sing(\G)\cap U_\a$ then $Z_\a(z)=0$. Moreover, the orbit of $Y_\a$ through $z$ is contained in $Sing(\G)$.
\end{itemize}
\end{cor}

The proof is a straightforward consequence of theorem \ref{t:4} and is left to the reader.
\vskip.1in

Our goal now is to prove the following:

\begin{lemma}\label{l:41}
Under the hypothesis of theorem \ref{t:3} assume that $n=p+2$. Then there is a one-dimensional foliation $\H$ of degree zero transverse to $\Si_o$ whoose leaves are $\G$-invariant.
\end{lemma}

{\it Proof.}
The foliation $\H$ will be constructed in homogeneous coordinates.
Let $\Pi\colon\C^{p+3}\setminus\{0\}\to\p^{p+2}$ be the canonical projection and $\wt\G=\Pi^*(\G)$.
Consider homogeneous coordinates $z=(z_0,...,z_{p+2})\in\C^{p+3}$ such that $\Pi^{-1}(\Si_o)\cup\{0\}=(z_o=0):=\wt\Si_o$.
In these homogeneous coordinates the foliation $\wt\G$ is defined by an integrable homogeneous p-form $\eta$ such that $i_R\eta=0$, where $R$ denotes the radial vector field on $\C^{p+3}$.
The foliation $\H$ will be defined in homogeneous coordinates by $R$ and a constant vector field $v$ such that $i_v\eta=0$.

The idea is to construct a formal series of vector fields of the form $V=\frac{\pa}{\pa z_0}+\sum_{j\ge0}z_0^j\,V_j$, where $V_j=\sum_{i=1}^{p+2}f_{ji}(z_1,...,z_{p+2})\,\frac{\pa}{\pa z_i}$,
the $f_{ji's}$ are holomorphic in some polydisc $Q$ of $\C^{p+2}$ containing the origin of $\C^{p+2}$ and such that $i_V\,\eta=0$.
If $v:=V(0)=\frac{\pa}{\pa z_0}+V_0(0)\ne0$ then $i_v\eta=0$ because the coefficients of $\eta$ are homogeneous of the same degree.
The constant vector field $v$ and $R$ induce a foliation $\H$ of degree zero on $\p^{p+2}$ tangent to $\G$.

Let us construct the series $V$.
The covering $(U_\a)_{\a\in A}$, given by corollary \ref{c:42}, induces the covering $\left(\wt{U}_\a=\Pi^{-1}(U_\a)\right)_{\a\in A}$ of $\wt\Si_o\setminus\{0\}$.
Without lost of generality, we can suppose that for any $\a\in A$ then $U_\a$ is contained in some affine chart $(z_{j(\a)}\ne0)$, where $j(\a)\ne0$.

\begin{claim}\label{cl:41}
{\rm There are collections of holomorphic vector fields $(\wt{Z}_\a)_{\a\in A}$ and $(\wt{Y}_\a)_{\a\in A}$, with $\wt{Z}_\a,\wt{Y}_\a\in\X(\wt{U}_\a)$ $\forall\a\in A$, with the following properties:
\begin{itemize}
\item[(i).] $D\Pi(z).\wt{Z}_\a(z)=Z_\a\circ\Pi(z)$ and $D\Pi(z).\wt{Y}_\a(z)=Y_\a\circ\Pi(z)$, $\forall z\in \wt{U}_\a$. In particular, $\wt{Y}_\a$ and $\wt{Z}_\a$ are tangent to $\wt\G|_{\wt{U}_\a}$: $i_{\wt{Y}_\a}\eta=0$ and $i_{\wt{Z}_\a}\eta=0$, $\forall\a$.
\item[(ii).] $\wt{Y}_\a$, $\wt{Z}_\a$ and $R$ generate $\wt\G$ in the sense that:
\begin{itemize}
\item{} if $z\in \wt{U}_\a\setminus Sing(\wt\G)$ then $T_z\wt\G=\left<\wt{Y}_\a(z),\wt{Z}_\a(z),R(z)\right>_\C$.
\item{} $z\in \wt{U}_\a\cap Sing(\wt\G)$ $\iff$ $\wt{Y}_\a(z)\wedge\wt{Z}_\a(z)\wedge R(z)=0$.
\end{itemize}
\item[(iii).] $\wt{Z_\a}$ is tangent to $\wt\Si_o$ along $\wt\Si_o\cap\wt{U}_\a$, $\forall\a\in A$. This means that
\[
\wt{Z}_\a(0,z_1,...,z_{p+2})\in\left<\frac{\pa}{\pa z_1},...,\frac{\pa}{\pa z_{p+2}}\right>_\O\,.
\]
\item[(iv).] $\wt{Y}_\a=g_\a(z)\,\frac{\pa}{\pa z_0}+V_\a$, where $V_\a\in\left<\frac{\pa}{\pa z_1},...,\frac{\pa}{\pa z_{p+2}}\right>_\O$ and $g_\a\in\O^*(\wt{U}_\a)$.
\end{itemize}
In particular, $Sing(\wt{Y}_\a)=\emp$ and $\wt{Y}_\a$ is transverse to $\wt\Si_o$ along $\wt\Si_o\cap\wt{U}_\a$, $\forall\a\in A$.}
\end{claim}

{\it Proof.}
Let us construct $\wt{Y}_\a$ and $\wt{Z}_\a$, $\a\in A$. Let $j\ne0$ be such that $U_\a\sub(z_j\ne0)$. Let us assume that $U_\a\sub(z_n=1)$, for instance, and that $Y_\a$ and $Z_\a$ are vector fields tangent to the affine plane $(z_n=1)$: $Y_\a=\sum_{i<n}g_i^\a(z_0,...,z_{n-1})\,\frac{\pa}{\pa z_i}$ and $Z_\a=\sum_{i<n}h_i^\a(z_0,...,z_{n-1})\,\frac{\pa}{\pa z_i}$, where $g_i^\a,h_i^\a\in\O(U_\a)$, $\forall\a$. Since $Y_\a$ is transverse to $\Si_o$ we have $g_0^\a\in\O^*(U_\a)$, $\forall\a$. The vector fields $\wt{Y}_\a$ and $\wt{Z}_\a$ are then constructed by extending $Y_\a$ and $Z_\a$ "radially":
we set $\wt{Y}_\a:=\sum_{i<n}\wt{g}_i^\a(z)\frac{\pa}{\pa z_i}$ and $\wt{Z}_\a:=\sum_{i<n}\wt{h}_i^\a(z)\frac{\pa}{\pa z_i}$, where $\wt{g}_i^\a(z)=z_0.\,g_i^\a(z_0/z_n,...,z_{n-1}/z_n)$ and
$\wt{h}_i^\a(z)=z_0.\,h_i^\a(z_0/z_n,...,z_{n-1}/z_n)$. We leave the proof of (i), (ii), (iii) and (iv) for the reader.
\qed

\vskip.1in

We now define a multiplicative cocycle of $3\times3$ matrices $(A_{\a\be})_{\wt{U}_\a\cap\wt{U}_\be\ne\emp}$. Since $cod(Sing(\G|_{U_\a}))\ge3$, we get $cod(Sing(\wt\G|_{\wt{U}_\a}))\ge3$, which implies
\[
cod\left(\{z\in\wt{U}_\a\,|\,\wt{Y}_\a(z)\wedge\wt{Z}_\a(z)\wedge R(z)=0\}\right)\ge3\,.
\]
From this and (ii) we get that, if $\wt{U}_\a\cap\wt{U}_\be\ne\emp$ then we can write
\[
\left\{
\begin{matrix}
\wt{Y}_\a(z)=a_{\a\be}(z)\,\wt{Y}_\be(z)+b_{\a\be}(z)\,\wt{Z}_\be(z)+c_{\a\be}(z)\,R(z)\\
\wt{Z}_\a(z)=d_{\a\be}(z)\,\wt{Y}_{\a\be}(z)+e_{\a\be}(z)\,\wt{Z}_\be(z)+f_{\a\be}(z)\,R(z)\\
\end{matrix}
\right.\,\,,\,\,\forall\,z\in\wt{U}_\a\cap\wt{U}_\be\,,
\]
where $a_{\a\be},...,f_{\a\be}\in\O(\wt{U}_\a\cap\wt{U}_\be)$.
The matrix
\[
A_{\a\be}:=
\left(
\begin{matrix}
a_{\a\be}&b_{\a\be}&c_{\a\be}\\
d_{\a\be}&e_{\a\be}&f_{\a\be}\\
0&0&1\\
\end{matrix}
\right)
\]
defines the transition
\begin{equation}\label{eq:mt}
\left(
\begin{matrix}
\wt{Y}_\a\\
\wt{Z}_\a\\
R\\
\end{matrix}
\right)=
A_{\a\be}.
\left(
\begin{matrix}
\wt{Y}_\be\\
\wt{Z}_\be\\
R\\
\end{matrix}
\right)\,\,.
\end{equation}
Of course, $A_{\a\be}=A_{\be\a}^{-1}$ and if $\wt{U}_\a\cap\wt{U}_\be\cap\wt{U}_\g\ne\emp$ then $A_{\a\be}\,A_{\be\g}\,A_{\g\a}=I$ on $\wt{U}_\a\cap\wt{U}_\be\cap\wt{U}_\g$.
Since $\wt{U}_\a\cap\wt{U}_\be$ is a neighborhood of $Q_\a\sub(z_0=0)$ in $\C^{p+3}$ we can write $A_{\a\be}$ as a power series in $z_0$:
\[
A_{\a\be}=\sum_{j\ge0}z_0^j\,A_{\a\be}^j\,\,,
\]
where $A_{\a\be}^j$ is a matrix with coefficients in $\O(Q_\a\cap Q_\be)$, $Q_\a=\wt{U}_\a\cap(z_0=0)$.
Now, the proof of Lemma \ref{l:41} can be reduced to the following:

\begin{lemma}\label{l:42}
The matrix cocycle $(A_{\a\be})_{\wt{U}_\a\cap\wt{U}_\be\ne\emp}$ is formally trivial: there exist a collection $(A_\a)_{\a\in A}$ of formal power series
\[
A_\a=\sum_{j\ge0}z_0^j\,A_\a^j\,\,,
\]
where
\begin{itemize}
\item[(a).] $A_\a^j$ is a matrix with coefficients in $\O(Q_\a)$, $Q_\a=\wt{U}_\a\cap(z_0=0)$, $\forall\a$, $\forall j\ge0$.
\item[(b).] $A_\a$ is invertible as a matrix formal series and its third line is $(0,0,1)$, $\forall\a$.
\item[(c).] if $Q_\a\cap Q_\be\ne\emp$ then $A_{\a\be}=A_\a^{-1}.\,A_\be$.
\item[(d).] $A_\a^0$ is triangular superior $\forall\a\in A$.
\end{itemize}
\end{lemma}
The proof of lemma \ref{l:42} will be done at the end of the section. Let us see how it implies lemma \ref{l:41}.
From (\ref{eq:mt}) we have
\[
\left(
\begin{matrix}
\wt{Y}_\a\\
\wt{Z}_\a\\
R\\
\end{matrix}
\right)=
A_\a^{-1}.\,A_\be.
\left(
\begin{matrix}
\wt{Y}_\be\\
\wt{Z}_\be\\
R\\
\end{matrix}
\right)
\,
\implies
\,
A_\a.\,
\left(
\begin{matrix}
\wt{Y}_\a\\
\wt{Z}_\a\\
R\\
\end{matrix}
\right)=
A_\be.\,
\left(
\begin{matrix}
\wt{Y}_\be\\
\wt{Z}_\be\\
R\\
\end{matrix}
\right)\,\,.
\]
Since the third line of $A_\a$ and $A_\be$ is $(0,0,1)$, it follows that there are formal series of vector fields $Y=\sum_{j\ge0}z_0^j\,Y_j$ and $Z=\sum_{j\ge0}z_0^j\,Z_j$ such that
\[
\left.\left(
\begin{matrix}
{Y}\\
{Z}\\
R\\
\end{matrix}
\right)\right|_{Q_\a\times(\C,0)}=
A_\a.\,
\left(
\begin{matrix}
\wt{Y}_\a\\
\wt{Z}_\a\\
R\\
\end{matrix}
\right)
\,\,,\,\,\forall\a\,.
\]
Note that $i_Y\eta=0$.
Since the coefficients of $\eta$ are homogeneous of the same degree, we obtain $i_v\eta=0$, where $v=Y(0)=Y_0(0)$.
Therefore, it is sufficient to see that $Y(0)=\sum_{j=0}^na_j\frac{\pa}{\pa z_j}$, where $a_0\ne0$.
This is a consequence of (d) in lemma \ref{l:42} and the fact that the $\frac{\pa}{\pa z_0}$ component of $\wt{Y}_\a$ does not vanishes at $(z_0=0)$, as the reader can check.
This finishes the proof of theorem \ref{t:3}.
\qed

\subsubsection{Proof of lemma \ref{l:42}.}
The restriction of $A_{\a\be}$ to $\wt\Si_o\cap\wt{U}_\a\cap\wt{U}_\be$ is triangular:
\begin{equation}\label{eq:22}
A_{\a\be}|_{\wt\Si_o\cap\wt{U}_\a\cap\wt{U}_\be}=
\left(
\begin{matrix}
a_{\a\be}&b_{\a\be}&c_{\a\be}\\
0&e_{\a\be}&f_{\a\be}\\
0&0&1\\
\end{matrix}
\right)\,\,.
\end{equation}
The cocycle defined by (\ref{eq:22}) is trivial, when restricted to a domain of $\wt\Si_o$ where we can apply Cartan's theorem \cite{ct}.
Fix two polydiscs $Q_1,\ov{Q}_2\sub\wt\Si_o$, where $Q_1=\{(z_1,...,z_{p+2})\,|\,|z_i|\le1\}$ and $\ov{Q}_2=\{(z_1,...,z_{p+2})\,|\,|z_i|\le1/2\}$, for instance.
The open set $H:=Q_1\setminus \ov{Q}_2$ is a Hartog's domain in $\wt\Si_o$, so that any $f\in\O(H)$ extends to a holomorphic function $\wt{f}\in\O(Q_1)$ (cf. \cite{si}).
By Cartan's theorem \cite{ct} we have $H^1(H,\O)=0$, because $n\ge3$. Since $H^2(H,\Z)=0$ we have also $H^1(H,\O^*)=0$.
Consider the Leray covering $(W_\a)_{\a\in A}$ of $H$ given by $W_\a=\wt{U}_\a\cap\wt\Si_o$. The restriction $A_{\a\be}|_{W_\a\cap W_\be}$ in (\ref{eq:22}) will be denoted by $B_{\a\be}$.
Since $B_{\a\be}$ is triangular, the entries $a_{\a\be}$ and $e_{\a\be}$ define multiplicative cocycles $(a_{\a\be})_{W_\a\cap W_\be\ne\emp}$ and  $(e_{\a\be})_{W_\a\cap W_\be\ne\emp}$, which are trivial: there are collections $(a_\a)_{\a\in A}$ and $(e_\a)_{\a\in A}$, $a_\a,e_\a\in\O^*(W_\a)$ such that $a_{\a\be}=a_\a^{-1}.\,a_\be$ and $e_{\a\be}=e_\a^{-1}.\,e_\be$ on $W_\a\cap W_\be\ne\emp$.
Hence, the cocycle $(B_{\a\be})_{W_\a\cap W_\be\ne\emp}$ is equivalent to a cocycle $(C_{\a\be})_{W_\a\cap W_\be\ne\emp}$, where
\[
C_{\a\be}=
\left(
\begin{matrix}
1&g_{\a\be}&h_{\a\be}\\
0&1&k_{\a\be}\\
0&0&1\\
\end{matrix}
\right)\,\,.
\]
By writing explicitly that $(C_{\a\be})_{W_\a\cap W_\be\ne\emp}$ is a multiplicative cocycle, we get that $(g_{\a\be})_{W_\a\cap W_\be\ne\emp}$ and $(k_{\a\be})_{W_\a\cap W_\be\ne\emp}$ are aditive cocycles.
In particular, there are collections $(g_\a)_\a$ and $(k_\a)_\a$ with $g_\a,k_\a\O(W_\a)$ such that $g_{\a\be}=g_\be-g_\a$ and $k_{\a\be}=k_\be-k_\a$ on $W_\a\cap W_\be\ne\emp$.
If we set
\[
M_\a=
\left(
\begin{matrix}
1&-g_\a&0\\
0&1&-k_\a\\
0&0&1\\
\end{matrix}
\right)
\]
then
\[
D_{\a\be}:=M_\a^{-1}\,C_{\a\be}\,M_\be=
\left(
\begin{matrix}
1&0&\ell_{\a\be}\\
0&1&0\\
0&0&1\\
\end{matrix}
\right)
\]
Using that $(D_{\a\be})_{V_\a\cap V_\be\ne\emp}$ is a multiplicative cocycle we obtain that $(\ell_{\a\be})_{W_\a\cap W_\be\ne\emp}$ is an aditive cocycle and $\ell_{\a\be}=\ell_\be-\ell_\a$ on
$W_\a\cap W_\be\ne\emp$. Finally, $L_\a^{-1}\,D_{\a\be}\,L_\be=I$, where
\[
L_\a=
\left(
\begin{matrix}
1&0&-\ell_\a\\
0&1&0\\
0&0&1\\
\end{matrix}
\right)
\]
as the reader can check.
From this, we obtain that there is a collection of invertible triangular superior matrices $(B_\a)_\a$ such that $B_{\a\be}=B_\a^{-1}B_\be$ on $W_\a\cap W_\be\ne\emp$.
Let
\[
B_\a=
\left(
\begin{matrix}
r_\a&s_a&t_\a\\
0&u_\a&v_\a\\
0&0&1\\
\end{matrix}
\right)\,\,.
\]
If $W_\a\cap W_\be\ne\emp$ then $r_\a\,\wt{Y}_\a+s_\a\,\wt{Z}_\a+t_\a\,R=r_\be\,\wt{Y}_\be+s_\be\,\wt{Z}_\be+t_\be\,R$ and $u_\a\,Z_\a+v_\a\,R=u_\be\,Z_\be+v_\be\,R$ on $W_\a\cap W_\be$. This defines two holomorphic vector fields $V_0$ and $Z_0$ on $H$ by
\[
V_0|_{W_\a}=r_\a\,\wt{Y}_\a+s_\a\,\wt{Z}_\a+t_\a\,R\,\,\text{and}\,\,Z_0|_{W_\a}=u_\a\,\wt{Z}_\a+v_\a\,R\,.
\]
Since $H$ is a Hartog's domain with holomorphic closure is the polydisc $Q_1$, $V_0$ and $Z_0$ can be extended to $Q_1$. We denote these extensions by the same symbols.
Moreover, we have $i_{V_0}\eta=i_{Z_0}\eta=0$.
We assert that $V_0(0)\ne0$. 

In fact, write $V_0(z)=\sum_{j=0}^{p+2}g_j(z)\,\frac{\pa}{\pa z_j}$, $z\in Q_1$. If $V_0(0)=0$ then $g_0(0)=0$ and so the analytic set $C:=\{z\in Q_1\,|\,g_0(z)=0\}$ must intersect the boundary $\pa Q_1$ of $Q_1$. If $z_0\in C\cap\pa Q_1$ then there is $\a\in A$ such that $z_0\in W_\a$. However, since $\wt{Z_\a}$ and $R$ are tangent to $\wt\Si_o$, we get $g_0(z_0)=r_\a(z_0).\,g_\a(z_0)\ne0$ (see (iv)), because $g_\a\in \O^*(\wt{U}_\a)$ and the matrix $B_\a$ is invertible.

Now, let us prove that there is a formal vector field $V=V_0+\sum_{j\ge1}z_0^j\,V_j$ such that $i_V\eta=0$. To do that we recall that $A_{\a\be}|_{W_\a\cap W_\be}=B_{\a\be}$ and $B_{\a\be}=B_\a^{-1}B_\be$.
Consider a collection of invertible matrices $(\wt{B}_\a)_{\a\in A}$, where $\wt{B}_\a$ is an extension of $B_\a$ to $\wt{U}_\a$.
Consider also the cocycle of matrices $\wt{A}_{\a\be}:=\wt{B}_\a.A_{\a\be}.\wt{B}_\be^{-1}$. This cocycle is equivalent to $A_{\a\be}$ and $\wt{A}_{\a\be}|_{W_\a\cap W_\be}=I$, $\forall\,W_\a\cap W_\be\ne\emp$.
Since $W_\a\cap W_\be=(z_0=0)\cap \wt{U}_\a\cap\wt{U}_\be$ we can write
\[
\wt{A}_{\a\be}=I+\sum_{j\ge1}z_0^j\,A_{\a\be}^j\,\,,
\]
where the entries of $A_{\a\be}^j$ are holomorphic in $W_\a\cap W_\be$. We claim that there are collections of power series of matrices of the form
\begin{equation}\label{eq:23}
A_\a=I+\sum_{j\ge1}z_0^j\,A_\a^j\,\,,\,\a\in A\,,
\end{equation}
such that the entries of $A_\a^j$ are holomorphic in $W_\a$ and $\wt{A}_{\a\be}=A_\a^{-1}\,A_\be$.
Given a power series in $z_0$, say $\s=\sum_{j\ge 0}z_0^j\,\s_j$, we will use the notation $J^k(\s)$ for the truncation $\sum_{0\le j\le k}z_0^j\,\s_j$.
The matrices $A_\a^j$ will be constructed by induction on $j\ge0$ in such a way that 
\[
\tag{$I_k$}\,\,\,J^k\left(\left(I+\sum_{1\le j\le k}z_0^j\,A_\a^j\right)^{-1}.\,\wt{A}_{\a\be}.\left(I+\sum_{1\le j\le k}z_0^j\,A_\be^j\right)\right)=I\,\,.
\]
Note that ($I_0$) is true and assume that we can construct collections $(A_\a^j)_{0\le j\le \ell-1}$, $\a\in A$, such that ($I_k$) is true for all $0\le k\le\ell-1$.
Set $\wt{A}_{\a\be}^\ell=I+\sum_{j=1}^\ell z_0^j\,A_{\a\be}^j$ and $C_\a^{\ell-1}=I+\sum_{j=1}^{\ell-1}z_0^j\,A_\a^j$. Since ($I_{\ell-1}$) is true, we get
\[
J^\ell\left[(C_\a^{\ell-1})^{-1}.\,\wt{A}_{\a\be}^\ell.\,C_\be^{\ell-1}\right]=I+z_0^\ell\,A_{\a\be}^\ell\,\,.
\]
Writing explicitly that the above expression is a multiplicative cocycle of matrices we get that $(A_{\a\be}^\ell)_{W_\a\cap W_\be\ne\emp}$ is an aditive cocycle.
Since $H^1(H,\O)=0$ we get collections $(A_\a^\ell)_{\a\in A}$ such that if $C_\a^\ell=I+\sum_{j=1}^\ell z_0^j\,A_\a^j$ then ($I_\ell$) is true.
In particular, the collection of formal series $C_\a=I+\sum_{j\ge1}z_0^j\,A_\a^j$, $\a\in A$, satisfies $C_\a^{-1}\,\wt{A}_{\a\be}\,C_\be=I$, so that
\[
\wt{A}_{\a\be}=C_\a.\,C_\be^{-1}\,\implies\,A_{\a\be}=\wt{B}_\a^{-1}.\,C_\a.\,C_\be^{-1}.\,\wt{B}_\be=(\wt{B}_\a^{-1}.\,C_\a).\,(\wt{B}_\be^{-1}.\,C_\be)^{-1}\,.
\]
This proves that the cocycle $(A_{\a\be})_{a\be}$ is formally trivial and finishes the proof of the existence of the constant vector field $v$ such that $i_v\eta=0$.

\section{Proofs of theorems \ref{t:5}, \ref{t:6} and \ref{t:31}.}\label{ss:5}

\subsection{Proof of theorem \ref{t:31}.}\label{ss:51}
Let $\fa$ be a two dimensional foliation on $\p^n$, $n\ge4$, having a curve $S$ in the singular set and that satisfies (I), (II) and (III) of section \ref{ss:35}.
The idea is to construct closed logarithmic 1-forms $\te_2,...,\te_{n-1}$, defined in a neighborhood $U$ of the curve $S$, such that $\te_2\wedge...\wedge\te_{n-1}$ defines the foliation $\fa|_U$.
By using an extension theorem of meromorphic functions (cf. \cite{ba} and \cite{ro}), each form $\te_j$ can be extended to a global closed meromorphic 1-form on $\p^n$, denoted again by $\te_j$, $2\le j\le n-1$. The fact that $\te_j|_V$ is logarithmic implies that $\te_j$ is also logarithmic: there are homogeneous polynomials in $\C^{n+1}$, say $g_1,...,g_r$, such that $\te_j\in\L^1_R(g_1,...,g_r)$, $\forall1\le j\le n-2$, $\implies$ $\fa\in\L_{td}^{n-2}(g_1,...,g_r)$.
The following result will be usefull:

\begin{thma}\label{t:32}
{\rm (Parametric linearization)} Let $(W_\tau)_{\tau\in(\C^k,0)}$ be a germ at $0\in\C^k$ of a holomorphic family of germs of holomorphic vector fields at $0\in\C^m$. Assume that:
\begin{itemize}
\item[(a).] The linear part $L_\tau=DW_\tau(0)$ is diagonal of the form
$L_\tau=\sum_{j=1}^m\rho_j(\tau)\,z_j\,\frac{\pa}{\pa z_j}$ in some local coordinate system $z=(z_1,...,z_m)$ around $0\in\C^m$.
\item[(b).] $\rho_1(0),...,\rho_m(0)$ are in the Poincaré domain and satisfy the non-resonance condition {\rm ($*$)} in {\rm (II)} of theorem \ref{t:31}.
\end{itemize}
Then there exists a holomorphic family of germs of biholomorphisms $(\Psi_\tau)_{\tau\in(\C^k,0)}$ such that $D\Psi_\tau(0)=I$ and
\[
\Psi_\tau^*(W_\tau)=L_\tau=\sum_{j=1}^m\rho_j(\tau)\,w_j\,\frac{\pa}{\pa w_j}\,\,.
\]
\end{thma}

Theorem \ref{t:32} is a parametric version of Poincaré's linearization theorem. Its proof can be found in \cite{ar} or \cite{ma}.
\vskip.1in
Let us continue the proof of theorem \ref{t:31}.
First of all, we will prove that there are n-2 closed logarithmic 1-forms $\te_2,...,\te_{n-1}$, defined in some neighborhood $W$ of $S\setminus F$, such that $\eta=\te_2\wedge...\wedge\te_{n-1}$ defines $\fa|_W$.

Fix $p\in S\setminus F$. Since $p\in K(\fa)$ there are local coordinates $(V,z=(z_1,...,z_n))$, with $p\in V$, such that
\begin{itemize}
\item[(i).] $z(p)=0$ and $S\cap V=(z_1=...=z_{n-1}=0)$.
\item[(ii).] $\fa|_V$ is defined by a holomorphic (n-2)-form $\om$ that can be written as $\om=i_Yi_X\nu$, where $X=\frac{\pa}{\pa z_n}$,
$Y=\sum_{j=1}^{n-1}Y_j(z_1,...,z_{n-1})\frac{\pa}{\pa z_j}$ is the normal type and $\nu=dz_n\wedge dz_1\wedge...\wedge dz_{n-1}$. 
\end{itemize}
Since the eigenvalues $\rho_1,...,\rho_{n-1}$ of $DY(0)$ satisfy the non-resonance conditions ($*$), by theorem \ref{t:32} (without parameters) we can assume that
$Y$ is linear
\[
Y=\sum_{j=1}^{n-1}\rho_j\,z_j\frac{\pa}{\pa z_j}\,,
\]
which implies
\[
\om=\sum_{j=1}^{n-1}(-1)^{j-1}\,\rho_j\,z_j\,dz_1\wedge...\wedge\wh{dz_j}\wedge...\wedge dz_{n-1}\,.
\]
In particular, the form $\eta_V:=\frac{1}{z_1...z_{n-1}}\,\om$ is logarithmic
\[
\eta_V=\sum_{j=1}^{n-1}(-1)^{j-1}\,\rho_j\,\frac{dz_1}{z_1}\wedge...\wedge\wh{\frac{dz_j}{z_j}}\wedge...\wedge\frac{dz_{n-1}}{z_{n-1}}\,.
\]
Note that $\eta_V$ can be decomposed as $\eta_V=\rho_1\,\te_V^2\wedge...\wedge\te_V^{n-1}$, where
\begin{equation}\label{eq:15.a}
\te_V^j=\frac{dz_j}{z_j}-\frac{\rho_j}{\rho_1}\,\frac{dz_1}{z_1}\,\,,
\end{equation}
because $i_Y\te_V^j=i_X\te_V^j=0$, $\forall\,2\le j\le n-1$.

The above argument implies that there exists a covering $\V$ of $S\setminus F$, by open sets, such that 
\begin{itemize}
\item[(iii).] For each $V\in \V$ there exists a coordinate system $z_V=(z_1,...,z_n)\colon V\to\C^n$ with $V=\{z\,|\,|z_j|<1\,,\,1\le j\le n\}$ and $S\cap V=(z_1=...=z_{n-1}=0)$.
\item[(iv).] If $\te_V^j$ is as (\ref{eq:15.a}), $2\le j\le n-1$, then the logarithmic form $\te_V^2\wedge...\wedge\te_V^{n-1}$ defines $\fa|_V$.
\item[(v).] The vector fields $X_V=\frac{\pa}{\pa z_n}$ and $Y_V=\sum_{j=1}^{n-1}\rho_j\,z_j\,\frac{\pa}{\pa z_j}$ generate $\fa|_V$. 
\end{itemize}
We assert that if $V,\wt{V}\in\V$ are such that $V\cap\wt{V}\ne\emp$ then $\te_V^j\equiv\te_{\wt{V}}^j$ on $V\cap\wt{V}$.

In fact, first of all let us remark that
\begin{itemize}
\item[(vi).] For all $j\in\{1,...,n-1\}$ the hypersurface $\Si_V^j:=(z_j=0)\sub V$ is invariant by $\fa|_V$. Moreover, if $\wh\Si^j\sub V$ is another smooth hypersurface which is $\fa|_V$ invariant and is tangent to $\Si_V^j$ along $S$ then $\wh\Si^j\sub \Si_V^j$.
\end{itemize}

Note that (vi) above is equivalent to the fact that the hyperplane $(z_j=0)$ is $Y$-invariant, $1\le j\le n-1$. Moreover, it is the unique smooth hypersurface which is $Y$-invariant and tangent to $(z_j=0)$.
This is well-known and is a consequence of the fact that $\rho_1,...,\rho_{n-1}$ satisfy ($*$) (see \cite{ar}).

Let $z_V=(z_1,...,z_n)$ and $z_{\wh{V}}=(\wh{z}_1,...,\wh{z}_n)$ be the coordinate systems of $V$ and $\wh{V}$ on which (iv), (v) and (vi) are true. We assert that $\wh{z}_j=u_j.\,z_j$ on $V\cap \wh{V}$, where $u_j(z)\ne0$ $\forall\,z\in V\cap \wh{V}$, $\forall$ $1\le j\le n-1$.

In fact, if we fix $1\le j\le n-1$, by (vi) we must have $\wh{z}_j=u.\,z_i$, where $u(z)\ne0$ $\forall\,z\in V\cap\wh{V}$, for some $1\le i\le n-1$. However, the fact that $\rho_\ell\ne\rho_j$ if $\ell\ne j$ implies that $i=j$, as the reader can check.
It follows that
\[
\te_{\wh{V}}^j=\frac{d\wh{z}_j}{\wh{z}_j}-\frac{\rho_j}{\rho_1}\,\frac{d\wh{z}_1}{\wh{z}_1}=
\frac{dz_j}{z_j}-\frac{\rho_j}{\rho_1}\,\frac{dz_1}{z_1}+dv=\te_V^j+dv\,\,,\,\,2\le j\le n-1\,,
\]
where $v=log(u_j)-\frac{\rho_j}{\rho_1}\,log(u_1)$ is holomorphic. Now, (iv) and (v) imply that
\[
i_{X_V}\,\te_{\wh{V}}^j=i_{X_V}\,\te_V^j=i_{Y_V}\,\te_{\wh{V}}^j=i_{Y_V}\,\te_V=0\,\,\implies\,\,X_V(v)=Y_V(v)=0\,\,.
\]
The first relation implies that $v(z)=v(z_1,...,z_{n-1})$, because $X_V=\frac{\pa}{\pa z_n}$. Since the eigenvalues of $Y_V$ are in the Poincaré domain $v(z_1,...,z_{n-1})$ is a constant and $dv=0$. Hence, $\te_{\wh{V}}^j=\te_V^j$ on $V\cap\wh{V}$, as asserted.
Therefore there are closed logarithmic 1-forms $\te_2,...,\te_{n-1}$, defined on $W=\bigcup_\V V$, such that $\fa|_W$ is defined by $\te_2\wedge...\wedge\te_{n-1}$, as asserted.
Let us prove that the forms $\te_j$ extend to a neighborhood of any point in $F$.

\vskip.1in

Given $p\in F$ let $\om$ be a germ of (n-2)-form defining the germ of $\fa$ at $p$. Let $(U,z=(x=z_1,...,z_n))$ be a coordinate system around $p$ as in (III), so that $z(p)=0$ and
$S\cap U=(z_1=...=z_{n-1}=0)$.
The rotational $X$ of $\om$ has eigenvalues $\la_1,...,\la_n$ and there exists $a\ne0$ such that $Re(a.\,\la_n)<0$ and $Re(a.\,\la_j)>0$, $\forall$ $1\le j\le n-1$. 
Since $p=0$ is an isolated singularity of $X$ there exists another germ of vector field $Z$ such that $\om=i_Z\,i_X\,\nu$, $\nu=dz_1\wedge...\wedge dz_n$.
The vector fields $X$ and $Z$ generate the germ of $\fa$ at $0$.

\begin{lemma}\label{l:32}
There are germs at $p$ of vector fields $\wt{X}$ and $\wt{Z}$ that generate the germ of $\fa$ at $p$ and a holomorphic coordinate system $(U_1,w=(w_1,...,w_n))$ around $p$, with the following properties:
\begin{itemize}
\item[(a)] $w(p)=0$ and $S\cap U_1=(w_1=...=w_{n-1}=0)$.
\item[(b)] $\wt{Z}(w)=\sum_{j=1}^{n-1}\rho_j\,w_j\,\frac{\pa}{\pa w_j}$. In particular, $\wt{Z}$ is the normal type of $\fa$ along $S\setminus F$.
\item[(c)] $\wt{X}=\sum_{j=1}^n\la_j\,w_j\,(1+\phi_j(w_n))\frac{\pa}{\pa w_j}$, where $\phi_j(0)=0$, $\forall\,1\le j\le n-1$.
\end{itemize}
\end{lemma}
{\it Proof.}
Let $W_u$ be the hyperplane of $T_p\p^n$ generated by the eigenspaces of $DX(p)$ associated to the eigenvalues $\la_1,...,\la_{n-1}$ and $W_s$ be the eigenspace associated to $\la_n$. Recall that we have assumed $W_s=T_p\,S$, which implies that $W_u$ is transverse to $S$ at $p$. 
The condition $Re(a.\,\la_n)<0$ and $Re(a.\,\la_j)>0$ implies that the vector field $a.\,X$ has an unique invariant smooth hypersurface $\Si_u$ tangent to $W_u$, which meets $S$ transversely at $p$.
This is a consequence of the existence of invariant manifolds for hyperbolic singularities of vector fields (see \cite{hps}). The hypersurface $\Si_u$ is the unstable manifold of the vector field $a.\,X$.
We assert that $\Si_u$ is also $Z$-invariant. For simplicity, we will assume $a=1$.

In the proof we will use the relation:
\begin{equation}\label{eq:16a}
[Z,X]=h\,X\,\,
\end{equation}
where $h\in\O_n$ and $h(0)=0$. Let us assume (\ref{eq:16a}) and prove that $\Si$ is $Z$-invariant.

Take representatives of $Z$, $X$ and $h$ defined in some small ball $B$ around $0$.
Let $Z_t$ and $X_t$ be the local flows of $Z$ and $X$, respectively.
Since $\Si_u$ is the unstable manifold of $X$ the real flow $X_t$ of $X$ satisfies $\underset{t\to-\infty}{\ell im}\,X_t(z)=0$ $\forall z\in \Si\cap B$.
Integrating (\ref{eq:16a}) we get
\[
Z_t^*(X)=\phi_t.\,X\,\,\text{, where}\,\,\phi_t(z)=exp\left(\int_0^th(Z_s(z))\,ds\right)\,\,.
\]
The above relation implies that $Z_t$ sends orbits of $X$ on orbits of $X$. 
Given $z\in \Si\cap B$ denote $O(z)=\{X_t(z)\,|\,t\le0\}$. Since $\underset{t\to-\infty}{\ell im}\,X_t(z)=0$ we get $\ov{O(z)}=O(z)\cup\{0\}$.
Let $\wt{O}(z)$ be the germ of $\ov{O(z)}$ at $0$. Note that $\wt{O}(z)\sub\Si$ and that $Z_t(\wt{O}(z))$ is a germ of curve through $0$ such that $Z_t(\wt{O}(z))\setminus\{0\}$ is an orbit of $X$.
This of course implies that $Z_t(\wt{O}(z))\sub \Si$. Hence, $\Si$ is $Z$-invariant.

\vskip.1in

{\it Proof of} (\ref{eq:16a}).
Since $d\om=i_X\,\nu$ and $\om=i_Z\,i_X\,\nu$ we have
\[
L_Z\om=i_Z(d\om)+d(i_Z\om)=\om\,\,\implies\,\,L_Z(d\om)=d\om\,\,\implies
\]
\[
i_X\,\nu=L_Z(i_X\,\nu)=i_{[Z,X]}\,\nu+i_X\,L_Z\nu=i_{[Z,X]}\,\nu+\nabla(Z)\,i_X\,\nu\,\,,
\]
where $\nabla(Z)=\sum_j\frac{\pa Z_j}{\pa z_j}$ is the divergence of $Z$. From this relation we get $[Z,X]=h\,X$, where $h=1-\nabla(Z)$.
We assert that $h(0)=0$.

In fact, let $X_1=DX(0)$ and $Z_1=DZ(0)$. Relation (\ref{eq:16a}) implies that
\[
[Z_1,X_1]=h(0).\,X_1\,\,.
\]
The above relation implies that if $h(0)\ne0$ then $X_1$ is nilpotent, so that $\la_1=...=\la_n=0$, a contradiction (see \cite{ln}).
In particular, we have proved that $X_1$ and $Z_1$ commute.

Let us continue the proof of lemma \ref{l:32}.
After a holomorphic change of variables, we can assume that $\Si_u\sub(z_n=0)$. Since $(z_n=0)$ is invariant for both vector fields, in the new coordinate system we can write the $n^{th}$ component of $X$ and $Z$ as $\la_n\,z_n\,(1+h_1(z))$ and $z_n\,f(z)$, respectively, where $h_1(0)=0$. If we set $\Psi:=-\frac{f(z)}{\la_n\,(1+h_1(z))}$ then the $n^{th}$ component of $\wt{Z}:=Z+\Psi\,X$ vanishes.
Moreover, $\om=i_{\wt{Z}}\,i_X\,\nu=i_Z\,i_X\,\nu$ and $[\wt{Z},X]=g\,X$, where $g=h-X(\Psi)$ and $g(0)=0$.
We assert that there are coordinates $(W,w=(w_1,...,w_{n-1},w_n))$ around $p\,$ such that
\begin{itemize}
\item[(i).] $w(p)=0$, $\Si_u\cap W=(w_n=0)$ and $S\cap W=(w_1=...=w_{n-1}=0)$.
\item[(ii).] $\wt{Z}=\phi(w_n).\,\sum_{j=1}^{n-1}\rho_j\,w_j\frac{\pa}{\pa w_j}$, where $\phi(0)\ne0$.
\end{itemize}

In fact, since the $n^{th}$ component of $\wt{Z}$ vanishes, the hyperplanes $\Si_c:=(z_n=c)$ are $\wt{Z}$-invariant. On the other hand, if $c\ne0$ then $\Si_c$ is transverse to $S=(z_1=...=z_{n-1}=0)$ and so $\wt{Z}|_{\Si_c}$ represents the normal type of $\fa$ in the section $\Si_c$. Therefore, the eigenvalues of $D\wt{Z}(0,c)$ are proportional to $\rho_1,...,\rho_{n-1}$. In other words there exists a $\phi\in\O_1$ such that the eigenvalues of $D\wt{Z}(0,c)|_{\Si_c}$ are $\phi(c).\,\rho_1,...,\phi(c).\,\rho_{n-1}$.
Considering $\wt{Z}$ as a 1-parameter family of germs of vector fields at $0\in\C^{n-1}$ and applying theorem \ref{t:32} to this family we get (ii) of the assertion. Now, we assert that there exists $\Phi\in\O_n$ such that if we set $\wt{X}:=e^\Phi.\,X$ then
\begin{equation}\label{eq:17a}
[\wt{Z},\wt{X}]=0\,\,.
\end{equation}

In fact, if $\Phi\in\O_n$ then
\[
[\wt{Z},\wt{X}]=[\wt{Z},e^\Phi.\,X]=e^\Phi.\,\wt{Z}(\Phi).\,X+e^\Phi.\,[\wt{Z},X]=e^\Phi(\wt{Z}(\Phi)+g)\,X\,.
\]
Therefore, we have to prove that $\wt{Z}(\Phi)=-g$ has a solution $\Phi\in\O_n$. Recall that $\wt{Z}=\phi(w_n).\,L$, where $L=\sum_{j=1}^{n-1}\rho_j\,w_j\frac{\pa}{\pa w_j}$.
Set $w=(x,w_n)$, $x=(w_1,...,w_{n-1})$. We can write
\[
-g(x,w_n)=\sum_\si b_\si(w_n).\,x^\si
\]
where $\si=(\si_1,...,\si_{n-1})\in\Z_{\ge0}^{n-1}$, $b_\si\in\O_1$ and $x^\si=w_1^{\si_1}...w_{n-1}^{\si_{n-1}}$.

Let $\si_0=(0,...0)$. We will prove below that $b_{\si_0}\equiv0$. Assuming this fact, the equation $\wt{Z}(\Phi)=-g$ has a formal solution $\Phi=\sum_\si c_\si(w_n)\,x^\si$ where
\[
c_\si(w_n)=\frac{b_\si(w_n)}{\phi(w_n)\,\left<\rho,\si\right>}\,\,,\,\,\left<\rho,\si\right>=\sum_{j=1}^{n-1}\rho_j\,\si_j\,\,.
\]
Since $\rho_1,...,\rho_{n-1}$ are in the Poincaré domain we have
\[
inf\{|\left<\rho,\si\right>|\,;\,\si\in\Z_{\ge0}^{n-1}\,,\,\si\ne(0,...,0)\}\ge C\,,
\]
where $C>0$. This implies that the formal series converges; $\Phi\in\O_n$.
\vskip.1in
{\it Proof that $b_{\si_0}(w_n)\equiv0$, or equivalently $g(0,w_n)\equiv0$.} First of all, the $n^{th}$ component of $[\wt{Z},X]$ is $\wt{Z}(X_n)$, where $X_n=\la_n\,w_n\,(1+h_1(x,w_n))$ is the $n^{th}$ component of $X$. Hence, $[\wt{Z},X]=g\,X$ implies that
\[
g.\,\la_n\,w_n\,(1+h_1)=\wt{Z}\left(\la_n\,w_n\,(1+h_1)\right)=\la_n\,w_n\,\wt{Z}(h_1)=\la_n\,w_n\,\phi(w_n)\,L(h_1)\,\implies
\]
\[
g(x,w_n)=\frac{\phi(w_n)\,\sum_{j=1}^{n-1}\rho_j\,w_j\,\frac{\pa h_1}{\pa w_j}}{1+h_1}\,\implies\,g(0,w_n)\equiv0\,\qed
\]

\vskip.1in

We can write $\wt{X}=\sum_{j=1}^nX^j(w)\frac{\pa}{\pa w_j}$, where $X^j(w)=\la_j\,w_j+h.o.t$, because $D\wt{X}(0)$ commutes with $L$. From $\wt{Z}=\phi(w_n).\,L$ we get that relation (\ref{eq:17a}) is equivalent to
\[
L\left(X^n\right)=0\,\,\text{and}\,\,\phi(w_n).\,L\left(X^j\right)=\wt{X}(\phi(w_n)\rho_j\,w_j)\,\,\text{if}\,\,1\le j\le n-1\,.
\]
From $L(X^n)=0$ we get $X^n(w)=\la_n\,w_n(1+\psi_n(w_n))$, $\phi_n(0)=0$, because the first integrals of $L$ are functions of $w_n$. In particular, if $1\le j\le n-1$ we get
\[
L\left(X^j\right)-\rho_j\,X^j=\frac{\phi'(w_n)}{\phi(w_n)}(1+\psi_n(w_n))\,\rho_n\,\rho_j\,w_j\,w_n\,,\,\,1\le j\le n-1\,\,.
\]
As the reader can check, the above relations imply that $\phi'(w_n)\equiv0$ and $L(X^j)=\rho_j\,X^j$, $1\le j\le n-1$.
Hence, $\phi$ is a non-zero constant, and we can suppose that $\phi=1$.
Finally, the solutions of $L(X^j)=\rho_j\,X^j$ with linear part $\la_j\,w_j$ are of the form $X^j(w)=\la_j\,w_j\,(1+\psi_j(w_n))$, $\psi_j(0)=0$. This finishes the proof of lemma \ref{l:32}.
\qed

\vskip.1in

Let us finish the proof that the forms $\te_j$, $2\le j\le n-1$, extend to a neighborhood of $p\in F$.
Define closed logarithmic 1-forms $\wt\te_j$, $2\le j\le n-1$, by
\[
\wt\te_j=\frac{dw_j}{w_j}-\frac{\rho_j}{\rho_1}\,\frac{dw_1}{w_1}\,-\,\zeta_j(w_n)\,\frac{dw_n}{w_n}\,\,,
\]
where
\[
\zeta_j(w_n)=\frac{\la_j\,(1+\phi_j(w_n))-\frac{\rho_j}{\rho_1}\,\la_1\,(1+\phi_1(w_n))}{\la_n\,(1+\phi_n(w_n))}\,\,.
\]
Note that $\zeta_j(0)\ne0$, because $\rho_1\,\la_j-\rho_j\,\la_1\ne0$. In particular, the pole divisor of $\wt\te_j$ contains $w_n$ with multiplicity one. 

The reader can check directly that $i_{\wt{Z}}\,\wt\te_j=i_{\wt{X}}\,\wt\te_j=0$, $\forall$ $2\le j\le n-1$, so that $\wt\eta:=\wt\te_2\wedge...\wedge\wt\te_{n-1}$ defines the germ of $\fa$ at $p$.
Taking representatives, we can assume that the $\wt\te_{j's}$ are defined in some polydisc $\wt{W}$ containing $p$ and with $F\cap\wt{W}=\{p\}$. We assert that $\wt\te_j=\te_j$ on $\wt{W}\cap W$, $2\le j\le n-1$.

In fact, fix a point $q\in S\cap\wt{W}\cap W$. We have seen that there are coordinates $(V,z=(z_1,...,z_{n-1},z_n))$ around $q$ such that $z(q)=0$, $S\cap V=(z_1=...=z_{n-1}=0)$, $\fa|_V$ is generated by the vector fields $X=\frac{\pa}{\pa z_n}$ and $Y=\sum_{j=1}^{n-1}\rho_j\,z_j\,\frac{\pa}{\pa z_j}$ and $\te_j|_V=\frac{dz_j}{z_j}-\frac{\rho_j}{\rho_1}\,\frac{dz_1}{z_1}$, $2\le j\le n-1$.
Note that $w_n|_V\in\O^*(V)$ and that $w_j|_V=v_j.\,z_j$, where $v_j\in\O^*(V)$, $1\le j\le n-1$.
This implies that $\wt\te_j|_V=\te_j|_V+df_j$, where $f_j$ is a primitive of the closed holomorphic form
\[
\frac{dv_j}{v_j}-\frac{\rho_j}{\rho_1}\,\frac{dv_1}{v_1}-\zeta_j(w_n)\,\frac{dw_n}{w_n}\,|_V\,\,.
\]
Finally, the fact that $i_X\te_j=i_Y\te_j=i_X\wt\te_j=i_Y\wt\te_j=0$ implies that $X(f_j)=Y(f_j)=0$ and so $f_j$ is a constant and $\wt\te_j|_V=\te_j|_V$, $2\le j\le n-1$.

\vskip.1in

We have proved that there are closed logarithmic 1-forms $\te_2,...,\te_{n-1}$ defined in a neighborhood $U$ of $S$ such that $\eta:=\te_2\wedge...\wedge\te_{n-1}$ defines $\fa|_U$.
By the extension theorem in \cite{ro} the form $\te_j$ can be extended to closed meromorphic 1-forms on $\p^n$, $2\le j\le n-1$, denoted by the same symbol.
The pole divisor of $\te_j$ must be reduced because the pole divisor of the restriction $\te_j|_U$ is reduced. Therefore $\te_j$ is logarithmic, $2\le j\le n-1$, and $\eta$ is t.d.l.f.
In particular, there exist $g_1,...,g_r$ such that $\eta\in\L^{n-1}_{td}(g_1,...,g_r)$.
This finishes the proof of theorem \ref{t:31}.
\qed

\subsection{Proof of theorem \ref{t:5} in the case of two dimensional foliations.}\label{ss:52}
We want to prove that for all $r\ge n$ and $d_1,...,d_r\ge1$ then $\ov{\L_{td}(d_1,...,d_r;2,n)}$ is an irreducible component of $\fol(D;2,n)$, where $D=\sum_jd_j-n+1$.
To avoid confusion we assume $d_1\le d_2\le ...\le d_r$.
Recall that the definition implies
\[
\L_{td}(d_1,...,d_r;2,n)=\bigcup_{\underset{1\le j\le r}{dg(f_j)=d_j}}\L_{td}^{n-2}(f_1,...,f_r)\,\,.
\]
Since $\ov{\L_{td}^{n-2}(f_1,...,f_r)}$ is irreducible for all polynomials $f_1,...,f_r$ with $deg(f_j)=d_j$, $1\le j\le r$, it is clear that
$\ov{\L_{td}(d_1,...,d_r;2,n)}$ is an irreducible algebraic subset of $\fol(D;2,n)$.
The proof that it is an irreducible component of $\fol(D;2,n)$ will be similar to the proof of (b) of theorem \ref{t:2} (see \S \ref{ss:35}).
The idea is to exhibit a foliation $\fa_0\in\L_{td}(d_1,...,d_r;2,n)$ such that for any germ of holomorphic deformation $t\in(\C,0)\mapsto\fa_t\in\fol(D;2.n)$, with $\fa_t|_{t=0}=\fa_0$, then $\fa_t\in\L_{td}(d_1,...,d_n;2,n)$ $\forall\,t\in(\C,0)$.

In order to do that, first of all let us fix homogeneous polynomials $f_1,...,f_r$ in $\C^{n+1}$ with $deg(f_j)=d_j$, $1\le j\le r$, such that the hypersurface $(f_1...f_r=0)\sub\C^{n+1}$ has a strictly ordinary singularity outside $0\in\C^{n+1}$. In particular, for any $J=(1\le j_1<...<j_{n-1}\le r)$ then the curve $S_J=\Pi(f_{j_1}=...=f_{j_{n-1}}=0)\sub\p^n$ is a smooth complete intersection.

From now on we fix $J=(1,2,...,n-1)$ and set $S_J=S$.
By claim \ref{cl:31} there exists $\eta_\tau=\te_\tau^2\wedge...\wedge\te^{n-1}_\tau\in\L_{td}^{n-2}(f_1,...,f_r)$ such that the foliation $\fa_{\eta_\tau}$ defined by $\eta_\tau$ satisfies (I), (II) and (III) of theorem \ref{t:31} along the curve $S$. The finite set of (I) is $F=S\cap\bigcup_{j\ge n}\Pi(f_j=0)$. 
\begin{rem}\label{r:51}
{\rm The parameter $\tau=(t_{ji})^{2\le j\le n-1}_{n\le i\le r}$ in claim \ref{cl:31} can be chosen in such a way that if $\fa_{\eta_\tau}\in\L_{td}(d_1',...,d_s';2,n)$, where $d_1'\le d_2'\le...\le d_r'$ then $d_i'=d_i$, $\forall i$. Recalling the definition of the $\te_\tau^{j's}$, an example in which $\fa_{\eta_\tau}$ belongs to two different $\L_{td's}$ is when $B_{jn}(\tau)=B_{jn+1}(\tau)$ for all $2\le j\le n-1$.
In this case, in the sum that defines $\te_\tau^j$ there are terms as below
\[
B_{jn}(\tau)\frac{df_n}{f_n}+B_{jn+1}(\tau)\frac{df_{n+1}}{f_{n+1}}=B_{jn}(\tau)\,\frac{d(f_n\,f_{n+1})}{f_n\,f_{n+1}}\,.
\]
In this case $\fa_{\eta_\tau}\in\L_{td}(d_1,...,d_r:2,D)\cap\L_{td}(d_1,...,d_{n-1},d_n\times d_{n+1},...,d_r;2,D)$.

For instance, if we choose the parameters $t_{ji}$ linearly independent over $\Z$ then the required property is true. From now on, we will assume this propety.}
\end{rem} 

Since $\fa_{\eta_o}$ satisfies property (III) of theorem \ref{t:31} along the curve $S$, all points of the finite set $F$ are n.d.g.K singularities of $\fa_{\eta_o}$. 
Fix any holomorphic germ of deformation $t\in(\C,0)\mapsto\fa_t\in\fol(D;2,n)$. 
The main fact that we will use is that the curve $S$ admits a $C^\infty$ deformation $t\in(\C,0)\mapsto S(t)$ such that $S(t)\sub Sing(\fa_t)$ and the foliation $\fa_t$ satisfies properties (I), (II), (III) of theorem \ref{t:31} along $S(t)$.

\begin{lemma}\label{c:51}
There exists a germ of $C^\infty$ isotopy $\Phi\colon(\C,0)\times S\mapsto\p^n$, such that, if we denote $S(t):=\phi(t,S)$, then:
\begin{itemize}
\item[(a).] $S(0)=S$ and $S(t)\sub Sing(\fa_t)$ is smooth $\forall s\in(\C,0)$. In particular, $S(t)$ is an algebraic complete intersection, $\forall t\in(\C,0)$.
\item[(b).] If $r=n-1$ then all poins of $S(t)$ are of Kupka type.
\item[(c).] If $r>n-1$ then any point $p\in F=S\cap\bigcup_{k\ge n}\Pi(f_k=0)$, $n-1<k\le r$, has a holomorphic deformation $t\in(\C,0)\mapsto P_p(t)$ such that $P_p(t)\in S(t)$ is a n.d.g.K singularity of $\fa_t$. Set $F(t):=\{P_p(t)\,|\,p\in F\}$.
\item[(d).] The points of $S(t)\setminus F(t)$ are in the Kupka set of $\fa(t)$. Moreover, if we denote by $Y_t$ the normal type of $\fa_t$ along $S(t)\setminus F(t)$ then the correspondence $t\in(\C,0)\mapsto Y_t$ is holomorphic.
\end{itemize}
\end{lemma}

{\it Proof.}
The argument for the proof of (c) uses the stability under deformations of the n.d.g.K points [theorem 3 of \cite{ln}]. 
The argument for the existence of the isotopy $\Phi$ is similar to [\cite{ln2}, lemma 2.3.3, p. 83] and uses essentially the
local stability under deformations of the Kupka set \cite{med} and of the n.d.g.K singular points \cite{ln}.
The fact that the deformed curve $S(t)$ satisfies (I), (II) and (III) for the foliation $\fa(t)$ is a consequence of the fact that these conditions are open.
We leave the details for the reader.
\qed

\vskip.1in
Let us finish the proof. We will assume that $\fa_{\eta_0}$ satisfies remark \ref{r:51}.
Lemma \ref{c:51} implies that the foliation $\fa(t)$ has a curve $S(t)$ in the singular set that satisfies (I), (II) and (III) of theorem \ref{t:31}. In particular, there are homogeneous polynomials $g_1(t),...,g_{s(t)}(t)$ such that $\fa_t\in\L_{td}^{n-2}(g_1(t),...,g_{s(t)}(t))$. Set $deg(g_j(t))=d_j(t)$. We assert that $s(t)=r$ and that we can assume $d_j(t)=d_j$, $1\le j\le r$.

In fact, since $D=\sum_{j=1}^{s(t)}d_j(t)-n+1$ we have $s(t)\le D+n-1$ and the number of possilities for the degrees $d_j(t)$, $1\le j\le s(t)$ is finite.
In particular, there is a germ of non-contable set $A\sub(\C,0)$ such that the functions $t\in A\mapsto s(t)$ and $t\in A\mapsto d_j(t)$, $1\le j\le s(t)$, are all constants, say
$s|_A=r'$ and $d_j|_A=d_j'$. In particular, $\fa(t)\in\L_{td}(d_1',...,d_{r'}';2;D)$ for all $t\in A$. Since $0$ is in the adherence of $A$ we get $\fa_{\eta_0}\in\L_{td}(d_1',...,d_{r'}';2,D)$.
Hence, $r'=r$ and $d_j'=d_j$, $1\le j\le r$, and $\fa(t)\in\L_{td}(d_1,...,d_r;2,D)$.
\qed

\subsection{Proof of theorem \ref{t:6}.}\label{ss:53}
Let $M\sub\p^n$ be a $m$-dimensional smooth algebraic submanifold, where $m<n$, and $\fa$ be a codimension p holomorphic foliation on $\p^n$, where $p+1\le m$. Assume that:
\begin{itemize}
\item[(I).] The set of tangencies of $\fa$ with $M$ has codimension $\ge2$ on $M$.
\item[(II).] $\fa|_M$ can be defined by a meromorphic closed p-form $\eta$.
\end{itemize}
We want to prove that $\eta$ admits a closed meromorphic extension $\wt\eta$ defining $\fa$ on $\p^n$.  In fact, this proof is similar to the proof of the extension theorem of \cite{CLS} (see also proposition 3.1.1 of \cite{ln2}).
The idea is to prove that $\eta$ admits a closed extension $\wh\eta$, defined in a neighborhood $U$ of $M$, such that $\fa|_U$ is represented by $\wh\eta$. After that, by \cite{ba} and \cite{ro}, the form $\wh\eta$ admits a meromorphic extension $\wt\eta$ to $\p^n$. Since $U$ is an open non-empty subset of $\p^n$, it is clear that $\wt\eta$ is closed and defines $\fa$ on $\p^n$.

Let $X=Sing(\fa|_M)$. Note that $X=Tang(\fa,M)\cup(Sing(\fa)\cap M)$, where $Tang(\fa,M)$ denotes the set of tangencies of $\fa$ and $M$. By (I) we have $cod_M(X)\ge2$. We begin by extending $\eta$ to a neighborhood of $M\setminus X$.

\vskip.1in

{\it 1. Extension to a neighborhood of $M\setminus X$.}
By definition, the foliation $\fa$ is transverse to $M$ at the points of $M\setminus X$. In particular, given $q\in M\setminus X$ there exists a local coordinate system around
$q$, $z=(z_1,...,z_n)\colon W\to\C^n$, with $z(W)$ a polydisc of $\C^n$, $z(q)=0\in\C^n$, and
such that
\begin{itemize}
\item[(i).] $M\cap W=(z_{m+1}=...=z_n=0)$.
\item[(ii).] The leaves of $\fa|_W$ are the levels $z_1=ct_1$,...,$z_p=ct_p$.
\end{itemize}

In particular, $\fa|_W$ is defined by the form $\Om_W=dz_1\wedge...\wedge dz_p$. Since $\fa|_{W\cap M}$ is also defined by $\eta|_{W\cap M}$ we must have $\eta|_{W\cap M}=f.\,\Om_W|_{W\cap M}$, where $f=f(z_1,...,z_m)$ is meromorphic on $W\cap M$. Since $\eta$ and $\Om_W$ are closed we get $df\wedge\Om_W=0$, which is equivalent to
\[
\frac{\pa f}{\pa z_j}=0\,,\,\forall\,p+1\le j\le m\,\,\implies
\]
$f(z_1,...,z_m)=f(z_1,...,z_p)$: $f$ depends only of $z_1,...,z_p$.
In particular, $\eta|_{W\cap M}$ admits an unique closed meromorphic extension to $W$ defining $\fa|_W$: $\wh\eta_W=f(z_1,...,z_p)\,dz_1\wedge...\wedge dz_p$.
This proves that $\eta|_{M\setminus X}$ admits an unique closed meromorphic extension $\wh\eta$ to a neighborhood $V$ of $M\setminus X$ representing $\fa|_V$.

\vskip.1in

{\it 2. Extension of $\wh\eta$ to a neighborhood of $M$.}
Since $cod_M(X)\ge2$, given $q\in X$ we can find a Hartog's domain $H\sub V$ such that $q\in\wh{H}$, the holomorphic closure of $H$ (for the details see \cite{ln2}).
Therefore, $\wh\eta$ admits a meromorphic extension to a neighborhood of $q$, by Levi's extension theorem \cite{si}. In particular, $\eta$ can be extended to a closed meromorphic p-form $\wt\eta$ defining $\fa$ on $\p^n$ by \cite{ba} and \cite{ro}.

Let us assume now that $\eta$ is logarithmic and let $(\wt\eta)_\infty=\wt{S}_1^{k_1}...\wt{S}_r^{k_r}$ be the decomposition of the pole divisor of $\wt\eta$ into irreducible components.
The pole divisor of $\eta$ will be then $(\eta)_\infty=(\wt\eta)_\infty\cap M$, which is reduced because $\eta$ is logarithmic. Hence, $k_1=...=k_r=1$ and $\wt\eta$ is logarithmic. 
\qed

\subsection{End of the proof of theorem \ref{t:5}.}\label{ss:54}
Recall that we want to prove that, if $k\ge3$, $n\ge5$ and $r\ge n-k+2:=p+2$ then $\ov{\L_{td}(d_1,...,d_r;k,n)}$ is an irreducible component of $\fol(D,k,n)$, where $D=\sum_jd_j-n+k-1$.
Fix $f_1,...,f_r$ homogeneous polynomials on $\C^{n+1}$ with the following properties:
\begin{itemize}
\item[(i).] $deg(f_j)=d_j$, $1\le j\le r$.
\item[(ii).] the hypersurface $(f_1...f_r=0)$ has s.o.s outside $0\in\C^{n+1}$.
\end{itemize}
Set $m=n-k+2$ and let $\p^m\simeq\Si\sub\p^n$ be a m-plane such that:
\begin{itemize}
\item[(iii).] If $\C^{m+1}\simeq E=\Pi^{-1}(\Si)\cup\{0\}\sub\C^{n+1}$ and $\wt{f}_j=f_j|_E$, $1\le j\le r$, then $(\wt{f}_1...\wt{f}_r=0)$ has s.o.s outside $0\in E$.
\end{itemize}
Such m-plane $E$ exists by transversality theory. In fact, it is sufficient to choose $E$ in such a way that for any sequence $I=(i_1<...<i_s)\in\s^r_s$, where $1\le s\le n-1$, then the algebraic smooth set
$\Pi(f_{i_1}=...=f_{i_s}=0)\sub\p^n$ meets transversely $\Si=\Pi(E)$ (see definition \ref{d:1} in \S \ref{ss:1}). We leave the details for the reader.

Since $m-2=n-k=p$, then for any $\fa\in\L^p_{td}(f_1,...,f_r)$ we have $\fa|_\Si\in\L^{m-2}_{td}(\wt{f}_1,...,\wt{f}_r)$, so that $\fa|_\Si$ is a two dimensional foliation.
Given a 1-form $\te=\sum_j\la_j\frac{df_j}{f_j}\in\L^1_\fa(f_1,...,f_r)$, we set $\wt\te=\sum_{j=1}^r\la_j\frac{d\wt{f}_j}{\wt{f}_j}$.

Choose $\fa_o\in\L^p_{td}(f_1,...,f_r)$ defined by a logarithmic form $\eta_o=\te_2\wedge...\wedge\te_{p+1}$, where $\wt\te_2,... \wt\te_{p+1}$ are as in claim \ref{cl:31} of \S \ref{ss:35}.
We assume also $\fa_o|_\Si$ satisfies remark \ref{r:51}: if $\fa_o|_\Si\in\L_{td}(d_1',...,d_s';2,m)$ then $s=r$ and $d_j'=d_j$, $1\le j\le r$.

Let $(\fa_t)_{t\in(\C,0)}$ be a germ of holomorphic 1-parameter family of foliations in $\fol(D,k,n)$ such that $\fa_t|_{(t=0)}=\fa_o$.
Consider the germ of 1-parameter family of two dimensional foliations $\wt\fa_t:=\fa_t|_\Si$, $t\in(\C,0)$. By the proof in \S\ref{ss:52} we get $\wt\fa_t\in\L_{td}(d_1,...,d_r;2,m)$, $\forall t\in(\C,0)$, so that $\wt\fa_t$ can be defined in homogeneous coordinates by a $m-2=n-k$ logarithmic form $\wt\eta_t\in\L^{m-2}(\wt{f}_{1\,t},...,\wt{f}_{r\,t})$, where $\wt{f}_{j\,t}|_{t=0}=\wt{f}_j$, $1\le j\le r$.
By theorem \ref{t:6} the foliation $\fa_t\in\fol(D,k,n)$ is logarithmic, $\forall t\in(\C,0)$, so that $\fa_t\in\L(d_1(t),...,d_{s_t}(t);k,n)$.
We assert that $s_t=r$ and $d_j(t)=d_j$, $1\le j\le r$.

In fact, since $\fa_t\in\L(d_1(t),...,d_{s_t}(t);k,n)$ we get $\wt\fa_t\in\L(d_1(t),...,d_{s_t}(t);2,m)$. Therefore, as in the proof of the two dimensional case, we have
$s_t=r$ and $d_j(t)=d_j$, $1\le j\le r$, $\forall t\in(\C,0)$. Finally, by corollary \ref{c:32} of \S\,\ref{ss:33} we get $\fa_t\in\L_{td}(d_1,...,d_r;k,n)$, $\forall t\in(\C,0)$.
This finishes the proof of theorem \ref{t:5}.
\qed

\section{Appendix. Proof of theorem \ref{l:23} (by Alcides Lins Neto)}\label{s:6}

Let $X$ be a germ at $0\in\C^n$ of an irreducible complete intersection with an isolated singularity at $0\in X$ and $dim_\C(X)=k$, where $2\le k\le n-1$. As before, set $X^*=X\setminus\{0\}$.
Let $m=n-k$ and $\I=\left<f_1,...,f_m\right>$ be the ideal defining: $X=(f_1=...=f_m=0)$. Since $0\in\C^n$ is an isolated singularity of $X$ then
\begin{equation}\label{eq:91}
df_1(z)\wedge...\wedge df_m(z)\ne0\,\,,\,\,\forall\,z\in X^*\,\,.
\end{equation}

If we fix representatives of $f_1,...,f_m$ in a polydisc $Q$ (denoted by the same letters), $0\in Q\sub\C^n$, we use the notation $X_0=Q$ and $X_s=\{z\in Q\,|\,f_1(z)=...=f_s(z)=0\}$, $1\le s\le m$.
We set also $X_s^*=X_s\setminus\{0\}$.

\begin{lemma}\label{l:61}
There are representatives of $f_1,...,f_m$ in a polydisc $Q$ such that:
\begin{itemize}
\item[(a).] $0$ is the unique singularity of $X_s$ in $Q$, $1\le s\le m$. In particular, $X_s^*$ is smooth of codimension $s$, $\forall\,1\le s\le m$.
\item[(b).] For all $0\le s\le m-1$ the function $f_{s+1}|_{X_s^*}$ is a submersion. In particular, $df_{s+1}(z)\ne0$ for all $z\in X_s^*$.
\end{itemize}
\end{lemma}

With lemma \ref{l:61} the proof of theorem \ref{l:23} is reduced to the following claim:
\vskip.1in
"Let $Q\sub\C^n$ be a polydisc with $0\in Q$. Let $X\sub Q$ be a connected complete intersection with a singularity $0\in X$, defined by $X=(f_1=...=f_{n-k}=0)$. Assume $2\le k\le n-1$ and:
\begin{itemize}
\item[(1).] $f_j$ has an isolated singularity at $0\in Q$, $\forall$ $1\le j\le n-k$.
\item[(2).] $\forall$ $I=(i_1,...,i_s)$, where $i_j\ne i_k$ if $j\ne k$, and $\forall$ $z\in (f_{i_1}=...=f_{i_s}=0)\setminus\{0\}$ then $df_{i_1}(z)\wedge...\wedge df_{i_s}(z)\ne0$. In particular:
\begin{itemize}
\item[(2.1).] $(f_{i_1}=...=f_{i_s}=0)\setminus\{0\}$ is smooth of codimension $s$. 
\item[(2.2).] $dim_{\,\C}(X)=k$.
\end{itemize} 
\end{itemize}

If $0\le\ell\le k-2$ then any $\ell$-form $\om_\ell\in\Om^\ell(X\setminus\{0\})$ admits an extension $\wt\om_\ell\in\Om^\ell(Q)$."

\vskip.1in

\begin{ex}
{\rm We would like to observe that the statement of theorem \ref{l:23} is not true for $k$ and $k-1$ forms. For instance, let $f\in\C[z_0,z_1,...,z_n]$, $n\ge3$, be a homogeneous polynomial of degree $\ge n+1$ and $X=(f=0)\sub\C^{n+1}$, so that $k-1=n-1$. Assume that $\Pi(X^*)\sub\p^n$ is smooth. It is known that there exists a non-vanishing holomorphic (n-1)-form on $X$, say $\a$. The (n-1)-form $\Pi^*(\a)$ is holomorphic on $X^*$ and has no holomorphic extension to any neighborhood of $0\in\C^{n+1}$.}
\end{ex}

In order to prove the above claim we will consider the situation below:

Let $Y$ be a connected complex manifold of dimension $n\ge3$ and $Z\sub Y$ be a codimension one complex codimension one submanifold defined by $f=0$, where $f\in\O(Y)$ and $0$ is a regular value of $f$. In particular, $Z$ is a smooth submanifold of $Y$.
For simplicity, we will use the notations $\Om^\ell$ for the sheaf of holomorphic $\ell$-forms on $Y$ and $Z$. Of course $\Om^0=\O$.

\begin{lemma}\label{lem:1}
In the above situation assume that $H^k(Y,\Om^\ell)=0$ for all $k$ and $\ell$ such that $k\ge1$ and $1\le k+\ell\le r+1$.
Then:
\begin{itemize}
\item[(a).] If $k\ge1$, $\ell\ge0$ and $r\ge1$ are such that $1\le k+\ell\le r$ then $H^k(Z,\Om^\ell)=0$.
\item[(b).] If $r\ge0$ and $0\le\ell\le r$ then any $\ell$-form $\om\in\Om^\ell(Z)$ can be extended to a $\ell$-form $\wt\om\in\Om^\ell(Y)$.
\end{itemize}
\end{lemma}

{\it Proof.}
We will use Leray's theorem (cf. \cite{gr} pg. 43). Let us consider Leray coverings $\U=(U_\a)_{\a\in A}$ and $\wt\U=(\wt{U}_{\wt\a})_{\wt\a\in\wt{A}}$ of $Z$ and $Y$ by open sets, respectively, such that: $A\sub\wt A$ and if $\a\in A$ then:
\begin{itemize}
\item[(i).] $\wt{U}_\a$ is the domain of a holomorphic chart $z^\a=(z_1,...,z_n)\colon \wt{U}_\a\to\C^n$, such that $\wt{U}_\a=\{z^\a\,|\,|z_j|<1,\,j=1,...,n\}$ and $f|_{\wt{U}_\a}=z_1$
\item[(ii).] $U_\a=\wt{U}_\a\cap Z$. In particular, $U_\a=\{z^\a\in \wt{U}_\a\,|\,z_1=0\}$.
\end{itemize}
Note that $\wt{U}_\a$ is biholomorphic to a polydisc of $\C^n$ and $U_\a$ to a polydisc of $\C^{n-1}$.
We assume also that:
\begin{itemize}
\item[(iii).] If $\a\in \wt A\setminus A$ then $\wt{U}_\a\cap Z=\emp$.
This implies that $A=\{\a\in\wt{A}\,|\,\wt{U}_\a\cap Z\ne\emp\}$.
\end{itemize}
Given $J=(j_0,...,j_k)\in \wt{A}^{k+1}$ (resp. $J\in A^{k+1}$) we set $\wt{U}_J=\wt{U}_{j_0}\cap...\cap \wt{U}_{j_k}$ (resp. $U_J=U_{j_0}\cap...\cap U_{j_k}$).
Note that by construction, if $J\in A^{k+1}$ is such that $U_J\ne\emp$ then $U_J=\wt{U}_J\cap Z$. Moreover, if $z^{\a_0}=(z_1,...,z_n)$ is a chart as in (i) then $U_J\sub \{z_1=0\}$.
\begin{claim}\label{claim:1}
{\rm Given $0\le \ell\le n-1$ and $J=(j_0,...,j_k)\in A^{k+1}$ such that $U_J\ne\emp$ then any $\ell$-form $\om$ on $U_J$ can be extended to a $\ell$-form on $\wt{U}_J$.
Moreover, if $\wt\om_1$ and $\wt\om_2$ are two extensions of $\om$ to $\wt{U}_J$ then:
\begin{itemize}
\item[(a).] $\wt\om_2-\wt\om_1=g.\,f$, $g\in\O(\wt{U}_J)$ if $\ell=0$.
\item[(b).] $\wt\om_2-\wt\om_1=\a\wedge df+f.\,\be$, where $\a\in \Om^{\ell-1}(\wt{U}_J)$ and $\be\in\Om^\ell(\wt{U}_J)$, if $\ell\ge1$.
\end{itemize}}
\end{claim}
{\it Proof.} 
Since $\wt{U}_\a$ is biholomorphic to a polydisc, for any $\a\in \wt{A}$ it follows that $\wt{U}_J$ is a Stein open subset of $Y$. Since $U_J=f^{-1}(0)\cap \wt{U}_J$ it follows that any holomorphic function $h\in\O(U_J)$ admits an extension $\wt{h}\in\O(\wt{U}_J)$ (cf. \cite{gr}). This proves the case $\ell=0$. When $\ell\ge1$, we consider the chart $z^{\a_0}=(z_1,...,z_n)\colon \wt{U}_{\a_0}\to\C^n$ where $f|_{U_{\a_0}}=z_1$ and $U_{\a_0}=\{z^\a\in\wt{U}_{\a_0}\,|\,z_1=0\}$, so that $U_J=\{z^{\a_0}\in \wt{U}_J\,|\,z_1=0\}$. In particular, any $\ell$-form $\om\in\Om^\ell(U_J)$ can be written as
\[
\om=\underset{I=(2\le i_1<...<i_\ell\le n)}\sum\,h_I.\,dz_{i_1}\wedge...\wedge dz_{i_\ell}\,\,\text{, where}\,\,h_I=h_I(z_2,...,z_n)\in\O(U_J)\,\,.
\]
By the case $\ell=0$ any function $h_{I}$ admits an extension $\wt{h}_I\in\O(\wt{U}_J)$. Therefore, $\om$ admits the extension
\[
\wt\om=\underset{I=(2\le i_1<...<i_\ell\le n)}\sum\,\wt{h}_I.\,dz_{i_1}\wedge...\wedge dz_{i_\ell}\in\Om^\ell(\wt{U}_J)\,\,.
\]

If $\wt\om_2$ and $\wt\om_1$ are two extensions of $\om$ to $\wt{U}_J$ then $(\wt\om_2-\wt\om_1)|_{z_1=0}=0$. Therefore, if $\ell=0$ then $\wt\om_2-\wt\om_1=g.\,z_1=g.\,f$ as in (a), whereas if $\ell\ge1$ then $\wt\om_2-\wt\om_1=\a\wedge dz_1+z_1.\,\be=\a\wedge df+f.\,\be$ as in (b).
\qed

\vskip.1in

Since $\Om^\ell$ is a holomorphic sheaf, by Leray's theorem we have $H^k(Z,\Om^\ell)=H^k(\U,\Om^\ell)$ and $H^k(Y,\Om^\ell)=H^k(\wt\U,\Om^\ell)$ for all $k\ge1$ and $\ell\ge0$.
Of course, $H^0(Z,\Om^\ell)=\Om^\ell(Z)$ and $H^0(Y,\Om^\ell)=\Om^\ell(Y)$.
Let us fix some notations (cf. \cite{gr}):
\begin{itemize}
\item[1.] $C^k(\U,\Om^\ell)$ (resp. $C^k(\wt\U,\Om^\ell)$) the $\O$-module of $k$-cochains of $\ell$-forms with respect to $\U$ (resp. with respect to $\wt\U$).
\item[2.] $\d=\d_k\colon C^k(*,\Om^\ell)\to C^{k+1}(*,\Om^\ell)$ the coboundary operator, where $*=\U$ or $\wt\U$. In this way, we have:
\[
H^k(*,\Om^\ell)=ker(\d_k)/Im(\d_{k-1})\,\,,\,\,k\ge0\,.
\]
\end{itemize}

Recall that
\[
C^k(\U,\Om^\ell)=\,\underset{J\in A^{k+1}}\prod\,\,\Om^\ell(U_J)\,\,,
\]
where $J=(\a_0,...,\a_k)\in A^{k+1}$ and $U_J=U_{\a_0}\cap...\cap U_{\a_k}$. In particular, a cochain in $\om_\ell^k\in C^k(\U,\Om^\ell)$ is of the form
\[
\om_\ell^k=\left(\om_J\right)_{J\in A^{k+1}}\,\,,\,\,\om_J\in \Om^\ell(U_J)\,\,.
\]
When $U_J=\emp$ by convenction we set $\om_J=0$.
Anagolously, a cochain $\wt\om_\ell^k\in C^k(\wt\U,\Om^\ell)$ is of the form
\[
\wt\om_\ell^k=\left(\wt\om_J\right)_{J\in \wt{A}^{k+1}}\,\,,\,\,\om_J\in \Om^\ell(\wt{U}_J)\,\,.
\]

\vskip.1in

{\it Restriction of cochains:} Given a cochain $\wt\om_\ell^k\in C^k(\wt\U,\Om^\ell)$, where $\wt\om_\ell^k=\left(\wt\om_J\right)_{J\in \wt{A}^{k+1}}$, its restriction to $Z$ is defined as
\[
\wt\om_\ell^k|_Z:=\left(\wt\om_J|_{U_J}\right)_{J\in A^{k+1}}\in C^k(\U,\Om^\ell)\,\,.
\]
Recall that if $J\in A^{k+1}$ then $U_J=\wt{U}_J\cap Z$.
\begin{rem}\label{rem:1}
{\rm Let $\wt\om_\ell^k,\,\wt\eta_\ell^k\in C^k(\wt\U,\Om^\ell)$ be two cochains with the same restriction to $Z$: $(\wt\eta_\ell^k-\wt\om_\ell^k)|_Z=0$. It follows from claim \ref{claim:1} that:
\begin{itemize}
\item[(a).] If $\ell=0$ then there exists a cochain $g_0^k=\left(g_J\right)_{J\in \wt{A}^{k+1}}\in C^k(\wt\U,\O)$ such that $\wt\eta_J-\wt\om_J=g_J.\,f$, for all $J\in \wt{A}^{k+1}$. In this case we will write
$\wt\eta_\ell^k-\wt\om_\ell^k=f.\,g_0^k$.
\item[(b).] If $\ell\ge1$ then there are cochains $\wt\a_{\ell-1}^k=\left(\wt\a_J\right)_{J\in \wt{A}^{k+1}}\in C^k(\wt\U,\Om^{\ell-1})$ and $\be_\ell^k=(\wt\be_J)_{J\in \wt{A}^{k+1}}\in C^k(\wt\U,\Om^\ell)$ such that $\wt\eta_J-\wt\om_J=\wt\a_J\wedge df+f.\,\wt\be_J$, for all $J\in \wt{A}^{k+1}$. In this case, we will write $\wt\eta_\ell^k-\wt\om_\ell^k=\wt\a_{\ell-1}^k\wedge df+f.\,\wt\be_\ell^k$.

We leave the details to the reader.
\end{itemize}}
\end{rem}

{\it Extension of cochains.} Claim \ref{claim:1} implies that given a cochain $\om_\ell^k\in C^k(\U,\Om^\ell)$ then there exists a cochain $\wt\om_\ell^k\in C^k(\wt\U,\Om^\ell)$ whoose restriction to $Z$ coincides with $\om_\ell^k$. We leave the details to the reader. The cochain $\wt\om_\ell^k$ will be called an {\it extension} of the cochain $\om_\ell^k$.

\vskip.1in

{\it Division of cochains.} Given a cochain $\be_\ell^k=(\be_J)_{J\in \wt{A}^{k+1}}\in C^k(\wt\U,\Om^\ell)$ we define the cochain
$\be_\ell^k\wedge df:=(\be_J\wedge df)_{J\in\wt{A}^{k+1}}\in C^k(\wt\U,\Om^{\ell+1})$. We would like to observe that, if $\ell\ge1$ and $df\wedge\be_\ell^k=0$ then there exists a cochain
$\be_{\ell-1}^k\in C^k(\wt\U,\Om^{\ell-1})$ such that $\be_\ell^k=\be_{\ell-1}^k\wedge df$. The proof is easy and is left to the reader. 

\vskip.1in

Let us assume the hypothesis of lemma \ref{lem:1}: $H^k(Y,\Om^\ell)=0$ if $k\ge1$ and $1\le k+\ell\le r+1$.

\begin{claim}\label{claim:2}
{\rm In the above situation, if $k\ge0$ and $\ell\ge0$ are such that $k+\ell\le r$ then any cocycle $\om_\ell^k\in C^k(\U,\Om^\ell)$ such that $\d\,\om_\ell^k=0$ admits an extension $\wt\om_\ell^k\in C^k(\wt\U,\Om^\ell)$ such that $\d\,\wt\om_\ell^k=0$.}
\end{claim}

{\it Proof.}
Let $\om^k_\ell\in C^k(\U,\Om^\ell)$ be such that $\d\,\om_\ell^k=0$. As we have seen before, $\om_\ell^k$ admits an extension $\wh\om_\ell^k\in C^k(\wt\U,\Om^\ell)$.
Then $\d\,\wh\om_\ell^k\in C^{k+1}(\wt\U,\Om^\ell)$ and so $\d\,\wh\om_\ell^k|_Z=\d\,\om_\ell^k=0$.

Let us assume first that $\ell=0$, so that $k+\ell=k\le r$. In this case, from remark \ref{rem:1} we obtain
$\d\,\wh\om_0^k=f.\,g_0^{k+1}$, where $g_0^{k+1}\in C^{k+1}(\wt\U,\O)$. Now, since $\d^2=0$, we have $f.\,\d\,g_0^{k+1}=0$, and so $\d\,g_0^{k+1}=0$. Since $k+1\le r+1$ the hypothesis implies that $H^{k+1}(\wt\U,\O)=0$ and so there exists a cochain $h_0^k\in C^k(\wt\U,\O)$ with $g_0^{k+1}=\d\,h_0^k$. Therefore,
\[
\d\,\wh\om_0^k=f.\,\d\,h_0^k\,\,\implies\,\,\d(\wh\om_0^k-f.\,h_0^k)=0\,\,.
\]
If we set $\wt\om_0^k=\wh\om_0^k-f.\,h_0^k$ then $\wt\om_0^k|_Z=\om_0^k$ and $\d\,\wt\om_0^k=0$, which proves in the case $\ell=0$.

Let us assume now that $\ell\ge1$. In this case, remark \ref{rem:1} implies that
\begin{equation}\label{equ:1}
\d\,\wh\om_\ell^k=\wh\a_{\ell-1}^{k+1}\wedge df+f.\,\wh\be_\ell^{k+1}\,\,,
\end{equation}
where $\wh\a_{\ell-1}^{k+1}\in C^{k+1}(\wt\U,\Om^{\ell-1})$ and $\wh\be_\ell^{k+1}\in C^{k+1}(\wt\U,\Om^\ell)$.
We assert that we can choose $\wh\a_{\ell-1}^{k+1}\in C^{k+1}(\wt\U,\Om^{\ell-1})$ and $\wh\be_\ell^{k+1}\in C^{k+1}(\wt\U,\Om^\ell)$ such that (\ref{equ:1}) is true and $\d\,\wh\be_\ell^{k+1}=0$.
Let us prove this assertion.

First we construct by induction a sequence of cochains
\[
\be_{\ell-j}^{k+j+1}\in C^{k+j+1}(\wt\U,\Om^{\ell-j})\,\,,\,\,j=0,...,\ell
\]
such that $\be_\ell^{k+1}=\wh\be_\ell^{k+1}$ and:
\begin{itemize}
\item[(I).] $\d\,\be_{\ell-j}^{k+j+1}\wedge df=0$, $\forall$ $j=0,...,\ell$.
\item[(II).] $\d\,\be_{\ell-j}^{k+j+1}=\be_{\ell-j-1}^{k+j+2}\wedge df$, $\forall$ $j=0,...,\ell-1$.
\end{itemize}
The construction is based in the division property. Since $\d^2=0$, relation (\ref{equ:1}) implies that
\[
\d\,\wh\a_{\ell-1}^{k+1}\wedge df+f.\,\d\,\wh\be_\ell^{k+1}=0\,\,\implies\,\,\d\,\wh\be_\ell^{k+1}\wedge df=0\,\,\implies\,\,\d\,\wh\be_\ell^{k+1}=\be_{\ell-1}^{k+2}\wedge df
\]
\[
\d\,\be_{\ell-1}^{k+2}\wedge df=0\,\,\implies\,\,\d\,\be_{\ell-1}^{k+2}=\be_{\ell-2}^{k+3}\wedge df\,\,\implies...\implies\,\,\d\,\be_{\ell-j}^{k+j+1}\wedge df=0
\]
\[
\d\,\be_{\ell-j}^{k+j+1}=\be_{\ell-j-1}^{k+j+2}\wedge df\,\,\implies...\implies\,\,\d\,\be_1^{k+\ell}\wedge df=0\,\,\implies\,\,\d\,\be_1^{k+\ell}=\be_0^{k+\ell+1}.\,df\,\,.
\]

Next, we will see that the sequence can be constructed in such a way that $\d\,\be_{j}^{k+\ell-j+1}=0$, $\forall$ $j=0,...,\ell$.
This involves another induction argument. 

{\it $1^{st}$ step: $j=0$.} From $\d\,\be_1^{k+\ell}=\be_0^{k+\ell+1}.\,df$ we get $\d\,\be_0^{k+\ell+1}=0$. Hence $\be_0^{k+\ell+1}\in ker(\d)$.

{\it $2^{nd}$ step.} Assume that we have constructed the sequence satisfying (I), (II) with $\d\,\be_{i}^{k+\ell-i+1}=0$ for $i=0,...,j-1\,\le\ell-1$ and let us prove that can assume that $\d\,\be_j^{k+\ell-j+1}=0$.

From (II) we have $\d\,\be_j^{k+\ell-j+1}=\be_{j-1}^{k+\ell-j+2}\wedge df$, where $\d\,\be_{j-1}^{k+\ell-j+2}=0$ by the induction hypothesis. Since $(k+\ell-j+2)+(j-1)=k+\ell+1\le r+1$ we have $H^{k+\ell-j+2}(\wt\U,\Om^{j-1})=0$ and so there exists a cochain $\g_{j-1}^{k+\ell-j+1}\in C^{k+\ell-j+1}(\wt\U,\Om^{j-1})$ such that $\be_{j-1}^{k+\ell-j+2}=\d\,\g_{j-1}^{k+\ell-j+1}$. Therefore, if we set $\wt\be_j^{k+\ell-j+1}=\be_j^{k+\ell-j+1}-\g_{j-1}^{k+\ell-j+1}\wedge df$ then
\[
\d\,\wt\be_j^{k+\ell-j+1}=\d\,\left(\be_j^{k+\ell-j+1}-\g_{j-1}^{k+\ell-j+1}\wedge df\right)=0\,\,.
\]
Moreover,
\[
\wt\be_j^{k+\ell-j+1}\wedge df=\be_j^{k+\ell-j+1}\wedge df=\be_{j+1}^{k+\ell-j}\,\,.
\]
Hence, if we replace $\be_j^{k+\ell-j+1}$ by $\wt\be_j^{k+\ell-j+1}$ in the sequence, then the new sequence still satisfies (I) and (II).

The induction process implies that there exists a cochain $\g_{\ell-1}^{k+1}\in C^{k+1}(\wt\U,\Om^{\ell-1})$ such that $\d\,(\wh\be_\ell^{k+1}-\g_{\ell-1}^{k+1}\wedge df)=0$. Hence, if we set $\wt\be_\ell^{\,k+1}=\wh\be_\ell^{k+1}-\g_{\ell-1}^{k+1}\wedge df$ and $\wt\a_{\ell-1}^{k+1}=\wh\a_{\ell-1}^{k+1}+f.\,h_{\ell-1}^{k+1}$ then (\ref{equ:1}) can be written as
\[
\d\,\wh\om_\ell^k=\wt\a_{\ell-1}^{k+1}\wedge df+f.\,\wt\be_\ell^{k+1}\,\,,\,\,\text{where}\,\,\d\,\wt\be_\ell^{k+1}=0\,\,.
\]
Since $H^{k+1}(\wt\U,\Om^\ell)=0$, there exists a cochain $\g_\ell^k\in C^k(\wt\U,\Om^\ell)$ such that $\wt\be_\ell^{k+1}=\d\,\g_\ell^k$. In particular, if we set $\ov\om_\ell^k=\wh\om_\ell^k-f.\,\g_\ell^k$
then $\ov\om_\ell^k|_Z=\wh\om_\ell^k|_Z=\om_\ell^k$ and
\begin{equation}\label{equ:2}
\d\,\ov\om_\ell^k=\wt\a_{\ell-1}^{k+1}\wedge df\,\,.
\end{equation}
If $\ell=1$ then $\wt\a_0^{k+1}\in \H^{k+1}(\wt\U,\O)$ and (\ref{equ:2}) implies that $\d\,\a_{\ell-1}^{k+1}=0$ and there exists a cochain $g_0^k\in C^k(\wt\U,\O)$ such that $\wt\a_0^{k+1}=\d\,g_0^k$.
In particular, the cochain $\wt\om_1^k=\ov\om_1^k-g_0^k.\,df$ satisfies $\d\,\wt\om_1^k=0$ and $\wt\om_1^k|_Z=\om_1^k$, proving claim \ref{claim:2} in this case.

Finally, when $\ell\ge2$ using (\ref{equ:2}) and an induction argument similar to that used in the case of $\wh\be_\ell^{k+1}$ it is possible to obtain a cochain $\g_{\ell-2}^{k+1}\in C^{k+1}(\wt\U,\Om^{\ell-2})$ such that $\d\,(\wt\a_{\ell-1}^{k+1}-\g_{\ell-2}^{k+1}\wedge df)=0$. Since $(\ell-1)+k+1=\ell+k\le r+1$ we have $H^{k+1}(\wt\U,\Om^{\ell-1})=0$, so that
$\wt\a_{\ell-1}^{k+1}-\g_{\ell-2}^{k+1}\wedge df=\d\,\eta_{\ell-1}^k$, where $\eta_{\ell-1}^k\in C^k(\wt\U,\Om^{\ell-1})$. 
From (\ref{equ:2}) we get
\[
\d\,\ov\om_\ell^{k+1}=\wt\a_{\ell-1}^{k+1}\wedge df=\d\,\eta_{\ell-1}^k\wedge df\,\,\implies\,\,\d\,(\ov\om_\ell^k-\eta_{\ell-1}^k\wedge df)=0\,\,.
\]
Hence, if we set $\wt\om_\ell^k=\ov\om_\ell^k-\eta_{\ell-1}^k\wedge df$ then $\d\,\wt\om_k^\ell=0$ and $\wt\om_\ell^k|_Z=\om_\ell^k$, which proves claim \ref{claim:2}.
\qed

Let us finish the proof of lemma \ref{lem:1}.

{\it Proof of (a).} By Leray's theorem it is suficient to prove that $H^k(\U,\Om^\ell)=0$, if $k\ge1$ and $k+\ell\le r$. If $\om_\ell^k\in C^k(\U,\Om^\ell)$ is such that $\d\,\om_\ell^k=0$ then by claim \ref{claim:2}, $\om_\ell^k$ admits an extension $\wt\om_\ell^k$ such that $\d\,\wt\om_\ell^k=0$. Since $k+\ell\le r<r+1$ then $H^k(\wt\U,\Om^\ell)=0$, so that $\wt\om_\ell^k=\d\,\wt\eta_\ell^{k-1}$ for some cochain $\wt\eta_\ell^{k-1}\in C^{k-1}(\wt\U,\Om^\ell)$. As the reader can check, this implies that $\om_\ell=\d\,\left(\wt\eta_\ell^{k-1}|_Z\right)$, which proves the assertion.

{\it Proof of (b).} Let $\om_\ell\in\Om^\ell(Z)$, where $\ell\le r$. We can associate to $\om_\ell$ a 0-cochain $\om_\ell^0=(\om_\ell|_{U_\a})_{\a\in A}$ with $\d\,\om_\ell^0=0$. By claim \ref{claim:2}, $\om_\ell^0$ admits an extension $\wt\om_\ell^0\in C^0(\wt\U,\Om^\ell)$ such that $\d\,\wt\om_\ell^0=0$. This is equivalent to say that there exists a section $\wt\om_\ell\in \Om^\ell(Y)$ such that
$\wt\om_\ell^0=(\wt\om_\ell|_{\wt{U}_\a})_{\a\in \wt{A}}$. Hence, $\wt\om_\ell$ extends $\om_\ell$ proving lemma \ref{lem:1}.
\qed

\vskip.1in

We are now in position to prove the statement of theorem \ref{l:23}. Let $0\in Q\sub\C^n$, $Q$ a polydisc, and $X=(f_1=...=f_{n-k}=0)$ be as in the statement of lemma 2.3.
Define a sequence of analytic complete intersections $X_0\sup X_1\sup...\sup X_{n-k}$, where $X_0=Q$ and $X_q=(f_1=...=f_q=0)$ if $1\le j\le n-k$, and set $X_q^*:=X_q\setminus\{0\}$, $0\le q\le n-k$.
The hypothesis implies the following:
\begin{itemize}
\item[(I).] $dim_{\,\C}(X_q)=k(q):=n-q$ and $X_q^*$ is smooth, $\forall$ $0\le q\le k$.
\item[(II).] $X_{q}=f_{q}^{-1}(0)\cap X_{q-1}$, $\forall$ $1\le q\le n-k$. Moreover, $0$ is a regular value of $f_q|_{X_{q-1}^*}$.
\end{itemize}
Recall that $k\ge 2$, so that $k(q)\ge 3$ if $q\le n-3$.
\begin{claim}\label{claim:3}
{\rm Let $p\ge1$, $\ell\ge0$ and $0\le q\le n-k-1$ be such that $1\le p+\ell\le k(q)-2$. Then $H^p(X_q^*,\Om^\ell)=0$.}
\end{claim}

{\it Proof.} The proof is by induction on $q=0,...,n-3$.
The case $q=0$ is consequence of a generalization of Cartan's theorem (cf. \cite{ct}): since $X_0=Q$ is Stein then (\cite{gr} pg. 133):
\[
H^p(X_0^*,\Om^\ell)=0\,\,,\,\,\forall\,p=1,...,n-2\,\,,\,\,\forall\,\,\ell\ge0\,\,.
\]
In particular, $H^p(X_0^*,\Om^\ell)=0$ if $p\ge1$, $\ell\ge0$ and $1\le p+\ell\le n-2=k(0)-2$.
The induction step is consequence of (a) of lemma \ref{lem:1}: let us assume that claim \ref{claim:3} is true for $q$, where $1\le q\le n-k-2$. Set $Y=X_q^*$, $Z=X_{q+1}^*$ and $f=f_{q+1}|_{X_{q+1}^*}$ in lemma \ref{lem:1}. The induction hypothesis implies that, if $p\ge1$ and $\ell\ge0$ are such that $1\le p+\ell\le k(q)-2$ then $H^p(X_q^*,\Om^\ell)=0$. In particular, (a) of lemma \ref{lem:1} implies that $H^p(X_{q+1}^*,\Om^\ell)=0$, $\forall$ $p\ge1$, $\ell\ge0$ such that $1\le p+\ell\le k(q)-3=k(q+1)-2$.
\qed

\vskip.1in

The extension property is consequence of (b) of lemma \ref{lem:1}. The idea is to use claim \ref{claim:3} and (b) of lemma 1 inductively. Let $\om_\ell\in\Om^\ell(X\setminus\{0\})$, where $\ell\le k-2$.
In the first step we set $Z=X_{n-k}^*=X\setminus\{0\}$, $Y=X_{n-k-1}^*$ and $f=f_{n-k}|_{X_{n-k-1}^*}$. From claim \ref{claim:3} we have $H^p(X_{n-k-1}^*,\Om^\ell)=0$ if $p\ge1$ and $\ell\ge0$ are such that $1\le p+\ell\le k(n-k-1)-2=k-1$. Hence, (b) of lemma \ref{lem:1} implies that if $0\le\ell\le k-2$ then any form $\om_\ell\in\Om(X_{n-k}^*)$ has an extension $\om_\ell^1\in\Om^\ell(X_{n-k-1}^*)$.
The induction step is similar: assume that $\om_\ell$ has an extension $\om_\ell^j\in\Om^\ell(X_{n-k-j}^*)$, where $1\le j\le n-k-1$. 
Since $\ell\le k-2<k-2+j=k(n-k-j)-2$, (b) of lemma \ref{lem:1} implies that $\om_\ell^j$ has an extension $\om_\ell^{j+1}\in\Om^\ell(X_{n-k-j-1}^*)$.
Finally, $\om_\ell$ has an extension $\om_\ell^{n-k}\in\Om^\ell(X_0^*)$, which by Hartog's theorem has an extension $\wt\om_\ell\in \Om^\ell(Q)$.
\qed

\bibliographystyle{amsalpha}

\end{document}